\documentclass[11pt,oneside]{amsart} 
\usepackage{amsmath}
\usepackage{amsthm}
\usepackage{amsfonts}
\usepackage{amssymb,amscd,epsf,verbatim}
\usepackage{mathrsfs}
\usepackage{graphicx}
\usepackage{latexsym}
\usepackage{lscape}
\usepackage[colorlinks=true]{hyperref}
\hypersetup{colorlinks, citecolor=blue, filecolor=black, linkcolor=red, urlcolor=green}
\usepackage{epstopdf}
\usepackage{tikz}
\usetikzlibrary{calc}
\usetikzlibrary{matrix,arrows,decorations.pathmorphing}
\usepackage{tikz-cd}
\usepackage{color}
\usepackage{multirow}  
\usepackage{scalefnt}
\usepackage{geometry}

\newcommand{\cl}{\overline}
\newcommand{\n}{\vartriangleleft}

\newcommand{\Z}{\mathbb{Z}}
\newcommand{\F}{\mathbb{F}}

\newcommand{\Q}{\mathbb{Q}}
\newcommand{\A}{\mathbb{A}}
\renewcommand{\P}{\mathbb{P}}
\newcommand{\C}{\mathbb{C}}
\newcommand{\T}{\mathbb{T}}
\newcommand{\mA}{\mathcal{A}}

\newcommand{\mC}{\mathcal{C}}

\newcommand{\mF}{\mathcal{F}}
\newcommand{\mG}{\mathcal{G}}
\newcommand{\mH}{\mathcal{H}}

\newcommand{\mJ}{\mathcal{J}}
\newcommand{\mM}{\mathcal{M}}

\newcommand{\mO}{\mathcal{O}}
\newcommand{\mP}{\mathcal{P}}

\newcommand{\mZ}{\mathcal{Z}}

\newcommand{\fa}{\mathfrak{a}}
\newcommand{\fb}{\mathfrak{b}}
\newcommand{\fc}{\mathfrak{c}}
\newcommand{\fd}{\mathfrak{d}}

\newcommand{\fp}{\mathfrak{p}}
\newcommand{\fq}{\mathfrak{q}}
\newcommand{\fm}{\mathfrak{m}}
\newcommand{\fX}{\mathfrak{X}}

\newcommand{\so}{\Rightarrow}
\renewcommand{\iff}{\Leftrightarrow}
\newcommand{\ra}{\rightarrow}

\newcommand{\sq}{\widetilde}
\newcommand{\minus}{\backslash}

\DeclareMathOperator{\im}{im}
\DeclareMathOperator{\Hom}{Hom}
\DeclareMathOperator{\End}{End}

\DeclareMathOperator{\SL}{SL}
\DeclareMathOperator{\GL}{GL}

\DeclareMathOperator{\spec}{Spec}

\DeclareMathOperator{\tr}{\textup{tr}}
\DeclareMathOperator{\FD}{D_{\cl \rho}^\square}
\DeclareMathOperator{\FDR}{R_{\cl \rho}^\square}
\DeclareMathOperator{\FDRx}{R_{\cl \rho}^{\square, \chi}}
\DeclareMathOperator{\FDRt}{R_{\cl \rho}^{\square, \tau}}
\DeclareMathOperator{\FDRtx}{R_{\cl \rho}^{\square, \tau, \chi}}
\DeclareMathOperator{\isom}{\;\xrightarrow{\: {}_{\sim} \:} \;}
\DeclareMathOperator{\ad}{\textup{ad}}
\DeclareMathOperator{\loc}{\textup{loc}}
\DeclareMathOperator{\Part}{\textup{Part}}

\DeclareMathOperator{\Igaln}{\textup{I}_n^{\textup{Gal}}}

\DeclareMathOperator{\Iwdn}{\textup{I}_n^{\textup{WD}}}

\newcommand{\bi}{\begin{itemize}}
\newcommand{\ei}{\end{itemize}}
\newcommand{\bt}{\begin{theorem}}
\newcommand{\et}{\end{theorem}}
\newcommand{\bbt}{\begin{theorem*}}
\newcommand{\eet}{\end{theorem*}}
\newcommand{\bp}{\begin{proposition}}
\newcommand{\ep}{\end{proposition}}
\newcommand{\bl}{\begin{lemma}}
\newcommand{\el}{\end{lemma}}
\newcommand{\bbl}{\begin{lemma*}}
\newcommand{\eel}{\end{lemma*}}
\newcommand{\bc}{\begin{corollary}}
\newcommand{\ec}{\end{corollary}}
\newcommand{\beg}{\begin{example}}
\newcommand{\eeg}{\end{example}}
\newcommand{\br}{\begin{remark}}
\newcommand{\er}{\end{remark}}
\newcommand{\bbr}{\begin{remark*}}
\newcommand{\eer}{\end{remark*}}
\newcommand{\bd}{\begin{definition}}
\newcommand{\ed}{\end{definition}}
\newcommand{\be}{\begin{enumerate}}
\newcommand{\ee}{\end{enumerate}}
\newcommand{\bex}{\begin{exercise}}
\newcommand{\eex}{\end{exercise}}
\newcommand{\bproof}{\begin{proof}}
\newcommand{\eproof}{\end{proof}}

\theoremstyle{theorem}
\newtheorem{theorem}{Theorem}[section]
\newtheorem{lemma}[theorem]{Lemma}
\newtheorem{proposition}[theorem]{Proposition}
\newtheorem{corollary}[theorem]{Corollary}
\newtheorem{deflem}[theorem]{Definition/Lemma}

\newtheorem{mainthm}{Theorem} 
\newtheorem{lemma*}[mainthm]{Lemma}
\newtheorem{proposition*}[mainthm]{Proposition}
\newtheorem{corollary*}[mainthm]{Corollary}

\theoremstyle{definition}
\newtheorem{definition}[theorem]{Definition}
\newtheorem{definition*}[mainthm]{Definition}
\newtheorem{conjecture}[theorem]{Conjecture}
\newtheorem{conjecture*}[mainthm]{Conjecture}
\newtheorem{remark}[theorem]{Remark}
\newtheorem{remark*}[mainthm]{Remark}

\theoremstyle{remark}
\newtheorem*{theorem*}{Theorem}
\newtheorem*{example}{Example}
\newtheorem*{notation}{Notation}

\title{The Breuil--M\'ezard conjecture for function fields}  
\author{Zijian Yao}
\address{Department of Mathematics, Harvard University.}
\email{zijian.yao.math@gmail.com}

\begin{document}

\maketitle

\begin{abstract}
Let $K$ be a local function field of characteristic $l$, $\F$ be a finite field over $\F_p$ where $l \ne p$, and $\cl \rho: G_K \ra \GL_n (\F)$ be a continuous representation. We apply the Taylor--Wiles--Kisin method over certain global function fields to construct a mod $p$ cycle map $\cl{\textup{cyc}}$, from mod $p$ representations of $\GL_n (\mO_K)$ to the mod $p$ fibers of the framed universal deformation ring  $\FDR$. This allows us to obtain a function field analog of the Breuil--M\'ezard conjecture. Then we use the technique of close fields to show that our result is compatible with the Breuil--M\'ezard conjecture for local number fields in the case of $l \ne p$, obtained by Shotton in \cite{S}. 
\end{abstract}

\tableofcontents

\section{Introduction}

\subsection{Introduction}

Let $K$ be a local field of residue characteristic $l$, and $L/\Q_p$ a finite extension with a uniformizer $\lambda$ and residue field $\F = \mO_L/\lambda$. Let $\cl \rho: G_K \ra \GL_n (\F)$ be a continuous mod $\lambda$ representation. The original Breuil--M\'ezard conjecture (cf. \cite{BM}, reformulated in \cite{K2, K3}) is stated for $n = 2$ and $K$ a finite extension of $\Q_l$ where $l = p$. It predicts the structure of $R^{\square, \tau, \mathbf{v}}_{\cl \rho}/\lambda$ in terms of the representation theory of $\GL_2 (\mO_K)$, where $R^{\square, \tau, \mathbf{v}}_{\cl \rho}$ is the potentially semistable deformation ring with Weil--Deligne inertia type $\tau$ and Hodge--Tate weights $\mathbf{v}$. 

More precisely, the conjecture takes the following form: for each $\tau, \mathbf{v}$ as above, one associates a smooth representation $\sigma(\tau) \otimes \sigma (\mathbf{v})$ of $\GL_2(\mO_K)$ via the theory of types. Let $\cl L_{\tau, \mathbf{v}} $ be the semi-simplified mod $\lambda$ reduction of an invariant lattice in $\sigma(\tau)\otimes \sigma( \mathbf{v})$. For each irreducible mod $\lambda$ representation $\cl \sigma$ of $\GL_2(\mO_K)$, let $ n_{\cl \sigma}(\tau, \mathbf{v})$ be the Jordan--Holder multiplicity of $\cl \sigma$ in $\cl L_{\tau, \mathbf{v}}$. The Breuil--M\'ezard conjecture then predicts that there exists non-negative integers $\mu_{\cl \sigma} (\cl \rho)$, such that for any $(\tau, \mathbf{v})$, the Hilbert--Samuel multiplicity $e(\textup{R}^{\square, \tau, \mathbf{v}}_{\cl \rho} /\lambda)$ is given by 
$$ e(\textup{R}^{\square, \tau, \mathbf{v}}_{\cl \rho} /\lambda) = \sum_{\cl \sigma} n_{\cl \sigma } (\tau, \mathbf{v}) \cdot \mu_{\cl \sigma} (\cl \rho).$$ 

The original conjecture is mostly proved by Kisin in \cite{K1} when $K = \Q_p$, using both the $p$-adic Langlands correspondence for $\GL_2 (\Q_p)$ and a global patching argument, and is used to deduce a modularity lifting theorem and the Fontaine--Mazur conjecture for $\GL_2$. In fact, the argument in \cite{K1} can be reversed to deduce the conjecture from modularity theorems. This is the strategy used by Gee--Kisin in \cite{GK} for more general $K$ and  Barsotti--Tate deformation rings. The conjecture 
has been generalized to all $n$ in a refined ``geometric'' form by Emerton--Gee in \cite{EG}, in terms of cycles. There the authors show that the general conjecture would again follow from modularity lifting theorems. 

Recently in \cite{S}, Shotton considers a simpler analog of the Breuil--M\'ezard conjecture in the case where $K/\Q_l$ and $l \ne p$, and proves this analog under mild hypotheses on $\cl \rho$ following the strategy of \cite{GK} and \cite{EG} -- namely using known modularity lifting type theorems. 

In this paper, we give a slight modification of the Breuil--M\'ezard conjecture for $\GL_n$ formulated  in \cite{S} and show that this conjecture holds in the function field setup. In the remaining of the introduction we briefly explain this formulation and describe our main results.  \\

The setup is as follows. On the Galois side let $\FDRx$ be the framed universal deformation ring of $\cl \rho$ with a fixed determinant $\chi: G_{K} \ra \mO_L^\times.$ For a Galois inertia type $\tau$, let $\textup{R}^{\square, \tau, \chi}_{\cl \rho}$ be the reduced quotient parametrizing $\cl L$-points of the Galois inertia type $\tau$. The deformation ring $\FDRx$ is reduced, flat over $\mO = \mO_L$ and $\spec \textup{R}^{\square, \tau, \chi}_{\cl \rho}$ is a union of irreducible components of $\spec \FDRx$. \footnote{Here we assume that $L$ is large enough in the sense of Remark \ref{remark:irreducible_component_inertia} (2).} Let $\mZ (\FDRx)$ be the free abelian group on the irreducible components of $\spec \FDRx$. Finally, let $\FDRx/\lambda$ be the mod $\lambda$ special fiber $\FDRx \otimes_{\mO_L} \F$ and similarly define $\mZ(\FDRx/\lambda)$. 

On the representation theory side, denote $\GL_n(\mO_K)$ by $\mathtt{K}$ and let $\textup{R}_{L}(\mathtt{K})$ be the Grothendieck group of the category of smooth representations of $\GL_n(\mO_K)$ valued in $L$. Let $B^{\textup{t}}(\GL_n)$ be the tempered Bernstein spectrum of $\GL_n (K)$, as defined in \cite{SZ}. From Subsection \ref{sec:Bernstein_decomposition}, the elements in $B^{\textup{t}} (\GL_ n)$ correspond bijectively to Galois inertia types. (For comparison, recall that each element in the usual Bernstein spectrum corresponds to a Weil--Deligne inertia type).

Roughly speaking, Shotton's formulation of the Breuil--M\'ezard conjecture can be stated in two parts, a representation theory part given by the ``generic multiplicity map'' and a Galois deformation part given by the ``component map''.

\begin{definition*}[Definition \ref{def:multiplicity_map}, Corollary \ref{cor:multiplicity_map}; Definition \ref{def:component_maps}] \label{def:intro_multiplicity} \indent 

\be
\item We define the generic multiplicity map 
$$m:  \textup{R}_{L} (\mathtt{K}) = \textup{R}_{L} (\GL_n(\mO_K))  \longrightarrow \underset{s \in B^{\textup{t}}(\GL_n) }{\oplus} \Z $$ 
by setting the $s$-component of $m(\sigma)$ to be the integer
$$ \dim \Hom_{\GL_n(\mO_K)} (\sigma, \pi_s|_{\GL_n(\mO_K)}), $$ 
where $\pi_s$ is any irreducible generic representation in the tempered Bernstein component of $s$ (see Remark \ref{remark:tempered_support_generic_rep}). This map is well-defined thanks to Lemma \ref{lemma:define_pi_tau} and Proposition \ref{lemma:image_m_direct_sum}. 
\item The component map $c = c_{\cl \rho}^{\chi}: \underset{s \in B^{\textup{t}} (\GL_n)}{\oplus} \Z  \longrightarrow \mZ(\FDRx)$ is given by sending $1_s$ to the irreducible components of $\spec \textup{R}^{\square, \tau_s, \chi}_{\cl \rho}$ in $\mZ (\FDRx)$. 
\ee
\end{definition*}

\noindent The composition $c \circ m$ is essentially the cycle map $\textup{cyc}$ defined in \cite{S} \footnote{with two minor differences: (i), we consider a fixed $\chi$, which seems to give finer information, and implies the version without specifying the determinant. See Remark \ref{remark:fix_determinant_imply_no_fix}.  (ii), for aesthetic reasons we have untwisted the cycle map in \cite{S} by an involution $\iota:  \textup{R}_{L} (\mathtt{K})  \ra  \textup{R}_{L} (\mathtt{K})$, which is given by taking dual representations. }. We can now state the Breuil--M\'ezard conjecture in the case where $l \ne p$ following \cite{S}. 

\begin{conjecture*}[Conjecture \ref{conj:BM}] \label{conj:intro}
Assume that $l \ne p$, then there exists a mod $\lambda$ cycle map 
$$ \cl{\textup{cyc}} : \textup{R}_{\F} (\mathtt{K}) = \textup{R}_{\F} (\GL_n(\mO_K))  \longrightarrow \mZ (\FDRx/\lambda) $$ making the following diagram commute: 
\[
\begin{tikzcd} [column sep=1.8em,row sep=1.8em]
  \textup{R}_{L} (\mathtt{K}) \arrow{rr}{m} \arrow[d, swap, "\textup{red}"] \arrow[drrr, swap, "\textup{cyc}"]  & & \underset{s \in B^{\textup{t}}(\GL_n)}{\bigoplus} \Z    \arrow{rd}{c_{\cl \rho}^{\chi}}     \\
  \textup{R}_{\F} (\mathtt{K}) \arrow[drrr, dashed, "\cl{\textup{cyc}}"] & & & \mZ(\textup{R}_{\cl \rho}^{\square, \chi}) \arrow[d, "\textup{red}"] \\ 
&  & & \mZ(\FDRx/\lambda)
\end{tikzcd}
\]
Note that by surjectivity of the reduction map on the left, such a $\cl{\textup{cyc}}$ is necessarily unique.  
\end{conjecture*}

From this slight change of perspective, the relation between conjecture \ref{conj:intro} and the classical form of the Breuil--M\'ezard conjecture (due to Breuil--M\'ezard \cite{BM} and Kisin \cite{K3}) becomes more transparent, with the theory of tempered types (and $\mathtt{K}$-types) providing a bridge between them. To be more precise, we recall the following (folklore) conjecture: there exists a section $\sigma_t$ of the generic multiplicity map $m$ which sends $1_s= (\cdots, 0,  1, 0, \cdots)$ to an irreducible representation of $\mathtt{K}$ for each $s \in B^t(\GL_n)$. Such a section is called a typical section in the article.  

The existence of $\sigma_t$ is closely related to the existence of types in the theory of admissible representations of compact subgroups of $\GL_n(K)$. Assuming $\sigma_t$ exists, our Breuil--M\'ezard conjecture then implies the classical formulation (cf. Subsection \ref{ss:BM_conj}, in particular \ref{sss:relation_with_classical_form}). Moreover, the existence of any section, not necessarily typical, would imply that Conjecture \ref{conj:intro} does not provide extra information about smooth representations of $\mathtt{K} = \GL_n (\mO_K)$ than the classical form. However, the cycle map $\cl{\textup{cyc}}$ provides more refined information than the Hilbert--Samuel multiplicities alone, so \textit{a priori} the classical form does not imply Conjecture \ref{conj:intro}. 

\begin{remark*}
Another perspective clarified by our formulation is the following: for inner forms $\GL_r(D)$ of $\GL_n(K)$ with maximal compact open subgroups $\mathtt{K}' = \GL_r(\mO_D)$, we can incorporate an inertial Jacquet-Langlands correspondence $\textup{JL}: \textup{Rep} (\mathtt{K}') \dashrightarrow \textup{Rep} (\mathtt{K})$ into the picture above. It can be shown that under the existence of $\textup{JL}$, which is compatible with the general Jacquet-Langlands correspondence and the theory of types, the Breuil--M\'ezard conjecture for $\GL_n$ essentially implies a certain form of Breuil--M\'ezard conjecture for $\GL_r(D)$. This will be the subject of a subsequent paper of the author, where we study the correspondence $\textup{JL}$ and the Breuil--M\'ezard conjectures for inner forms of $\GL_n$. 
\end{remark*}

\begin{remark*}
The diagram above also suggests that the Breuil--M\'ezard conjecture is closely related to a ``mod $p$ inertia Langlands correspondence''. 
In this way, the Breuil--M\'ezard conjectures (both in the case of $l = p$ and $l \ne p$) can be thought of as the compatibility between inertia Langlands correspondence and mod $p$ inertia Langlands correspondence. 
\end{remark*}
 
Now we turn to the main result of this article, which asserts that the Breuil--M\'ezard conjecture holds in the local function field case where $l \ne p$, under mild Taylor--Wiles type assumptions. 

\begin{mainthm}[Theorem \ref{thm:Main}] \label{thm:intro}
Let $K$ be a local function field over $\F_l$, and $\cl \rho: G_K \ra \GL_n (\F)$ a mod $p$ Galois representation with $l \ne p$. Assume that $p \nmid n (n-1)(l-1)$, then Conjecture \ref{conj:intro} holds. 
\end{mainthm}

Roughly speaking, to prove the theorem it suffices to construct an exact functor $M_\infty$, which sends a representation $\sigma$ of $\mathtt{K}$ with $\mO_L$-coefficients to a module $M_\infty(\sigma)$ over $\FDRx$, satisfying i). it is compatible under reduction mod $\lambda$, and ii). its support is (up to a constant) given by $\textup{cyc} (\sigma)$. The precise statement is given in Proposition \ref{prop:main_prop}. The main ingredient of the construction of the desired $M_\infty$ is the Taylor--Wiles--Kisin patching method (following the strategy of \cite{K3} and \cite{S}). To be slightly more precise, we first realize $\cl \rho$ as a local component of the mod $\lambda$ reduction of a global representation $r: G_{F, S} \ra \GL_n(\mO_{\cl L})$ for some global field $F$ over $\F_l$. The representation $r$ is carefully chosen to satisfy a list of properties, the most important one being that it comes from an automorphic form $\pi$ of a central division algebra $D$ over $F$, where $D_v^\times \cong \GL_n(K)$ for some place $v$. Now we consider a suitable global Galois deformation ring $R_{\cl r, S}^{\square_T, \eta}$ with local constraints, equipped with a surjective map $R_\infty \twoheadrightarrow R_{\cl r, S}^{\square_T, \eta}$, where $R_\infty$ is essentially (a finite product of) $\FDRx$ adjoining some formal variables. Once we have $r = r (\pi)$, we construct a space $S(\sigma)$ of automorphic forms with prescribed local behavior at $v$, specified by the representation $\sigma$ of the compact group $\mathtt{K} \subset D_v^\times \cong \GL_n(K).$ This way we obtain an $R_{\cl r, S}^{\square_T, \eta}$-module $S(\sigma)$. Then applying the standard Taylor--Wiles--Kisin method (in the function field setup), we build a module $S_\infty (\sigma)$ over $R_\infty$, by patching a family of modules $S_{N}(\sigma)$ over global deformation rings $R^{\square_T, \eta}_{T}$, for $N \in \Z_{\ge 1}$, where deeper and deeper tame ramifications are allowed at certain auxiliary Taylor--Wiles primes $Q_N$. Finally we check that the patched module $M_\infty(\sigma)$ has the desired support property because of the prescribed local conditions given by $\sigma$. 
  
Our proof of Conjecture \ref{conj:intro} in the function field case is simplified compared to the number field case because of the non-existence of infinite places. Moreover, L. Lafforgue's result on the full global Langlands correspondence for $\GL_n$ in the function field case is available. This together with the absence of $p$-adic Hodge theory makes the patching argument more transparent in our setup as well. We carry out the now standard computation of the numerical criterions needed for the patching argument in Section \ref{sec:Global} and \ref{sec:Patching}, in order to convince the readers that everything indeed works out.

\begin{remark*}
The condition $p \nmid n (n-1)$ is assumed mostly for simplicity, and can be removed by more careful patching method introduced in \cite{Thorne}. For this we will need to use the modified notion of ``adequate'' and modify the patching arguments accordingly. We choose not to work with this generality so the proofs presented in \ref{sec:Patching} are less technically involved. The condition $p \nmid l - 1$, on the other hand, is more serious. This is because in the argument we will need a global field $F$ such that at some place $v_0$, the size $l_0$ of the residue field $k_0$ satisfies $p \nmid l_0 - 1$, in order for there to exist a subgroup of $(D \otimes \A_{F})^\times$ which is ``sufficiently small'' in the sense of  Definition \ref{def:s_small}. This is easily achieved if $F$ is a number field without any extra condition, but if $F/\F_l (t)$, it seems necessary to require $p \nmid l - 1$. 
\end{remark*}

Finally, we compare the Breuil--M\'ezard conjecture over local function fields and local number fields, and show that they are compatible in a sense which can be made precise. In particular, we show that the conjecture for function fields can be deduced from the case of number fields. 

\begin{mainthm} (Theorem \ref{thm:number_function})
 Fixed $n, l, p$ where $l \ne p$. If the Breuil--M\'ezard conjecture (Conjecture \ref{conj:intro}) holds for all local number fields $K' /\Q_l$ and all mod $p$ representations $\cl \rho': G_{K'} \ra \GL_n(\F)$, then it also holds for all local function fields $K/\F_l$ and $\cl \rho : G_{K} \ra \GL_n(\F)$. 
\end{mainthm}

This is proved using the theory of Deligne--Kazhdan on pairs of ``close fields'', to compare the ``moduli'' of the representation theory side as well as the deformation rings in both cases. Two local fields $K'/\Q_l$ and $K/\F_l (\!(t)\!)$ are $m$-close to each other if $\mO_K' /(\varpi_K')^m \cong \mO_K /(\varpi_K)^m$ where $\omega_K'$ and $\omega_K$ are their uniformizers.  A typical example of $m$-close fields are $\Q_l (l^{1/m})$ and $\F_l (\!(t^{1/m})\!)$.


Note that, by using the main result of \cite{S}, we can deduce a slightly weaker form of Theorem \ref{thm:intro}, but with the restriction on $p\nmid (l - 1)$ removed. 
\begin{corollary*}[Corollary \ref{cor:number_function}] \label{maincor:1}
Suppose that $p > 2$, then Conjecture \ref{conj:intro} holds (with $\FDR$ replacing $\FDRx$ everywhere), for all local function fields $K/\F_l$ where $l \ne p$. 
\end{corollary*}

In fact, we can run this argument in the reverse direction, and recover part of Shotton's result in \cite{S} (now with specified characters), for representations which are trivial on large compact open subgroups. More precisely, let $\textup{R}_{L}^{(m)}(\mathtt{K})$ be the subgroup of  $\textup{R}_{L}(\mathtt{K})$ generated by irreducible smooth representations of $\mathtt{K}$ which factor through $\mathtt{K}_m := \GL_n(\mO_K/\varpi_K^m)$.   

\begin{corollary*} \label{mainthm:3}
Let $p \ne l$ as usual and suppose that $p \ne n(n-1)(l-1).$ Let $K'/\Q_l$ be a local number field with sufficiently large ramification index $e$. Let $m$ be an integer such that $2^{n-1} \cdot m \le e$.  Then there exists a cycle map $\cl{\textup{cyc}}: \textup{R}_{\F}^{(m)}(\mathtt{K}) \ra \mZ(\FDRx/\lambda)$ making the following diagram commute:
\[
\begin{tikzcd}
  \textup{R}^{(m)}_{L} (\mathtt{K}) \arrow[d, swap, "\textup{red}"] \arrow[r, dashed, "\cl{\textup{cyc}}"]   &\mZ(\FDRx)  \arrow[d, "\textup{red}"] \\
  \textup{R}^{(m)}_{\F} (\mathtt{K}) \arrow[r, dashed, "\cl{\textup{cyc}}"]  & \mZ(\FDRx/\lambda) 
\end{tikzcd}
\]
\end{corollary*} 
 
\subsection{Organization of the paper}

In Section \ref{sec:G_D} we recall Galois deformations and define the deformation spaces $\spec \FDRx \subset \spec \FDR$ as well as their subspaces $\spec \FDRtx, \spec \FDRt$ with fixed inertia types. We show that $\FDRx$ is reduced and flat over $\mO$. 
In Section \ref{sec:Rep} we review the relevant representation theory for $\GL_n (K)$ and $\GL_n (\mO_K)$, including the identification between the tempered Bernstein spectrum $B^t(\GL_n)$ and  the set $\Igaln$ of Galois inertia types.  We then define the generic multiplicity map. This allows us to define cycle maps $\textup{cyc}$ and $\textup{cyc}^{\chi}$ and formulate the Breuil--M\'ezard conjecture in Section \ref{sec:M_T}. In this section we also give an outline of the proof of the main theorem. In Section \ref{sec:Global} we compute the dimension of global deformation rings and apply the global Langlands correspondence over function fields to realize $\cl \rho$ as a suitable local component of a global modular representation $r$. The main result of this section is Proposition \ref{thm:global_realization}. In Section \ref{sec:Patching} we modify the now standard Taylor--Wiles--Kisin patching method over the function field, and apply it to our setup to prove Theorem \ref{thm:Main}. Finally, in Section \ref{sec:Recover_Number}, we apply Deligne--Kazhdan's theory of close fields to our problem between the local function field setup and local number field setup, and deduce Theorem \ref{maincor:1} and \ref{mainthm:3} in the introduction. 

\subsubsection*{Acknowledgement} 
The idea of looking at the function field analog of the Breuil--M\'ezard conjecture was suggested to me by Mark Kisin, and it is a great pleasure to thank him for introducing this area of mathematics to me and for his encouragement throughout this project. I also want to thank Jack Shotton and Richard Taylor for their suggestions to look at the relation between the case of number fields and the case of function fields, which lead to the discussion in Section \ref{sec:Recover_Number}.

\section{Galois deformations} \label{sec:G_D}

\subsection{Galois deformation rings}  We first fix some notations used in the article. 

\be
\item Let $K$ be a nonarchimedean local field of residue characteristic $l$. 
\footnote{All statements in this subsection hold for any non-archimedean local field of residue characteristic $l$, where $l \ne p$. We will only specialize to function fields starting from Section \ref{sec:Global}.}
\item $F$ will denote a global field, and $S$ a finite set of places of $F$. 
\item Let $L$ be a finite extension of $\Q_p$, where $p \ne l$, with the ring of integers $\mO = \mO_L $, and residue field $\F = \mO / \lambda $.  
\item For this section $G$ is a topological group satisfying Mazur's $\Phi_p$-finiteness condition. In the remaining of the article we will let $G = G_{F, S}$ or $ G = G_K$.  
\item Let $\cl \rho$ be a continuous representation  
$$ \cl \rho: G \ra \GL_n (\F). $$ 
\ee

\subsubsection{Deformation rings}    \label{def:deformation_functor} 

Let $\mC = \mC_{\mO}$ be the category of \textit{local Artin $\mO$-algebras with residue field $\F$}.  Recall the Galois   
deformation functors $\textup{D}_{\cl \rho}^{\square} $ and $\textup{D}_{\cl \rho}$.
\bi
\item The framed deformation functor $\textup{D}_{\cl \rho}^{\square}: \mC \longrightarrow Sets$ is  
\begin{align*}
 \textup{D}_{\cl \rho}^\square (A) &: = \{\rho: G \ra \GL_n (A) : \rho \textup{ reduces to } \cl \rho \} 
\end{align*}
\item Likewise, the (unframed) deformation  $\textup{D}_{\cl \rho}: \mC \longrightarrow Sets$ is  
$$ \textup{D}_{\cl \rho}  (A) : = \textup{D}_{\cl \rho}^\square (A)  / \textup{conj.\;by } \ker (\GL_n(A) \ra \GL_n (\F))  $$
\ei
 $\FD$ is pro-representable by a framed deformation ring $\FDR$ with the universal framed deformation (cf. \cite{K1})
$$\rho^{\square}: G \ra \GL_n (\FDR).$$
 Under the assumption that $\cl \rho$ is absolutely irreducible, or more generally that $\End_{\F[G]} \cl \rho = \F$, the functor $D_{\cl \rho}$ is pro-representable by the universal deformation  
$$\rho^{\textup{univ}}: G \ra \GL_n (\textup{R}_{\cl \rho}).$$ 
In this case, we have  
$ \FDR  \cong \textup{R}_{\cl \rho} [\![ x_{ij} ]\!]_{i, j \in \{1, ..., n\}} / (x_{11}).$
In particular, $\FDR$ and $\textup{R}_{\cl \rho}$ ``share the same singularities". We choose to work with the framed deformation ring $\FDR$ in this article. 

\subsubsection{Fixing the determinant} 

Let $\chi: G \ra \mO^\times$ be a continuous character such that $\chi \equiv \det \cl \rho \mod \lambda$, then there exists a quotient $\FDRx$ of $\FDR$ parametrizing framed deformations of $\cl \rho$ with determinant $\chi$.  From the functoriality of framed deformations, there is a natural map $ \Lambda := \textup{R}_{\det \cl \rho}^{\square} \ra  \FDR$, and $\spec \FDRx$ is the fiber of $\spec \FDR$ over the point $\chi: \Lambda \ra \mO$ of $\spec \Lambda$. 
As notation suggests, $\Lambda$ is independent of $\cl \rho$. 

\subsubsection{The generic fiber $\spec \FDRx [1/p]$}  

For this discussion let $G = G_K$. The following lemma is a slight variant of Lemma 1.3.2 in \cite{BLGGT}, taking into account the determinant condition.   
Recall that a geometric point on $\spec \FDR$ is called \textit{non-degenerate} if  its corresponding admissible irreducible representation of $\GL_n(K)$ via the local Langlands correspondence is generic \textup{[}cf. Subsection \ref{sss:LLC}\textup{]}. 

\bl \label{lemma:BLGGT_132}
The generic fiber $\spec \FDRx[1/p]$ of $\spec \FDRx$  is generically formally smooth of dimension $n^2 - 1$. More precisely, 
\be
\item If $x: \spec \cl L \ra \spec \FDRx$ is a non-degenerate geometric point, then $x$ lies on a unique irreducible component of $\spec \FDRx$, and $\spec \FDRx[1/p]$ is formally smooth at $x$.
\item The non-degenerate geometric points are Zariski-dense in $\spec \FDRx$.
\ee
\el

\bproof 
The proof in \cite{BLGGT} goes through with fixed determinant $\chi$ -- with $\textup{ad}^{0} \cl \rho$ replacing $\textup{ad} \cl \rho$ everywhere. \eproof  
 
\subsection{Galois and Weil--Deligne inertia types} 

From now on we fix a mod $\lambda$ representation $\cl \rho: G_K \ra \GL_n (\F)$. 
In this subsection we define $\textup{R}^{\square, \tau}_{\cl \rho}$ (resp. $\textup{R}^{\square, \tau, \chi}_{\cl \rho}$) for a Galois inertia type $\tau$. Under the condition that $L$ is large enough, they are unions of irreducible components in $\spec \FDR$ (resp. $\spec \FDRx$). 

\bd  Let $V$ be a vector space over $\cl \Q_p$ of dimension $n$. We distinguish the following two notions of inertia types. \footnote{The inertia types used in \cite{BM}, \cite{K3} and \cite{EG} are called WD inertia types here, while the inertia types defined in \cite{S} are called Galois inertia types in our setup.} 
\be
\item A Weil--Deligne (WD) inertia type $\xi$ is an isomorphism class of representations $\xi: I_K \ra \GL (V)$ which extends to a Frobenius-semisimple Weil--Deligne representation of $W_K$. In particular, $\xi$ has finite image and carries no information about the monodromy operator $N$. 
\item  A Galois inertia type $\tau$ is an isomorphism class of continuous representation $\tau: I_K \ra \GL(V)$ which extends to a continuous representation $\rho: W_K \ra \GL(V)$, where the image of Frobenius is semisimple.  
\ee
We denote by $\Iwdn$ (resp. $\Igaln$) the set of WD (resp. Galois) inertia types of dimension $n$. 
\ed

\br \label{remark:Galois_inertia_type_finer}
Fix a tame character $t =t_p: I_K \ra \Z_p^\times$ for $K$ and a lift of arithmetic Frobenius $\varphi \in G_K$. We have an equivalence of categories between continuous Frobenius-semisimple representations and Frobenius-semisimple Weil--Deligne representations (of the same dimension), given by
$$\rho \mapsto \textup{WD}(\rho) = (\rho', N), \quad \textup{where } \rho' (\sigma \varphi^m) := \rho(\sigma \varphi^m) \exp (- t(\sigma) N),$$
where $N$ is the (unique) monodromy operator such that $\rho(\sigma) = \exp (N t(\sigma))$ for $\sigma \in I _K' \subset I_K$ for some finite extension $K'/K$. For a Galois inertia type $\tau$ which extends to $\rho: W_K \ra \GL(V)$, let $\textup{WD}(\rho) = (\rho', N)$, then the restriction $\rho' |_{K}$ gives rise to a WD inertia type, which is well defined up to isomorphism classes and we denote this isomorphism class by $\textup{WD}(\tau)$. This gives a surjective map 
$$\textup{WD}: \Igaln \ra \Iwdn. $$
\er 

The Galois (resp. WD) inertia type $\tau(x)$ (resp. $\xi (x)$) of a geometric point $x: \spec \cl L \ra \spec \FDR$ is the Galois (resp. WD) inertia type of its associated Galois (resp. Weil--Deligne) representation. 

\bd Let $\tau$ be a Galois inertia type. We define $\spec \FDRt$ \textup{\big(}resp.  $\spec \FDRtx$\textup{\big)} to be the reduced subscheme of $\spec \FDR$  \textup{\big(}resp. of $\spec \FDRx$\textup{\big)} obtained by taking the Zariski closure of $\cl L$-points of Galois inertia type $\tau$. For a Weil--Deligne inertia type $\xi$, we define $\textup{R}_{\cl \rho}^{\square, \xi}$ and $\textup{R}_{\cl \rho}^{\square, \xi, \chi}$ similarly. 
\ed 

\br As the next lemma shows, an $\cl L$-point of $\spec \FDR$ (resp. of $\spec \FDRx$) lies on $\spec \textup{R}_{\cl \rho}^{\square, \xi}$  \textup{\big(}resp.  $\spec \textup{R}_{\cl \rho}^{\square, \xi, \chi}$\textup{\big)} if and only if it has WD inertia type $\xi$.  The subscheme $\spec \FDRt$ \textup{\big(}resp.  $\spec \FDRtx$\textup{\big)} on the other hand, may contain $\cl L$-points which are not of Galois inertia type $\tau$, so having a fixed Galois type is really a generic condition (even on the generic fiber) and in particular not a local deformation condition. 
\er 

Nevertheless we have the following rigidity lemma.
\bl \label{lemma:BLGGT_134}
Let $\cl \rho$ be as above. 
\be
\item The Weil--Deligne inertia type $\xi$ is constant on geometrically connected components of $\spec \FDRx[1/p]$. 
\item The Galois inertia type is generically constant on geometrically irreducible components of $\spec \FDRx [1/p]$. More precisely, let $x$ and $y$ be non-degenerate geometric points on the same irreducible component of $\spec \FDRx \otimes_{\mO} \cl \Q_p$, then $\tau (x) = \tau(y)$. 
\ee
\el
The corresponding claims for $\spec \FDR[1/p]$ without fixing the character $\chi$ still hold, and are proved in Lemma 1.3.4 of \cite{BLGGT} (see also part (4) of Proposition 3.6 of \cite{S}). 

\bproof  The proof of (1) is the same as the standard argument for the case without fixing $\chi$ and is standard. 
Part (2) follows from a slight modification of Lemma 1.3.4 of \cite{BLGGT}, where we replace $\FDR$ by $\FDRx$. Let $\fm$ be a maximal ideal of $\FDRx \otimes \cl \Q_p$, then following the same proof we may exhibit a surjection 
$$\phi: (\FDRx \otimes \cl \Q_p)^{\wedge}_{\fm} \twoheadrightarrow \cl \Q_p [\![x_1,..., x_{n^2 - 1}]\!].$$
The rest of the proof there only uses the fact that $\dim \FDRx [1/p] = n^2 - 1$. Therefore, for any two points $x$ and $y$ on the same irreducible component of $\spec \FDRx \otimes \cl \Q_p$, such that neither $x$ or $y$ is contained in any other irreducible component, we have $\tau (x) = \tau(y)$. This in particular applies to non-degenerate points by part (1) of Lemma \ref{lemma:BLGGT_132}. 
\eproof

\subsection{Flatness} 

The following lemma will be useful for us. Shotton \cite{S} attributes the proof to Helm, which we briefly reproduce in the beginning of the proof of Proposition \ref{lemma:deformation_character_flat} because later we will need it in Section \ref{sec:Recover_Number}. 

\bl[Theorem 2.5, \cite{S}]  \label{lemma:deformation_ring_flat} 
The space of framed deformations $\spec \FDR$ is a reduced  complete intersection, and flat over $\spec \mO$  of relative dimension $n^2$. 
\el

Before we state a corollary, let us record a simple commutative algebra fact. 

\bl \label{lemma:flat_over_Z_p}
Let $R$ be a reduced noetherian ring over $\mO$, then $\spec R$ is flat over $\spec \mO$ if and only if every irreducible component of $\spec R$ is flat over $\spec \mO$. \footnote{In general an irreducible component inside a flat family may fail to be flat over the base. }
\el 
\bproof  
First suppose $R$ is flat over $\mO$ and let $\wp$ be a minimal prime of $R$, then we claim that $R/\wp$ is $\lambda$-torsion free. Suppose for the contrary there exists $x \notin \wp$ such that $\lambda x \in \wp$, then $\lambda \in \wp$. Let $\wp_1, ..., \wp_r$ be the other minimal primes of $R$ distinct from $\wp$, then $\cap \wp_i \ne 0$ for otherwise $\wp \supset \wp_i$ for some $i$. Let $y \in \cap \wp_i$ be a nonzero element, then $\lambda y \in (\cap \wp_i) \cap \wp$ is $0$ since $R$ is reduced, which contradicts the fact that $R$ is $\lambda$-torsion free.

Now suppose for all minimal primes $\wp \subset R$, $R/\wp$ is $\lambda$-torsion free. Let $x \in R$ be element such that $\lambda x = 0$, then $x \in \wp$ for all minimal primes $\wp$, therefore by reducedness  $R$ is $\lambda$-torsion free. 
\eproof 

\bc \label{cor:inertia_component}
Let $\tau$ \textup{(}resp. $\xi$\textup{)} be a Galois  \textup{(}resp. WD \textup{)} inertia type. 
\be
\item The irreducible components of $\FDR$ are flat over $\mO$. 
\item Suppose that the coefficient field $L$ is sufficiently large, then each nonempty $\FDRt$  \textup{(}resp. $\textup{R}_{\cl \rho}^{\square, \xi}$\textup{)} is a union of irreducible components of $\FDR$. 
\item Under the assumptions above, each nonempty $\FDRt$  \textup{(}resp. $\textup{R}_{\cl \rho}^{\square, \xi}$\textup{)} is flat over $\mO$ of dimension $n^2$. 
\ee
\ec

\bproof 
Part (1) follows from Lemma \ref{lemma:deformation_ring_flat} and Lemma \ref{lemma:flat_over_Z_p}.  For part (2), we prove the claim for Galois inertia type $\tau$, the argument for the case of Weil--Deligne inertia type works the same. We define $\spec \FDRt [1/p]$ to be the reduced subscheme of $\spec \FDR [1/p]$ by taking the Zariski closure of the $\cl L$-points of type $\tau$. 
We claim that $\FDRt [1/p]$ is a union of irreducible components of $\FDR [1/p]$, because on each irreducible component of $\FDR [1/p]$, all the non-degenerate points have the same Galois inertia type and they are Zariski dense by Lemma \ref{lemma:BLGGT_132} and Lemma \ref{lemma:BLGGT_134}. Therefore $\FDRt$ is the Zariski closure of a union of irreducible components of $\FDR [1/p]$. Taking Zariski closure commutes with taking finite unions, therefore the claim follows. 
For part (3), note that both $\FDRt$ and $\textup{R}_{\cl \rho}^{\square, \xi}$ are reduced quotients of $\FDR$ by definition, now apply Lemma \ref{lemma:deformation_ring_flat}. 
\eproof

\br  \label{remark:irreducible_component_inertia}  
In the main theorem of the article, the field $L$ will be irrelevant as long as it is sufficiently large so that each irreducible component is geometrically irreducible, which we always assume. See Remarks \ref{remark:well_defined_chi_multiplicity_2}, \ref{remark:relation_with_fixing_determinant}, and Subsection \ref{sss:invariance_base_change}.  In particular, the residue field $\F$ of $L $ always has cardinality $\ge 5$ (since $p \ne 2$ by assumption, this condition simply excludes $\F_3$), which is convenient for us in the proof of Proposition \ref{cor:global_realization}. 
\er

\bl  \label{lemma:adding_variable_to_deformation_ring} Let $\cl \rho$ and $\chi$ be as before, and assume that $p \nmid n$. 
\be
\item Let $\xi$ be a WD inertia type such that that the closed subscheme $\spec \textup{R}_{\cl \rho}^{\square, \xi, \chi} \subset \spec \textup{R}_{\cl \rho}^{\square, \xi}$ is non-empty, then we have 
$$ \textup{R}_{\cl \rho}^{\square, \xi} \cong \textup{R}_{\cl \rho}^{\square, \xi, \chi}  [\![T]\!] $$
where the universal deformation is given by $\rho^{\square, \xi} (g) = \rho^{\square, \xi, \chi} (g) \cdot \textup{diag}(1 + T)$. 
\item Let $\tau$ be a Galois inertia type such that $\spec \FDRtx$ is non-empty, then 
$$ \FDRt \cong  \FDRtx  [\![T]\!].$$
\ee 
\el 
 
\bproof  
From definition there is a natural map $\textup{R}_{\cl \rho}^{\square, \xi} \ra  \textup{R}_{\cl \rho}^{\square, \xi, \chi}  [\![T]\!]$. To show that it is an isomorphism, note that for any $L'$ point $\rho: G_K \ra \GL_n(L')$ of $\FDR$, the character $ \mu = \chi \det (\rho)^{-1}$ is unramified and reduces to identity mod $\lambda'$ where $\lambda'$ is a uniformizer of $\mO_{L'}$. Since $p \nmid n$, $n$ is invertible in $\cl \F$, so there exists an $n^{th}$ root $\nu$ of $\mu$ such that $\mu = \nu^n$ and $\nu$ reduces to identify mod $\lambda'$. The representation $\rho \otimes (\nu \circ \det) $ therefore has determinant $\chi$. This proves part (1) of the lemma, part (2) is similar.  
\eproof

\bc \label{cor:deformation_ring_character_union}
Suppose that $p \nmid n$ and that $L$ is sufficiently large. Let $\tau$ be a Galois inertia type, assume that  $\spec \textup{R}_{\cl \rho}^{\square, \tau, \chi}$ is non-empty. Then $\spec \textup{R}_{\cl \rho}^{\square, \tau, \chi} $ is a union of irreducible components of  $\spec \FDRx$, and is flat over $\spec \mO$ of relative dimension $n^2 - 1$.
\ec  
 
\bproof 
The first claim follows from the same proof of Corollary \ref{cor:inertia_component} (2).  Then Lemma \ref{lemma:adding_variable_to_deformation_ring} tells us that $\spec \textup{R}_{\cl \rho}^{\square, \tau, \chi} [\![T]\!]$ is $\lambda$-torsion free. 
\eproof 
 
However, this does not directly imply that $\FDRx$ is flat. Nevertheless, following the same strategy of the proof of Lemma \ref{lemma:deformation_ring_flat}, we have the following 

\bp \label{lemma:deformation_character_flat} 
 Suppose $p \nmid n$ and fix a character $\chi$. $\spec \FDRx$ is a reduced complete intersection and flat over $\spec \mO$ of relative dimension $n^2 - 1$. 
\ep

\bproof (1).  First let us prove the lemma assuming that $\cl \rho$ is tamely ramified. Let $\sq P_{K} \n I_K$ be the kernel of the tame character $t_p: I_K \twoheadrightarrow \prod_{r \ne l}\Z_r \twoheadrightarrow \Z_p$, we may replace $G_K \ra \GL_n(\F)$ by $\cl \rho: G_K/\sq P_K \ra \GL_n (\F)$.  This allows us to view $\FDRx$ as the completion of the a certain local ring of a certain moduli space which we now describe. 

Let $R$ be any commutative ring, let $r \ge 2$ be an integer and let $\alpha \in R^\times$ be a unit. We define $\mM(n, r, \alpha)_R$ to be the following moduli space over $R$: for $S$ an algebra over $R$, the $S$-points of $\mM(n, r, \alpha)_R$ are 
$$\{(\Sigma, \Phi) \in \GL_n (S) \times \GL_n (S)  \quad \vline \quad   \Phi \Sigma \Phi^{-1} = \Sigma^r, \det (\Phi) = \alpha\}.$$
Therefore, $\mM(n, r, \alpha)_R$ is a closed sub-scheme of $\GL_{n, R} \times_{R} \GL_{n, R}$ cut out by $n^2 + 1$ equations. Let 
$$\pi: \mM(n, r, \alpha)_R \ra \GL_{n, R} $$ be the projection given by $(\Sigma, \Phi) \mapsto (\Sigma).$ Now in our setup, when $R = \mO$ and $r = \#|k|$, we claim that $\mM(n, r, \alpha)_{\mO}$ is a local complete intersection, and is flat over $\mO$ of relative dimension $n^2 - 1$. 
  
For notational simplicity let $\fX = \mM(n, r, \alpha)_\mO$ and $\fX' = \GL_{n, \mO} \times \GL_{n, \mO}$. Let $\fX_s \subset \fX'_s$ be the corresponding special fibers over $\spec \F$. $\fX \subset \fX'$ is a closed sub-scheme cut out by $n^2 + 1$ equations inside the regular scheme $\fX'$ of relative dimension $2 n^2$ over $\mO$. Upper semi-continuity of the map $x \mapsto \dim (\fX_{f(x)})$ says that $\{x \in \fX: \dim \fX_{f(x)} = n^2 \}$ is open in $\fX$. Note that  every point in $\fX$ specializes to a point in the special fiber $\fX_s$ (since every point in $\fX'$ specializes to a point in $\fX_s'$). Therefore, to prove that $\fX$ is a local complete intersection it suffices to show that the geometric special fiber $\fX_{\cl s}  = \mM (n, r, \cl \alpha)_{\cl \F}$ has dimension $n^2 - 1$.   
 
Next we show that $\mM(n, r, \cl \alpha)_{\cl \F}$ is a local complete intersection of dimension $n^2 - 1$. For this we consider the morphism $\pi: \mM(n, r, \cl \alpha)_{\cl \F} \ra \GL_{n, \cl \F}$ and let $\Sigma_0 \in \im (\pi)$ be a point in the image of $\pi$. Let $Z_{\Sigma_0}$ be the centralizer of $\Sigma_0$ in $\SL_{n, \cl \F}$ and $C_{\Sigma_0}$ the $\SL_{n}$ conjugacy class of $\Sigma_0$, note that $C_{\Sigma_0} \cong \GL_{n, \cl \F}/Z_{\Sigma_0}$. It is clear that $C_{\Sigma_0} \subset \im (\pi)$ and that the pre-image $\pi^{-1} (C_{\Sigma_0})$ is a $Z_{\Sigma_0}$-torsor over $C_{\Sigma_0}$.  In particular,  $\pi^{-1} (C_{\Sigma_0})$ has dimension  $\dim C_{\Sigma_0} + \dim Z_{\Sigma_0} = n^2 - 1$.  Finally, note that there exists an integer $m \ge 1$, such that for any $\Sigma \in \im (\pi)$,  the eigenvalues of $\Sigma$ are $(r^m - 1)^{th}$ roots of unity. This implies that the number of $\SL_n$ conjugacy classes in $\im (\pi)$ is finite. To see this, note that since $p \nmid n$, one can always conjugate a given matrix to its Jordan normal form by elements in $\SL_n(\cl \F)$. Thus $\dim (\mM (n, r)_{\cl \F}) = n^2 - 1$.  This completes the proof that $\fX$ is a complete intersection, in particular Cohen-Macaulay, and that $\fX$ is flat over $\mO$ (since the fibral dimensions are constant over the regular ring $\mO$).

(2). The proof of the general case relies on the following tame reduction lemma of \cite{CHT}.  Let $\delta: \sq P_K \ra \GL_{d_\delta} (\F)$  be an irreducible representation of $\sq P_K$. Consider the open subgroup  
$$G_{\delta}:= \{g \in G_K: \delta^{g} \sim \delta\}$$ of $G_K$ of finite index. 
 $\delta$ has a unique deformation to $\sq \delta: \sq P_K \ra \GL _{d_\delta} (\mO)$, which extends to $G_\delta$. In this discussion we fix this extension and again denote it by $\sq \delta$. For any finite $\mO$-module $M$ with $G_K$ action, we define the $\mO[T_\delta] := \mO[G_\delta/\sq P_K]$ module $M_\delta: = \Hom_{\sq P_K} (\sq \delta, M).$ Lemma 2.4.12. of \cite{CHT} states that 
$$M \cong \bigoplus_{[\delta] \in \textup{Conj}} \textup{Ind}_{G_\delta}^{G_K} (\sq \delta \otimes M_\delta) $$ 
where the direct sum is over $G_K$ conjugacy classes of irreducible $\F[\sq P_K]$-modules. (cf. \cite{CHT} 2.4.12. and \cite{Notes} for a detailed proof).  As above, write $\cl \rho_{\delta} = \Hom_{\sq P_K} (\sq \delta, \cl \rho)$, then there is an isomorphism of deformation functors 
$$ \textup{D}_{\cl \rho} \isom \prod_{[\delta] \in \textup{Conj}} \textup{D}_{\cl \rho_\delta},$$ where the inverse map is given by $(\rho_\delta)_{[\delta]} \mapsto \oplus_{[\delta]} \textup{Ind}_{G_\delta}^{G_K} (\sq \delta \otimes \rho_\delta)$. By carefully taking into account the effect of lifting frames, we in fact have 
$$\FDR \cong \big( \widehat{\otimes}_{[\delta] \in \textup{Conj}} \textup{R}_{\cl \rho_{\delta}}^{\square} \big) [\![x_1, ..., x_{n^2-\sum n_{\delta}^2} ]\!]. $$
We need to add the determinant condition. To proceed we fix the following choice:
\bi
\item Let $\varphi \in G_K$ be a lift of the arithmetic Frobenius, also fix a choice of $\sigma \in I_K$ whose image is a topological generator of $I_K/\sq P_K \cong \Z_p$.  
\item Fix a representative $\delta $ in each conjugacy class $[\delta]$. Let $K_\delta := (\cl K)^{G_\delta}$ which is a finite extension over $K$, and let $e_{\delta}$ (resp. $f_\delta$) be the inertia degree (resp. residue degree) of $K_{\delta}/K$.
\item $I_{K_\delta}/\sq P_K \subset I_K/\sq P_K$ is of index $e_\delta$, so $\sigma_\delta:= \sigma^{e_\delta}$ is a topological generator of the tame quotient in $I_{K_{\delta}}$.  Next, choose $0 \le i_\delta \le e_\delta -1$ such that $\varphi_\delta := \sigma^{i_\delta} \cdot \varphi^{f_\delta} \in G_\delta$. $\varphi_\delta$ is a lift of arithmetic Frobenius in $G_\delta$. 
\item Let $\{\sigma^i \varphi^j\}_{\substack{0 \le i \le e_\delta -1 \\ 0 \le j \le f_\delta - 1}}$ be a coset representative of $G_K/G_{\delta}$. 
\ei
This allows us to compute the determinant of $\textup{Ind}_{G_\delta}^{G_K} (\sq \delta \otimes \rho_\delta)$ for an $n_\delta$-dimensional representation $\rho_\delta$, by using the formula (cf. \cite{CR})
$$ \det \textup{Ind}_{G_\delta}^{G_K} (\sq \delta \otimes \rho_\delta) (\varphi) = \textup{Sgn}_\delta \cdot \prod_{s \varphi t^{-1} \in G_\tau} \det ((\sq \delta \otimes \rho_{\delta} )(s \varphi t^{-1})) $$
where the product is taken over $s, t \in \{\sigma^i\varphi^j\}$ in the coset representative chosen above, and $\textup{Sgn}_\delta = (\textup{Sign of action of } \varphi \textup{ on } G_K/G_\delta)^{n_\delta}$. Write $s = \sigma^i \varphi^j$ and $t = \sigma^m \varphi^k$, then $s \varphi t^{-1} \in G_\delta$ precisely when one of the following occurs 
\be
\item $k = j + 1$, $0 \le j \le f_\delta - 2$. In which case, $0 \le i = m \le e_\delta - 1$. 
\item $k=0$ and $j = f_\delta - 1$; and $i \equiv m \cdot r^n + i_\delta \mod e_\delta$, for $0 \le i, m \le e_\delta - 1$. 
\ee
Therefore, when we compute the determinant of the induced representation at $\varphi$, only $s, t$ in the second case above contributes, and we have 
\begin{align*}
 & \det \textup{Ind}_{G_\delta}^{G_K} (\sq \delta \otimes \rho_\delta) (\varphi)  \\ 
 = \quad & \textup{Sgn}_\delta \cdot \prod_{s \varphi t^{-1} \in G_\tau} \det ((\sq \delta \otimes \rho_{\delta} )(s \varphi t^{-1})) \\
 = \quad  & \textup{Sgn}_{\delta} \prod_{i - m r^n \equiv i_\delta \mod e_\delta} \Big((\det \sq \delta (\sigma^{i - m r^n} \varphi^{f_\delta}))^{n_\delta} (\det \rho_{\delta} (\sigma^{i - m r ^n - i_{\delta}}))^{d_{\delta}} (\det \rho_{\delta} (\varphi_\delta))^{d_{\delta}} \Big)  \\
 = \quad & \beta_{\delta} \cdot (\det \rho_{\delta} (\sigma_\delta))^{c_\delta d_\delta} \cdot (\det \rho_{\delta} (\varphi_\delta))^{e_\delta d_\delta} 
\end{align*}
where $c_{\delta} = \frac{1}{e_\delta} \sum (i - m r^n - i_\delta)$, where the sum is taken over $i, m \in \{0, ..., e_\delta - 1\}$ such that $i \equiv m r^n + i_\delta \mod e_\delta$; and 
$$\beta_{\delta} := \textup{Sgn}_{\delta} \prod_{i - m r^n \equiv i_\delta \textup{ mod } e_\delta} (\det \sq \delta (\sigma^{i - m r^n} \varphi^{f_\delta}))^{n_\delta}  \in \mO^\times. $$
Both $c_\delta \in \Z$ and $\beta_\delta \in \mO^\times$ only depends on $\delta$.  

Now let us define a subscheme $\mM(n, r, \cl \rho, \alpha)_R \subset \prod_{[\delta]} (\GL_{n_\delta, R}\times_R \GL_{n_\delta, R})$ where the product is indexed by $G_K$ conjugacy classes of irreducible $\F[\sq P_K]$-modules,  to consist of tuples  $(\Sigma_\delta, \Phi_\delta)_{[\delta]}$ satisfying 
\be
\item $\Phi_{\delta} \Sigma_{\delta} \Phi_{\delta}^{-1} = \Sigma_{\delta}^{r_\delta}$ 
for each $[\delta]$. 
\item Let  $\beta (\cl \rho) := \prod_{[\delta] \in \textup{Cong}} \beta_{\delta}$, then we require that
$$\beta(\cl \rho) \cdot \prod_{[\delta]}\Big(\det(\Sigma_\delta)^{c_\delta d_\delta} \cdot \det(\Phi_\delta)^{e_\delta d_\delta} \Big)= \alpha.$$
\ee
In other words, 
        \begin{multline*} \mM(n, r, \cl \rho, \alpha)_R:= \Big\{ (\Sigma_\delta, \Phi_\delta)_{[\delta]} 
\:\: | \:\: \Phi_{\delta} \Sigma_{\delta} \Phi_{\delta}^{-1} = \Sigma_{\delta}^{r_\delta}, \\
         \textup{ and } \prod_{[\delta]} \Big(\det(\Sigma_\delta)^{c_\delta d_\delta} \cdot \det(\Phi_\delta)^{e_\delta d_\delta} \Big) = \alpha \cdot \beta(\cl \rho)^{-1} \Big\}.
        \end{multline*}
As the tame case, we want to show that $\fX(\cl \rho)_{\cl \F} := \mM(n, r, \cl \rho, \alpha)_{\cl \F}$ has dimension $\sum_{[\delta]} n_\delta^2 - 1$. The proof is similar: we first define $G_{\underline{n}} \n \prod_{[\delta]} \GL_{n_\delta} (\cl \F)$ to be the normal subgroup consisting of elements $(g_{[\delta]})$ such that $\prod_{[\delta]} \det(g_{[\delta]}) = 1$. $G_{\underline{n}}$ acts on $\prod_{[\delta]} \GL_{n_\delta} (\cl \F)$ by pairwise conjugation. Now consider
$$ \pi:  \mM(n, r, \cl \rho, \alpha)_{\cl \F} \ra \prod_{[\delta]} \GL_{n_\delta} (\cl \F).$$
Let $\Sigma_{0} = (\Sigma_{0, \delta})_{[\delta]}$ be an element in $\im (\pi)$, let $Z_{\Sigma_0}$ be the centralizer of $\Sigma_0$ in $G_{\underline{n}}$, and let $C_{\Sigma_0}$ be the $G_{\underline{n}}$-conjugacy class of $\Sigma_0$. Then the determinant condition on the pre-image $\pi^{-1} (C_{\Sigma_0})$ becomes 
$$ \prod_{[\delta]} \det (\Phi_{\delta})^{e_\delta d_\delta} = \gamma_{0} $$
where $\gamma = \alpha \cdot \beta(\cl \rho)^{-1} \cdot \prod_{[\delta]} (\det \Sigma_{0, \delta})^{- c_{\delta} d_\delta}$, which only depends on the conjugacy class $C_{\Sigma_0}$ (and $\cl \rho$). This shows that $\pi^{-1} (C_{\Sigma_0})$ is again a $Z_{\Sigma_0}$-torsor over $C_{\Sigma_0}$, therefore has dimension $-1 + \sum_{[\delta]} n_\delta^2 $. Now the same argument as the tame case tells us that $\mM(n, r, \cl \rho, \alpha)_{\mO}$ is a complete intersection, and flat over $\mO$. 

(3). Finally, given our choices of $\sigma, \varphi \in G_K$ and $\sigma_\delta, \varphi_\delta \in G_\delta$, the representation $\cl \rho$ and character $\chi$ provides us a closed point $x_{\cl \rho}$ on $\mM(n, r, \cl \rho, \alpha)_{\mO}$, where $\alpha := \chi (\varphi)$. The completion of $\mM(n, r, \cl \rho, \alpha)_{\mO}$ at $x_{\cl \rho}$ is precisely $\FDRx$, since on each $\pi^{-1}(C_{\Sigma_0})$, $\chi |_{I_K}$ is fixed. Therefore, $\FDRx$ is a complete intersection and flat over $\mO$. Now reducedness follows from generically reducedness (which implies the $R_0$ condition) and being Cohen-Macaulay (which in particular implies Serre's $S_1$ condition). 
\eproof

\subsection{Cycles and reduction} 

\subsubsection{Cycles}  \label{sss:cycles}

Let $X$ be a quasi-compact noetheorian scheme locally of finite type\footnote{All schemes in this subsection are quasi-compact, noetherian, and locally of finite type. }.   An (algebraic) cycle on $X$ is a formal sum 
$\sum n_Z [Z]$, taken over integral closed subschemes of $X$, such that the set $\{Z: n_Z \ne 0\}$ is finite.  Let $\mZ_d(X)$ be the group of $d$-dimensional algebraic cycles on $X$. If $X$ is equi-dimensional, then we often suppress the subscript and write $\mZ (X) = \mZ_{\dim X} (X)$. 
Let $\mF$ be a coherent sheaf on $X$, then $\textup{Supp} (\mF)$ is a closed subscheme of $X$ with its reduced subscheme structure. Let $Z \subset \textup{Supp} (\mF)$ be an irreducible component with generic point $\xi$, then the multiplicity 
$$ e(\mF, Z) := \textup{length}_{\mO_{X, \xi}} \mF_{\xi} < \infty$$
is finite. For an integral closed subscheme $Z'$ not contained in $\textup{Supp} (\mF)$, we set $e(\mF, Z') = 0$. Recall that 
\be 
\item For $d \ge \dim \textup{Supp} (\mF)$, the $d$-dimensional cycles associated to $\mF$ is 
$$[\mF]_{d} := \sum e(\mF, Z) [Z] \in \mZ_d (X)$$
where the sum is over all integral closed subscheme of dimension $d$. Again for equi-dimensional $X$, we write $Z(\mF) = [\mF]_{\dim X}$. 
\item 
For a closed subscheme $Y \subset X$ of dimension $\le d$, define $[Y]_d := [\mO_Y]_{d}. $
\ee

\subsubsection{Reduction of cycles}  Now suppose that $X$ is a reduced scheme over $\mO = \mO_L$ of equi-dimension $d$.  Let 
$$X_{s} = X/\lambda  := X \times_{\spec \mO} \spec \F$$
be its mod $\lambda$ special fiber. Let $i: X_s \ra X$ be the closed immersion.

\bd  \label{def:reduction_of_cycles}
The reduction map (on irreducible components of $X$) 
$\textup{red}: \mZ (X) \ra \mZ(X_s)$ is given as follows. For each irreducible component $Y$ (which is a $d$-cycle of $X$), define $\textup{red}: [Y] \mapsto Z(i^* \mO_Y) = [i^* \mO_Y]_{d - 1}. $
\ed

In general, if $f: Y \ra X$ is a not flat, then pullbacks of cycles along $f$ usually do not behave well. Therefore we will need the following lemma
 
\bl \label{lemma:reduction_cycle_mod_p}
 Let $X$ be a reduced scheme over $\mO$ of equi-dimension $d$. Furthermore assume that $X$ is flat over $\mO$. Let $\mF$ be a $\lambda$-torsion free coherent sheaf on $X$, then 
$$ \textup{red} ([\mF]_{d}) = [i^* \mF]_{d - 1} \in \mZ (X_s).$$ In other words, $\textup{red} (Z(\mF)) = Z (\mF/\lambda) \in \mZ (X_s)$. 
\el

\bproof 
This follows from Proposition 2.2.13 of \cite{EG}. By taking $f = \lambda$ in their proposition, we know that if $[\mF]_d =\sum_{i \in I} n_i  Z_i  $ then
$$[i^* \mF]_{d - 1} =  \sum_{i \in I} n_i  (Z_i \cap X_s ) =  \textup{red} ([\mF]_{d}). $$ 
The flatness of $X /\mO$ is needed to ensure that $\lambda$ is a non-zero divisor in each stalk of $\mO_X$, in order to apply their proposition. 
\eproof 

\subsubsection{Product of cycles} We will need the following lemma on irreducible components in completed tensor product. 

\bl \label{lemma:tensor_product_cycle}
Let $R$ and $S$ be two complete local Noetherian $\mO$-algebras with geometrically integrable irreducible components. Further assume that $R$ and $S$ are equi-dimensional, and that they either (1). both $\lambda$-torsion free or (2). both $\lambda$-torsion. Then $\mZ(R \hat \otimes S) = \mZ (R) \otimes \mZ (S)$. 
\el

\bproof 
This is Lemma 2.10 in \cite{S}. 
\eproof

\section{Representation theory of $\GL_n (K)$}  \label{sec:Rep}

Throughout the article we fix the following notations 
\bi
\item Let $\textup{Rep}^{\infty}_n = \textup{Rep}^{\infty}_n (K)$ be the \textit{category} of smooth representations of $\GL_n(K)$.
\item Let $\mA_n = \mA_{n} (K)$ be the \textit{set of isomorphic classes} of irreducible admissible representations of $\GL_n(K)$. 
\item  Let $\mA_{n}^{sc} \subset \mA_{n}$ (resp. $\mA_{n}^{d.s.}, \mA_{n}^{t}, \mA_{n}^{g}$) be the subset of classes of supercuspidal (resp. square integrable, tempered, generic) representations.  
\ei

\subsection{Local Langlands for $\GL_n$} \label{section:BZ_classification}
 
\subsubsection{The Bernstein--Zelevinsky classification.}  \label{sss:BZ} 
We first fix notations used in the Bernstein--Zelevinsky classification. 

Let $(\pi, V) \in \mA_{n}$, for any $s \in \C$ let $\pi(s)$ be the twist of $\pi$ by the character $|\!\cdot\!|^s \circ \det$.  Recall that a (Zelevinsky) segment in $\mA_{n}$ is an ordered collection of $m$ supercuspidal representations of $\GL_n(K)$ of the form: 
$$\Delta (\pi, m) = [\pi, \pi(1), ..., \pi (m-1)].$$
Two segments  $\Delta_1 = \Delta (\pi_1, m_1)$ and $ \Delta_2= \Delta (\pi_2, m_2)$ in $\mA_{n_1}^{sc}$  and $\mA_{n_2}^{sc}$ are linked if $\Delta_1 \nsubseteq \Delta_2$, $\Delta_2 \nsubseteq \Delta_1$, and $\Delta_1 \cup \Delta_2$ is a segment.  $\Delta_1$ precedes $\Delta_2$ if $\Delta_1$ and $\Delta_2$ are linked, and $\pi_1 < \pi_2 $. 
 For a segment $\Delta= \Delta (\pi, m)$ in $\mA_n$, we define an admissible representation of $\GL_{mn}$ by
$$\pi(\Delta) : = \pi_1 \times \cdots \times \pi_m,$$ where 
$\pi_1 \times \cdots \times \pi_m := \textup{Ind}_{P_{(n,...,n)}}^{\GL_{mn}} (\otimes \pi_i)$ is the normalized parabolic induction, where $P_{(n, ..., n)}$ denotes  the standard parabolic corresponding to the partition $mn = n + \cdots + n$.  Recall from \cite{BZ} that for each segment $\Delta$,  the representation $\pi (\Delta)$ admits a unique irreducible quotient-representation $Q(\Delta).$ If $\Delta_1, ..., \Delta_r$ are segments such that for each $ i < j,  \:  \Delta_i$  does not precede $\Delta_j$, then the representation 
$$Q(\Delta_1)\times \cdots \times  Q(\Delta_r) \in \textup{Rep}^{\textup{adm}} ({\GL_{\sum n_i m_i}}) $$
admits a unique irreducible quotient $Q(\Delta_1, ..., \Delta_r)$.  Any irreducible representation $(\pi, V)$ of $\GL_n(K)$ is isomorphic to one of the form $\pi \cong  Q(\Delta_1, ..., \Delta_r),$ where $\Delta_i = (\pi_i, m_i) \textup{ with } n = \sum n_i m_i $, 
for a unique collection of segments $\Delta_1, ... \Delta_{r}$ up to permutation. 

 Moreover, the representation $Q(\Delta_1) \times \cdots \times Q(\Delta_r) $ is irreducible if and only if none of the segments $\Delta_i$ and $\Delta_j$ are linked. The (quotient) representation  $Q(\Delta_1, ..., \Delta_r)$ is generic if and only if none of the segments are linked; it is tempered if and only if it is generic and the central character of each $\pi_i (\frac{m_i-1}{2})$ is unitary; it is square integrable if and only if it is tempered and there is only $1$ segment (so $\pi \cong Q(\Delta)$). 

\subsubsection{Steinberg representations} \label{sss:St_rep} 
 
 The Steinberg representation of $\GL_n(K)$ is 
$$\textup{St}_n:=Q(\Delta) \:\:\: \textup{ where } \Delta = \Delta(|\cdot|^{\frac{1-n}{2}}, n).$$ 
More generally, we can twist $\textup{St}_n$ by a quasi-character $\chi$ and define $\textup{St}_n(\chi):=Q(\Delta), \textup{ where } \Delta = \Delta(\chi |\cdot|^{\frac{1-n}{2}}, n).$

The (twisted) Steinberg representations $\textup{St}_n (\chi)$ correspond to special representations $\textup{Sp}_n (\chi)$ on the Galois side. Let us recall the definition. 
\bd \label{def:special_reps}
Let $V$ a finite dimensional vector space over $\cl \Q_p$.  Let $(r, V)$ be a representation of $W_K$, we define the special representation $\textup{Sp}_m (r)$ of $W_K$ on $V^m$, to be the Weil--Deligne representation
$$\textup{Sp}_m (r):= r \oplus r(1) \oplus ... \oplus r(m-1)$$ 
where $r(i)$ is the representation $r$ twisted by $\chi_{\textup{cyc}}^{-1}$, 
and $N \in \End(V^m)$ is the operator sending each $r(i)$ isomorphically to $r(i+1)$ for $i \le m-2$ and $N = 0$ on $r(m-1).$  
\ed

\subsubsection{The Local Langlands correspondence (\textup{LLC})} \label{sss:LLC} Let $\mG_n (K)$ be the set of isomorphism classes of n-dimensional Frobenius-semisimple Weil--Deligne representations of $W_K$ over $\cl \Q_p$.\footnote{we have implicitly fixed an isomorphism $\iota: \cl \Q_p \isom \C$.} By the work of \cite{LRS} for the function field case (and \cite{HT} for the number field case), there exists a unique collection of bijections
$$\textup{rec}_n: \mA_n (K)\isom \mG_n (K)$$ satisfying a list of compatibilities conditions \footnote{we will not make use of the other compatibilities: namely it is compatible with local class field theory, with taking contragradients, and with $L$- and $\epsilon$-factors for every pair  $\pi_1, \pi_2$.}, in particular: 
\bi
\item for any quasi-character $\chi$,  $ \textup{rec}_n(\pi \chi) = \textup{rec}_n(\pi) \otimes \textup{rec}_1 (\chi) $; and 
\item for any $\pi \in \mA_n$ with central character  $\omega_\pi$, 
$\det \circ \textup{ rec}_n (\pi) = \textup{rec}_1 (\omega_\pi)$.  
\ei 
Under the local Langlands correspondence,  
\be
\item $\pi \in \mA_n$ is supercuspidal if and only if $\textup{rec}_n (\pi)$ is irreducible; 
\item $\pi$ is essentially square integrable if and only if $\textup{rec}_n (\pi)$ is indecomposable.  
\item Let $\pi = Q(\Delta_1, ..., \Delta_r)$ with $\Delta_i = \Delta (\pi_i, m_i)$, then 
$$\textup{rec}_n (\pi)  = \bigoplus_{i = 1}^{r} \textup{Sp}_{m_i}( \textup{rec}_{n_i} (\pi_i) ).$$ 
\ee

\subsection{Inertia correspondence} \label{sec:Bernstein_decomposition}  
In this subsection we recall the notion of Bernstein decomposition and a tempered Bernstein decomposition due to \cite{SZ}. Roughly speaking, each Bernstein component corresponds to a WD inertia type, while each tempered Bernstein component corresponds to a Galois inertia type. 

\subsubsection{Bernstein components and decomposition}  

Let $P \subset \GL_n(K)$ be (the $K$-points of) a parabolic subgroup with Levi factor $L$.  Let $\varsigma$ be an irreducible supercuspidal representation of $L$. A Bernstein component 
$$R^{(L, \varsigma)}  =R^{(L, \varsigma)} (\GL_n(K)) \subset \textup{Rep}^{\infty}_{n}$$ is the full subcategory consisting of $\pi$ such that every irreducible subquotient of $\pi$ is a subquotient of $\textup{I}_{L}^{\GL_n} (\varsigma \otimes \eta):= \textup{Ind}_P^{\GL_n(K)}(\varsigma \otimes \eta)$, where $\eta$ is an unramified character of $L$ and $\textup{I}_{L}^{\GL_n} := \textup{Ind}_P^{\GL_n(K)} \circ j^*$ is the normalized parabolic induction, with $j: P \ra L$ the map to the Levi quotient.  
Two Bernstein components $R^{(L_1, \varsigma_1)}$ and $R^{(L_2, \varsigma_2)}$ are equal if and only if they are inertially equivalent, meaning that there exists an unramified character $\eta$ of $L_2$  and $g \in \GL_n(K)$ such that 
$$g L_1 g^{-1} = L_2, \qquad {\;}^g\varsigma_1 = \varsigma_2 \eta $$
The Bernstein spectrum of $\GL_n(K)$ is 
$$B(\GL_n(K)) := \{ (L, \varsigma)\}/ \textup{inertial equivalence}.$$
A well-known theorem of Bernstein \cite{Be} states that the category $\textup{Rep}^{\infty}_{n}$ of smooth representations of $\GL_n(K)$ decomposes as a direct product 
$$\textup{Rep}^{\infty}_{n} = \prod_{(L, \varsigma) \in B(\GL_n)} \textup{R}^{(L, \varsigma)}. $$ In particular, this provides a partition on $\mA_n$. 

\bl \label{lemma:Bernstein_component_inertia}
Two irreducible admissible representations $\pi$ and $\pi'$ of $\GL_n(K)$ belong to the same Bernstein component $R^{(L, \varsigma)} $ if and only if the underlying $W_K$ representations under the local Langlands correspondence satisfy $\textup{rec} (\pi) |_{I_K} \cong \textup{rec} (\pi') |_{I_K};$ in other words, if and only if they give rise to the same WD inertia type.   
\el

\bproof 

The Bernstein--Zelevinsky classification (cf. Subsection \ref{sss:BZ}) allows us to reduce to the supercuspidal cases.  For this we need to show that two supercuspidal representations $\pi$ and $\pi'$ of $\GL_{m} (K)$  differ by a twist of an unramified character $\eta$ if and only if $\textup{rec}_m (\pi)|_{I_K} \cong \textup{rec}_m (\pi')|_{I_K}$. This follows from Subsection \ref{sss:LLC} and the following consequence of Schur's lemma: 

\textup{Sublemma. } Let $H \n G$ be a normal subgroup with an abelian quotient, and let $\rho$ and $\rho'$ be  two finite dimensional representations  of $G$ over an algebraically closed field $k$, then $\rho|_{H} \cong \rho'|_{H}$ if and only if  there exists a character $\psi: G \ra G/H \ra k^\times$ such that $\rho \cong \rho'\psi$.  
\eproof

\subsubsection{Tempered Bernstein decomposition} The Bernstein decomposition has been refined by Schneider and Zink \cite{SZ} in the case of tempered representations. Consider the tempered Bernstein spectrum $B^t (\GL_n)$, which by definition is the set 
$$B^t (\GL_n) := \{(L, \vartheta)\}/ \textup{Conjugation},$$ 
where $L$ is a Levi of $\GL_n(K)$ and $\vartheta$ is an irreducible square integrable representation of $L$.  Note that we no longer quotient by inertia equivalence since after twisting by unramified characters $\vartheta$ may not be unitary. From the Bernstein--Zelevinsky description of tempered and square integrable representations, the normalized parabolic induction gives a bijection between 
$$ B^t(\GL_n) \isom \mA^{t}_n. $$ 
Let $\textup{Rep}^t_n \subset \textup{Rep}^\infty_n$ be the full subcategory of tempered representations. For $(L, \vartheta) \in B^t(\GL_n)$, let $R^{(L, \vartheta)}$ to be the full subcategory of $\textup{Rep}^t_n$ whose irreducible sub-quotients are induced from $(L, \vartheta)$. Then the Bernstein decomposition refines to 
$$\textup{Rep}^t_n = \prod_{(L, \vartheta) \in B^t(\GL_n)} R^{(L, \vartheta)}.$$

\subsubsection{Inertia types and Bernstein spectrums}
Recall from Section \ref{sec:G_D} that $\Iwdn$ (resp. $\Igaln$) is the set of WD (resp. Galois) inertia types. Let $\textup{I}^0$ be the set of Galois types which extend to an irreducible representation of $W_K$.  

\bd By a finite partition we mean a finite ordered collection of (non-strictly) decreasing non-negative integers.  Let $\Part$ be the set of finite partitions. Let the degree of $(a_1, a_2, ... ) \in \Part$ be $\sum a_i$, hence $x \in \Part$ is a partition of $\deg x$. 
\be
\item Define $\mP^\infty = \{\xi:  \textup{I}^0 \ra \Z_{\ge 0}  \: | \: \xi \textup{ has finite support}\}.$ \item Define $\mP^t = \{\tau:  \textup{I}^0 \ra \Part  \: | \: \tau \textup{ has finite support}\}.$
\item For each $\xi \in \mP^\infty$ (resp. $\tau \in \mP^t$), define its degree to be 
$$\deg (\xi) := \sum_{\xi_0 \in \textup{I}^0} \xi(\xi_0) \qquad \Big( \textup{resp. } \deg (\tau) := \sum_{\xi_0} \deg (\tau(\xi_0)) \Big).$$
\item Define $\mP^\infty_n$ (resp. $\mP^t_n$) to be the subset of $\mP^\infty$ (resp. $\mP^t$) consisting of elements of degree $n$. 
\ee
\ed

We have the following maps of sets. 
\bi
\item $B(\GL_n) \ra \mP^\infty_n$. Let $(L, \varsigma)$ be an element in the Bernstein spectrum, for each $\tau_0 \in \textup{I}^0$, which corresponds to a unique supercuspidal representation $[\pi_0]$ up to inertia equivalence, we associate to it the multiplicity of $[\pi_0]$ in the inertia equivalence class of the supercuspidal support of $(L, \varsigma).$
\item $B^t (\GL_n) \ra \mP^t_n$. Let $(L, \vartheta)$ be an element of $B^t(\GL_n)$, so $\vartheta$ corresponds to a collection of square integrable representations $Q(\Delta(\pi_i, m_i))$ with $\sum m_i = n$. For each $[\pi_0]$ as above, we send it to the partition consisting of all $m_i$ (in decreasing order) such that $[\pi_i] = [\pi_0]$ (meaning that they differ by an unramified twist). 
\item $\mP^\infty_n \ra \Iwdn$. For $\xi \in \mP_n^\infty$,we associated a WD inertia type by $\underset{\tau_0 \in \textup{I}^0}{\oplus} \tau_0^{\xi (\tau_0)}$.
\item $\mP^t_n \ra \Igaln$. For $h \in \mP^t_n, $ we associate a Galois inertia type as follows: for each $\tau_0 \in \textup{I}^0$, choose an irreducible extension $\sq \tau_0$ of $\tau_0$ to $W_K$, and then take the Galois representation associated to the Weil--Deligne representation 
$$\bigoplus_{\tau_0 \in \textup{I}^0} \Big( \bigoplus_i \textup{Sp}_{h(\tau_0)(i)} (\sq \tau_0) \Big)$$
where $h(\tau_0)(i) $ is the $i^{th}$ term in the partition $h(\tau_0)$. 
\ei

\bl \label{lemma:types_and_centers}
These maps fit together into the following commutative diagram, where all horizontal maps are bijections. 
\[
\begin{tikzcd}
B^t(\GL_n) \arrow[r, "\sim"] \arrow[d, two heads, "\mathfrak{Z}"] &\mP^t_n \arrow[r, "\sim"] \arrow[d, two heads, swap, "\deg"]  & \Igaln \arrow[d, two heads, swap, "\textup{WD}"]
  \\
B(\GL_n) \arrow[r, "\sim"] &\mP^\infty_n \arrow[r, "\sim"] & \Iwdn 
\end{tikzcd} 
\]
\el 
\bproof 
The bottom horizontal maps are all bijections by Lemma \ref{lemma:Bernstein_component_inertia} and its proof. For the second bijection on the top row, note that by Remark \ref{remark:Galois_inertia_type_finer}, the associated Galois representations of two Weil--Deligne representations 
$$ \bigoplus_{i = 1}^{r} \textup{Sp}_{m_i}(\rho_{i}), \quad  \bigoplus_{i = 1}^{r'} \textup{Sp}_{m_i'}(\rho_{i}')$$
have the same restrictions on $I_K$ if and only if $r =r'$ and up to reordering $m_i = m_i'$ and $\rho_i|_{I_K} \cong \rho_i'|_{I_K}$. This shows that all the horizontal maps are bijections. 

The map $\mathfrak{Z}: B^t(\GL_n) \ra B(\GL_n)$ is the map sending $(L, \vartheta)$ to its supercuspidal support, and the diagram commutes by unwinding definitions. 
\eproof 

\br \label{remark:tempered_support_generic_rep}
$B^t(\GL_n)$ can also be regarded as the ``generic Bernstein spectrum", in the sense that it gives a stratification of generic representations.  A generic representation has the form $\pi = Q(\Delta_1, ..., \Delta_r)$, $\Delta_i = \Delta (\pi_i, m_i)$, and for each $i$ there exists a unique unramified character $|\cdot|^{s_i}$ for some $s_i \in \C$ such that $Q(\Delta (\pi_i |\cdot|^{s_i}, m_i))$ has unitary central character (hence square integrable). 
Therefore we have a map 
$\mA^{g}_n \ra B^t(\GL_n)$ which sends $\pi$ to $(L = L_{m_1,..., m_r}, \vartheta)$, where $\vartheta$ is the product of $Q(\Delta(\pi_i |\cdot|^{s_i}, m_i))$ described above. Unwinding definitions, the map 
$\mA^g_n \ra B^t(\GL_n) \ra \Igaln$ is given by taking the Galois representation associated to $\textup{rec}_n (\pi)$ and restrict to $I_K$. We summarize this discussion with the following nice diagram 
\[
\begin{tikzcd}
\mA^t_n  \arrow[r,  "\subset"]    \arrow[rd, swap, "\cong"]  & \mA^g_n   \arrow[r,  "\subset"]  \arrow[d, two heads] & \mA_n \arrow[d, two heads]
  \\
& B^t(\GL_n) \arrow[r, two heads, "\mathfrak{Z}"] & B(\GL_n) 
\end{tikzcd} 
\]
where the middle vertical arrow is the map defined above, and the other vertical arrow is given by the usual stratification of irreducible representations by the Bernstein spectrum. Finally, for any $s \in B^t(\GL_n)$ let $R^s \subset \mA^g_n$ be the pre-image of $s$. As in the classical case, we say two irreducible generic representations lie in the same tempered Bernstein component if they belong to the same $R^s$. 
\er
 
\subsubsection{Inertia correspondence for generic representations} 

Let  $\mathtt{K} = \GL_n(\mO_K)$ be the standard maximal compact open subgroup of $\GL_n(K)$. 

\bl \label{lemma:define_pi_tau}

Let $\pi, \pi'$ be two generic irreducible admissible representations of $\GL_n(K)$ in the same tempered Bernstein component of $\GL_n(K)$, or equivalently $ \textup{rec}_n (\pi)|_{I_K}  \cong  \textup{rec}_n (\pi')|_{I_K}. $
Then they restrict to isomorphic representations on $\mathtt{K}$.   
\el

\bproof   
For any (standard) parabolic $P \subset \GL_n(K)$ and a representation $\sigma$ of $P$, we have 
$(\textup{Ind}_P^{\GL_n(K)} \sigma)|_{\mathtt{K}} = \textup{Ind}_{P \cap \mathtt{K}}^{\mathtt{K}} (\sigma|_{P \cap \mathtt{K}}). $
Since $\pi$ is generic, by  Subsection \ref{sss:BZ} and Subsection \ref{sss:LLC}, $\pi = Q(\Delta_1) \times \cdots \times Q(\Delta_r)$ and $\textup{rec}_{n} (\pi) = \bigoplus_{i =1}^{r} \textup{rec}_{n_i m_i} (Q(\Delta_i))$ and similarly for $\pi'$. 

Similar to the proof of Lemma \ref{lemma:types_and_centers}, it suffices to prove the lemma under the assumption that $\pi = Q(\Delta(\lambda, m))$ and $\pi' = Q(\Delta(\lambda', \pi))$, where $\lambda$ and $\lambda'$ are supercuspidal representations of $\GL_{d}(K)$. By Lemma \ref{lemma:Bernstein_component_inertia}, $\lambda = \lambda' \otimes \eta$,  where $\eta$ is an unramified character of $\GL_d(K).$ The character $\eta$ is of the form $\eta = \chi \circ \det$ where $\chi$ is an unramified quasi-character $K^\times \ra \C^\times$. Therefore, $\textup{rec}_n(\pi) = \textup{rec}_n(\pi') \otimes \textup{rec}_1(\chi)$, from which it follows that $\pi = \pi' \chi$. 
\eproof 

\br
The generic condition is needed. For example, for $\GL_2(K)$, the trivial representation $\textup{triv}$ and an unramified principal series $\chi_1 \times \chi_2$ lie in the same Bernstein component, but they have different restrictions to $\mathtt{K} = \GL_n(\mO_K)$.   
\er

\subsection{The generic multiplicity map} \label{sec:G_M}

Next we define the generic multiplicity map $m$. This will allow us to define cycle maps $\textup{cyc}$ and $\textup{cyc}^{\chi}$ in Section \ref{sec:M_T}, where we also formulate the Breuil--M\'ezard conjectures and our main theorem.   

\subsubsection{Multiplicities}  Recall that $\mathtt{K} = \GL_n(\mO_K)$, and fix the following convention.  
\bi
\item  Let $\textup{Rep}_{L}  (\mathtt{K}) = \textup{Rep}^{\textup{f.l.}}_{L} (\mathtt{K})$ \textup{\big(}resp. $\textup{Rep}_{\F}  (\mathtt{K})$ \textup{\big)} be the category of smooth representations of $\mathtt{K}$ of finite length over $L$ \textup{\big(}resp. over $\F$ \textup{\big)}. 
\item Let $\textup{R}_{L} (\mathtt{K})$  \textup{\big(}resp. $\textup{R}_{\F} (\mathtt{K})$ \textup{\big)}  be the Grothendieck group of the abelian category  $\textup{Rep}_{L}  (\mathtt{K})$ \textup{\big(}resp. $\textup{Rep}_{\F}  (\mathtt{K})$ \textup{\big)}. 
\ei

\br \label{remark:semisimplicity_rep}  
 For smooth representations of $\mathtt{K}$, finite length is the same as finite dimensional. Moreover, the category $\textup{Rep}_{L} (\mathtt{K})$ is semi-simple, since finite dimensional smooth representations have open kernel, $\mathtt{K}$ acts each $\sigma \in \textup{Rep}_L(\mathtt{K})$ through a finite quotient.  
On the other hand, however, $\textup{Rep}_{\F}(\mathtt{K})$ is not semi-simple. 
 
\er

\bd  \label{def:pi_tau}
Let $\tau \in \Igaln$ be a Galois inertia type, we define an isomorphism class of admissible representation $\pi_\tau$ of $\mathtt{K}$ as follows. 
Let $\sq \pi$ be any generic irreducible admissible representation of $\GL_n (K)$ such that $\textup{rec}_n (\sq \pi) |_{I_F} \cong \tau,$ then let $\pi_\tau := \sq \pi |_{\mathtt{K}}.$
This is well defined up to isomorphism by Lemma \ref{lemma:define_pi_tau}. 
\ed 

The following definition is a variant of Definition 4.2 of \cite{S}. 

\bd   Let $\sigma \in  \textup{Rep}_{L}  (\mathtt{K})$, and let $\tau$ be an inertia type. Fix $\chi: G_K \ra \mO^\times$ a lifting of $\det \cl \rho$. 

\be
\item Define the multiplicity $m (\sigma, \tau)$ to be 
$$m (\sigma, \tau) := \dim  \Hom_{\cl L [\mathtt{K}]} (\sigma, \pi_\tau) = \dim  \Hom_{ \mathtt{K}} (\sigma \otimes_{L} \cl L, \pi_\tau).$$
\item Define the $\chi$-multiplicity $m (\sigma, \tau, \chi)$ to be 
$$m (\sigma, \tau, \chi) := \begin{cases} m (\sigma, \tau)  & \textup{if there exists } \sq \pi \textup{ in Def. } \ref{def:pi_tau} \textup{ with } \textup{rec}_1(\omega_{\sq \pi}) \cong \chi
\\  0  & \textup{otherwise.} \end{cases}$$ 
Here $\omega_{\sq \pi}$ denotes the central character of $\sq \pi$. 
\ee 
\ed 
In the definition above,  $\Hom_{\mathtt{K}} (\sigma_0, \pi|_{\mathtt{K}})$ is a finite dimensional $\cl L$-vector space. Indeed, the action of $\mathtt{K}$ factors through $\mathtt{K}/\mathtt{K_0}$ where $\mathtt{K_0} < \mathtt{K}$ is a compact open subgroup of finite index. 
Since $\sigma_0|_{\mathtt{K}_0}$ is the trivial representation of $\mathtt{K}_0$, we have 
$$\dim_{\cl L} \Hom_{\mathtt{K}} (\sigma_0, \pi|_{\mathtt{K}}) \le \dim_{\cl L}\Hom (\sigma_0 |_{\mathtt{K_0}}, \pi|_{\mathtt{K_0}})^{\mathtt{K_0}} = \dim_{\cl L} (\pi^{\oplus \dim \sigma_0})^{\mathtt{K_0}}.$$
Now $\pi^{\mathtt{K}_0}$ is finite dimensional.  
 
\br \label{remark:well_defined_chi_multiplicity}  
For a fixed $\tau$ and $\chi$, the multiplicities $m(\sigma, \tau)$ and $m(\sigma, \tau, \chi)$ are both additive on short exact sequences in $\textup{Rep}_{L}  (\mathtt{K})$ by semi-simplicity. Therefore it is well defined on the Grothendieck group $\textup{R}_{L} (\mathtt{K})$. Now fix $\sigma \in \textup{R}_{L} (\mathtt{K})$ and a Galois  inertia type $\tau$. Let $\chi, \chi'$ be two liftings of $\det \cl \rho$.  $$m (\sigma, \tau, \chi) =  m (\sigma, \tau, \chi') \:\: \iff \: \: \chi |_{I_K} \cong  \chi' |_{I_K}. $$
The forward implication follows directly from definition of $m(\sigma, \tau, \chi)$ and compatibility between determinants and central characters (\ref{sss:LLC}). The other direction follows from the same argument near the end of the proof of Lemma \ref{lemma:define_pi_tau}. 
\er
\br \label{remark:well_defined_chi_multiplicity_2}  
 The coefficient field $L$ is irrelevant in the definition of $m(\sigma, \tau)$.  
\er 

\subsubsection{The generic multiplicity map} \label{sss:multiplicity_map} 

In this discussion we organize the information of the multiplicities into a \textit{multiplicity map} $m_{L}$ from $\textup{R}_{L} (\mathtt{K}) $ to the direct sum $\underset{s \in B^t}{\oplus} \Z$ indexed by the tempered Bernstein spectrum of $\GL_n(K)$. Recall from Remark \ref{remark:tempered_support_generic_rep}, each tempered Bernstein component contains at least one irreducible generic representation, among which there is a unique tempered one. 

\bd \label{def:multiplicity_map} 
The \textit{generic multiplicity} map is a homomorphism
$$m_L:  \textup{R}_{L} (\mathtt{K}) \longrightarrow \underset{s \in B^t(\GL_n)}{\prod} \Z, \qquad \sigma \longmapsto \big( m (\sigma, \tau_s) \big)_s $$
where $\tau_s \in \Igaln$ is given by Lemma \ref{lemma:types_and_centers}.  In other word, for any $\sigma \in \textup{R}_{L} (\mathtt{K})$, define the $s$-component of $m_L (\sigma) \in \prod_{s \in B} \Z$ to be $m(\sigma, \tau_s)$.  
\ed 

By Remark \ref{remark:semisimplicity_rep}, $m_L$ is compatible with base change of the coefficient field $L$. We will write $m = m_{\cl L}$ for the generic multiplicity map on $ \textup{R}_{\cl L} (\mathtt{K})$. The next proposition is probably known to experts, but we cannot find any published reference so we provide a proof here.

\bp \label{lemma:image_m_direct_sum}
Let $\sigma$ be a smooth irreducible representation of $\mathtt{K}$, then there are only finitely many $s \in B^t(\GL_n)$ such that $m(\sigma, \tau_s) \ne 0$. 
\ep

Before proceeding to the proof of the proposition, we need two lemmas on supercuspidal representations. Let $G^{\circ} \n \GL_n(K)$ be the kernel of 
$$\GL_n (K) \xrightarrow{\det} K^\times \xrightarrow{\textup{val}} \Z $$ where $\textup{val}$ is the normalized valuation sending a uniformizer $\varpi_K$ to $1 \in \Z$. $G^{\circ}$ therefore contains all compact open subgroups of $\GL_n(K)$.  Let $H$ be the subgroup $Z G^{\circ} \subset \GL_n(K)$ where $Z = Z(\GL_n (K))$ is the center of $\GL_n (K)$.  In our setup, $Z \cong K^\times \hookrightarrow \GL_n (K)$ and  $H \n \GL_n (K)$ is a normal open subgroup of index $n$. 

\bl \label{lemma:finitely_many_irreducible_supercuspidal_G_0}
Let $\mathtt{K}_0$ be a compact open subgroup of $\GL_n(K)$. There are only  finitely many isomorphism classes of irreducible supercuspidal representations of $G^{\circ}$ with a $\mathtt{K_0}$ fixed vector.
\el

\bproof 
This is the reformulation of the corollary of Proposition 3.21 in \cite{Be}. 
\eproof 

\bl
Let $\mathtt{K}_0$ be a compact open subgroup of $\GL_n(K)$. Then up to twisting by an unramified central character, there are finitely many isomorphism classes of irreducible supercuspidal representations $\pi$ of $H = Z G^{\circ}$ satisfying the following conditions:  
\be
\item $Z$ acts by a character $\chi$; 
\item $\pi$ contains a $\mathtt{K}_0$ fixed vector, namely $\pi^{\mathtt{K_0}} \ne 0 $. 
\ee 
\el

\bproof 
Let $\mathfrak S_{0}$ be the finite set of isomorphism classes of irreducible supercuspidal representations of $G^{\circ}$ with nonzero $\mathtt{K}_0$ fixed vectors as in Lemma \ref{lemma:finitely_many_irreducible_supercuspidal_G_0}.  Now consider an irreducible supercuspidal representation $\pi$ of $\GL_n(K)$ satisfying the two conditions in the lemma. Since $Z$ acts by a character, the restriction of $\pi$ to $G^{\circ}$ as a supercuspidal $G^{\circ}$ representation is still irreducible. Therefore $\pi|_{G^{\circ}} = \pi_0$ for some $\pi_0 \in \mathfrak S_0$. Let $\pi'$ be another such representation of $\GL_n(K)$ with $\pi' |_{G^{\circ}} = \pi_0$. Let $\chi' := \pi'|_{Z}$.  Since $Z \cap G^\circ = \mO_K^\times \subset K^\times$, $\chi |_{\mO_K^\times} = \pi_0 |_{\mO_K^\times} = \chi'|_{\mO_K^\times}$, therefore $\chi = \chi' \eta$ for an unramified character $\eta$ of $Z = K^\times$. This shows that for any $ h = z g \in Z G^{\circ}$ with $z \in Z$ and $g \in G^{\circ}$, 
$\pi' (h)  = \eta (z) \pi (h)$, hence the lemma.  
\eproof 

\bc \label{cor:finite_supercuspidal_case}
For each compact open subgroup  $\mathtt{K}_0 < \mathtt K = \GL_n(\mO_K)$, there are a finite number of isomorphism classes of irreducible supercuspidal representations $\pi$ of $\GL_n(K)$ up to twisting by an unramified character, such that $\pi^{\mathtt K_0} \ne 0$.  
\ec

\bproof 

Let $\pi$ be such a supercuspidal representation with $\pi^{\mathtt{K}_0} \ne 0$. Consider the restriction $ \pi_H := \pi |_{H}$, by Clifford theory (for restriction of possibly infinite dimensional representation to a normal subgroup of finite index $\iota$), $\pi_H = \pi_{H, 1} \oplus \cdots \oplus \pi_{H, d}$ where $d \le \iota$ and all the $\pi_{H, i}$'s are conjugates to each other by some elements in $G$. In particular, we know for each factor $\pi_{H, i}$, the center $Z$ acts by the same character $\chi$ and that each $(\pi_{H, i})^{\mathtt{K}_0} \ne 0$. Now apply the Lemma above. 
\eproof

\bproof[Proof of Proposition \ref{lemma:image_m_direct_sum}] 

First we observe that Proposition \ref{lemma:image_m_direct_sum} is equivalent to the following statement: 
\bi
\item[(A)] \textup{For each normal open subgroup } $\mathtt{K}_0 \n \mathtt{K}$ \textup{ of finite index}, there are only finitely many isomorphism classes of irreducible tempered representations $\pi$ such that $\pi^{\mathtt{K}_0} \ne 0$ (as a representation of $\mathtt{K}/\mathtt{K}_0$). 
\ei
To see this, note that $\sigma$ has open kernel $\mathtt{K}_0 \n \mathtt{K}$ of finite index, so the representation factors through $\mathtt{K} /\mathtt{K}_0$. The equivalence then follows from Remark \ref{remark:tempered_support_generic_rep}
and the fact that $\Hom_{\mathtt{K} } (\sigma, \pi_s|_{\mathtt{K} }) = \big(\Hom(\sigma, \pi_s|_{\mathtt{K} })\big)^{\mathtt{K}} = \Hom_{\mathtt{K}/\mathtt{K}_0} (\sigma, \pi_s^{\mathtt{K}_0}|_{\mathtt{K} }).$ 
Both directions of this equivalence of claims will be used in the proof. First we reduce to the square integrable case. 
\bi
\item[(B)]    If for each $m \ge 1$ and each normal open subgroup $\mathtt{K}_0 \n \mathtt{K}$ of finite index, there are finitely many isomorphism classes of irreducible square integrable representations $\pi$ of $\GL_m (K)$ such that $\pi^{\mathtt K_0} \ne 0$, then the Proposition holds.  
\ei
  
\noindent \textit{Proof of }\textup{Claim (B). } Fix $n$ and $\sigma$ as in the proposition. 
Now consider an irreducible tempered representation 
$\pi \cong Q(\Delta_1) \times \cdots \times Q(\Delta_r) = \textup{I}_{\textup{L}}^{\GL_n} \big( \otimes Q(\Delta_i) \big)$ of $\GL_n (K)$, where $\textup{L} \cong \prod \GL_{m_i n_i} (K)$ is the standard Levi corresponding to the partition $n = \sum m_i n_i$. Let $P$ be the corresponding parabolic. By the same argument as in Lemma \ref{lemma:define_pi_tau}, we have 
$$\Big( \textup{Ind}_P^{\GL_n(K)}  \big( \otimes Q(\Delta_i)  \big) \Big) |_{\mathtt{K}} = \textup{Ind}_{P \cap \mathtt{K}}^{\mathtt{K}} \Big( \big(\otimes Q(\Delta_i) \big)_{P \cap \mathtt{K}} \Big).$$
Therefore, by adjunction of $\textup{J}_{\mathtt{K}_{\textup{L}}}^{\mathtt{K}} \dashv \textup{I}_{\mathtt{K}_{\textup{L}}}^{\mathtt{K}} $, where $\textup{J}$ denotes the Jacquet functor, and the fact that $\mathtt{K}_{\textup{L}} = \mathtt{K} \cap \GL_n(K) \cong  \prod \mathtt K_i :=  \prod_{i} \GL_{m_i n_i} (\mO_K)$, we have 
$$\Hom_{\mathtt{K}} (\sigma, \pi|_{\mathtt{K}}) = \Hom_{\mathtt{K}_L} (\textup{J}_{\mathtt{K}_L}^{\mathtt{K}} \sigma, \otimes Q(\Delta_i)|_{\mathtt{K}_L}).$$
Note that $\textup{J}_{\mathtt{K}_L}^{\mathtt{K}} \sigma$ is a finite direct sum of smooth irreducible representations of $\mathtt{K}_{L}$, so it suffices to show that for any smooth irreducible representation $\sigma_L$ of $\mathtt{K}_L$, there are only finitely many square integrable representations of $L$ of the form $Q(\Delta_1) \otimes \cdots \otimes Q(\Delta_r)$ such that $\Hom_{\mathtt{K}_L} (\sigma_L,  \otimes Q(\Delta_i)|_{\mathtt{K}_L}) \ne 0$. Equivalently (as being observed above), for any tuple $(\mathtt{K}_{i, 0} \n \mathtt K_i)$ of normal open subgroups of $\mathtt K_i = \GL_{m_i n_i} (\mO_K)$, there are finitely many such $ \otimes Q(\Delta_i)$ such that 
$$(Q(\Delta_1) \otimes \cdots \otimes Q(\Delta_r))^{\mathtt K_1 \times \cdots \times \mathtt K_r} \cong Q(\Delta_1)^{\mathtt K_1} \otimes \cdots \otimes Q(\Delta_r)^{\mathtt K_r} \ne 0. $$
But this follows from the assumption of Claim (B).  

\bi
\item[(C)] If for each $m \ge 1$ and each irreducible representation $\sigma$ of $\GL_m(\mO_K)$, the multiplicity $m(\sigma, s) \ne 0$ for only finitely many supercuspidal Bernstein components, then the assumption of Claim (B) holds.  
\ei
 
\noindent  \textit{Proof of }\textup{Claim (C). }  It suffices to show the following \textit{a priori} slightly stronger form of the assumption of Claim (B): for each given $\sigma$ as in Claim (C), there are a finite number of admissible representations of $\GL_n(K)$ of the form $\pi (\Delta) = \pi' \times \cdots \times \pi' (m-1)$ up to twisting $\Delta$ (more precisely twisting the $\pi' (i)$'s) by an unramified character, such that 
 $$\Hom_{\mathtt{K}} (\sigma, \pi (\Delta)|_{\mathtt{K}}) =  \Hom_{\mathtt K_{\Delta}} (\textup{J}_{\mathtt{K}_{\Delta}}^{\mathtt{K}} \sigma, \underset{0 \le i \le m-1}{\otimes} \pi' (i) |_{\mathtt{K}_\Delta}) \ne 0.$$
Here $\mathtt K_{\Delta}$ is the intersection of $\mathtt{K}$ with the standard Levi $L_{\Delta} = \prod \GL_{n/m}(K)$.  This follows from the same analysis as the proof of Claim (B) (note that each supercuspidal Bernstein component contains a unique one with unitary central character).  

Finally, the proposition follows from Corollary \ref{cor:finite_supercuspidal_case} (which implies the assumption of Claim (C)) and the reduction steps above.  
\eproof 

As a consequence, we obtain Definition \ref{def:intro_multiplicity} from the introduction:
\bc \label{cor:multiplicity_map}
The image of the generic multiplicity map $m$ is contained in the subgroup $\bigoplus_{s \in B^t} \Z \subset \prod_{s \in B^t} \Z$. Therefore, we get 
$$m:  \textup{R}_{\cl L} (\mathtt{K}) \longrightarrow  \bigoplus_{s \in B^t(\GL_n)} \Z. $$ The same statement holds for $m_L$ instead of $m$. 
\ec

\subsection{Tempered $\mathtt{K}$-type} In this subsection we reformulate the notion of tempered $\mathtt{K}$-types in \cite{SZ} as a certain section of the generic multiplicity map. \footnote{Their notion of tempered $\mathtt{K}$-type, however, is slightly different from the notion \cite{BK}}

\subsubsection{\empty}
For each $\xi \in B(\GL_n)$, let $B^t(\xi)$ be the pre-image $\mathfrak{Z}^{-1} (\xi) \subset B^t(\GL_n)$. In terms of Lemma \ref{lemma:types_and_centers}, $B^t (\xi)$ corresponds bijectively to $\mP^t (\xi)$ which are elements $\tau \in \mP^t_n$ that map to $\xi \in \mP^\infty_n$. There is a partial order on each $\mP^t(\xi)$ (hence on $B^t(\xi)$), induced from the dominance order on $\Part$, namely $\tau \succeq \tau'$ if and only if for all $\xi_0 \in \textup{I}^0,$  $\sum_{i} \tau(\xi_0) (i) \ge \sum_i \tau'(\xi_0) (i)$ for all $i \ge 1$. For the discussion that follows it is convenient to write $\underset{s \in B^t(\GL_n)}{\oplus} \Z = \underset{\xi \in B(\GL_n)}{\oplus} \big( \underset{\tau \in B^t(\xi)}{\oplus} \Z\big) $. For each $s \in B^t(\GL_n)$, let $1_s \in \underset{s \in B^t(\GL_n)}{\oplus} \Z $ be the element with entry $1$ at the $s$-component and $0$ everywhere else.

\bp[Proposition 6.2 \cite{SZ}] \label{prop:SZ_K_types}
 For $\tau \in B^t(\GL_n)$, there exists a smooth irreducible $\cl L$-representation $\sigma (\tau)$ of $\mathtt{K}$, such that 
\be
\item Let $\xi = \mathfrak{Z}(\tau) \in B(\GL_n)$, then $m (\sigma (\tau))$ is contained in $\underset{s \in B^t(\xi)}{\oplus} \Z \subset \underset{s \in B^t(\GL_n)}{\oplus} \Z.$
\item The $\tau$-component of $m(\sigma(\tau))$ is $1$; and more generally the $s$-component of $m(\sigma(\tau))$ is nonzero if and only if $\tau \succeq s$.  
\ee
\ep

\bproof 
This is translated from Proposition 6.2 of \cite{SZ}, or Theorem 3.7 of \cite{S} which is the same translation in terms of Galois types. 
\eproof 

\bc 
Let $\sigma: \underset{\tau \in B^t(\GL_n)}{\oplus} \Z \ra  \textup{R}_{\cl L} (\mathtt{K})$ be the map given by $1_\tau \mapsto \sigma(\tau)$ as the proposition above. Then the composition with the generic multiplicity 
$$m\sigma:  \bigoplus_{\tau \in B^t(\GL_n)}  \Z  \longrightarrow \bigoplus_{\tau \in B^t(\GL_n)}  \Z  $$ is an isomorphism. 
\ec

\bproof 
It suffices to show that for each $\xi \in B(\GL_n)$, the restriction of $m\sigma$ to 
$$m\sigma(\xi):  \bigoplus_{\tau \in B^t(\xi)}  \Z  \longrightarrow \bigoplus_{\tau \in B^t(\xi)}  \Z  $$ 
is an isomorphism. The basis of the finite free $\Z$-module $\underset{\tau \in B^t(\xi)}{\oplus} \Z $ is index by $B^t(\xi)$, which we put a total order $\succ$ as follows: start from the maximal element of the dominance partial order on $B^t(\xi)$ and follow this partial order as much as possible, until two elements are not comparable, then we order them randomly. As a result $\tau \succeq \tau'$ always implies $\tau \succ \tau'$ but not conversely. For example, for the following elements in $\Part$ with the dominance partial order
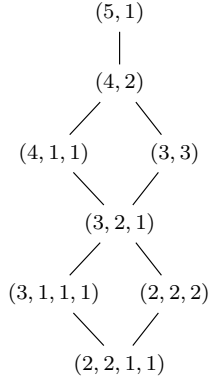
\[
{\scalefont{0.7} \begin{tikzcd} [column sep=-2em,row sep=1.5em]
& (5,1) \arrow[d, no head] \\ 
& (4,2)  \arrow[ld, no head] \arrow[rd, no head]  \\
(4, 1, 1)  \arrow[rd, no head]  &  & (3, 3)  \arrow[ld, no head]  \\
& (3, 2,1)  \arrow[ld, no head] \arrow[rd, no head]  \\
(3, 1, 1, 1)  \arrow[rd, no head]  &  & (2, 2, 2)  \arrow[ld, no head]  \\
& (2, 2,1, 1) 
\end{tikzcd}}
\]
One way to put the order $\succ$ is 
$$(5,1) \succ (4,2) \succ (3, 3) \succ (4,1, 1) \succ (3,2,1) \succ (3,1,1,1) \succ (2,2, 2) \succ (2,2,1,1).$$ 
Then by Proposition \ref{prop:SZ_K_types}, $m\sigma (\xi)$ is an upper triangular matrix with $1$'s on the diagonal, therefore invertible. This shows that $m \sigma$ is invertible. 
\eproof 

\bc \label{cor:SZ_section}
The map 
$$\sigma_v:= \sigma \circ (m\sigma)^{-1}: \bigoplus_{\tau \in B^t(\GL_n)}  \Z  \longrightarrow  \textup{R}_{\cl L} (\mathtt{K})$$
is a section of the generic multiplicity map. 
\ec 

\br  \label{remark:virtial_typical}
Let us define the following two types of sections of $m$: 
\be
\item A section $t$ is \textup{typical} if $t(1_\tau)$ is an irreducible representation of $\mathtt{K}$ for every $\tau \in B^t(\GL_n)$. 
\item A section is \textup{virtual} if otherwise. Note that $r$ is a virtual section of $m$ if and only if $r$ sends $1_\tau$ to a virtual representation for some $\tau \in B^t(\GL_n).$ 
\ee
The section $\sigma_v$ from Corollary \ref{cor:SZ_section} is a virtual section. 
\er 

\subsubsection{Tempered types}
Let $\sigma$ be a smooth irreducible representation of $\mathtt{K}$, let $s = (L, \varsigma) \in B(\GL_n)$ be an element of the Bernstein spectrum.  $\sigma$ is a type 
\footnote{Our definition of type agrees with the one given in \cite{Pas}, which is slightly different from \cite{BK}. We refer to latter notion as Bushnell-Kutzko types, which is a pair $(J, \sigma)$ where $J$ is a compact open subgroup and $\sigma$ a smooth irreducible representation of $J$ -- namely it allows more general compact open subgroups than $\mathtt{K} = \GL_n (\mO_K)$} 
for $(L, \varsigma)$  if for all smooth irreducible representations $\pi \in \textup{Rep}^\infty_n$, $\pi |_{\mathtt{K}} \textup{ contains } \sigma \iff \pi \in R^{(L, \varsigma)}.$ 
Suppose only the direction $\so$ holds, then $\sigma$ is called a $\mathtt{K}$-type. Let $\pi \in R^{(L, \varsigma)}$ be a smooth irreducible representation in $R^{(L, \varsigma)}$, a minimal $\mathtt{K}$-type for $\pi$ is a $\mathtt{K}$-type of minimal dimension occurring in $\pi |_{\mathtt{K}}$. 
The existence and uniqueness of types has been extensively studied, centering around the following folklore conjecture: for each smooth irreducible representation $\pi$ of $\GL_n(K)$,  there exists a unique minimal $\mathtt{K}$-type for $\pi$; and it occurs in $\pi|_{\mathtt{K}}$ with multiplicity $1$.  This conjecture has been proven for $\GL_2$ by Henniart as an appendix to \cite{BM}. For supercuspidal components, the existence of types essentially follows from \cite{BK}; the uniqueness and multiplicity one of types is the main result of \cite{Pas}. In view of the generic multiplicity map, we make the following definitions 
\bd \indent 

\be
\item Let $\tau \in B^t (\GL_n)$, a tempered type for $\tau$ is an irreducible representation $\sigma$ of $\mathtt{K}$ such that $m(\sigma) = 1_\tau$. 
\item Let $\xi \in B(\GL_n)$ be an element in the Bernstein spectrum, define $1_\xi$ to be $1_\xi := \sum_{\tau \in B^t(\xi)} 1_\tau  \in \underset{\tau \in B^t(\xi)}{\oplus} \Z$. A generic type for $\xi$ is an irreducible representation $\sigma$ of $\mathtt{K}$ such that $m (\sigma) = 1_\xi$. 
\ee
\ed 
The folklore conjecture above in particular implies that a unique minimal generic type exists for all $\xi$.  
 
\begin{conjecture} \label{conj:2}
 For each $\tau \in B^t(\GL_n)$, there exists a tempered type $\sigma$ for $\tau$.  
\end{conjecture}
From the definition,  Conjecture \ref{conj:2} is equivalent to the existence of  a typical section
$$\sigma_t: \bigoplus_{\tau \in B^t (\GL_n)} \Z \longrightarrow  \textup{R}_{\cl L} (\mathtt{K}).$$ 
One might even hope for a preferred such choice, which corresponds to the uniqueness of minimal tempered types.

\section{The Breuil--M\'ezard conjecture} \label{sec:M_T}  

\subsection{Statement of the conjecture} \label{ss:BM_conj} 
Recall from Section \ref{sec:G_D} that we start with a mod $\lambda$ representation $\cl \rho: G_K \ra \GL_n(\F)$, where $K$ is a local field of residue characteristic $l$, $\F = \mO_L/\lambda$ a finite field of characteristic $p \ne l$. We further assume that $L$ is large enough in the sense of Corollary \ref{cor:inertia_component} -- this ensures that the Galois inertia type is constant on each connected component of $\spec \textup{R}_{\cl \rho}^{\square}$ (or its closed subscheme $\spec \textup{R}_{\cl \rho}^{\square, \chi}$). 

\subsubsection{The cycle maps} We first define the cycle maps, in terms of the generic multiplicity map defined in Definition \ref{def:multiplicity_map} and the following definition.   

\bd  \label{def:component_maps}
Let $\cl \rho$ and $\chi$ be as before, and fix $L$ large enough. We define the component maps $c_{\cl \rho}$ and $c_{\cl \rho}^{\chi}$
 $$ c_{\cl \rho}: \bigoplus_{\tau \in B^t (\GL_n)} \Z \ra \mZ(\textup{R}_{\cl \rho}^{\square}), \qquad    c_{\cl \rho}^{\chi}: \bigoplus_{\tau \in B^t(\GL_n)}\Z \ra  \mZ(\textup{R}_{\cl \rho}^{\square, \chi}) $$ 
as follows. We identify $\tau \in B^t(\GL_n) = \Igaln$ by Lemma \ref{lemma:types_and_centers}, and define 
$$c_{\cl \rho}: 1_\tau \mapsto Z(\textup{R}_{\cl \rho}^{\square, \tau}), \qquad c_{\cl \rho}^{\chi}: 1_\tau \mapsto Z(\textup{R}_{\cl \rho}^{\square, \tau, \chi}).$$ Here we follow the notation from Subsection \ref{sss:cycles}, $Z(\textup{R}_{\cl \rho}^{\square, \tau}) = [\textup{R}_{\cl \rho}^{\square, \tau}]_{n^2+1} \in \mZ(\FDR)$, and similarly for $Z(\textup{R}_{\cl \rho}^{\square, \tau, \chi})$.  
\ed

\bd[Cycle maps] Retain notations from the definition above. We define the cycle map  
$$\textup{cyc}_L: \textup{R}_{L} (\mathtt{K}) \longrightarrow \mZ(\FDR), \quad \textup{cyc}_L^\chi: \textup{R}_{L} (\mathtt{K}) \longrightarrow \mZ(\FDRx)$$ to be the composition 
\bi
\item $\textup{cyc}_L := c_{\cl \rho} \circ m_L$.
\item  $\textup{cyc}_L^{\chi} := c_{\cl \rho}^\chi \circ m_L$.
\ei
where $m_L$ is the generic multiplicity map and $c_{\cl \rho}, c_{\cl \rho}^\chi$ are defined above. 
\ed 

\br \label{remark:relation_with_fixing_determinant} \indent 

\be
\item Let $\mathscr{C}$ be the set of inertia types appearing in the deformation ring $\textup{R}_{\det \cl \rho} = \textup{R}_{\det \cl \rho}^{\square}$ of the character $\det \cl \rho$, namely that $\textup{R}_{\det \cl \rho} [1/p]$ is the disjoint union $ \underset{c \in \mathscr{C}}{\sqcup} \textup{R}_{\det \cl \rho}^{c} [1/p]$, where each  $\textup{R}_{\det \cl \rho}^{c}$ parametrizes liftings of $\det \cl \rho$ with inertia type $c$. For each $c \in \mathscr{C}$, choose a lifting $\chi$ of $\det \cl \rho$ of with inertia type $c$. These form a set $[\mathscr{C}]$, which is the set of $\cl L$-points of $\textup{R}_{\det \cl \rho}$, one for each $c \in \mathscr{C}$. If $\chi$ and $\chi'$ lie on the same component $\textup{R}_{\det \cl \rho}^{c}$, then $m(\sigma, \tau, \chi) = m (\sigma, \tau, \chi')$ by Remark \ref{remark:well_defined_chi_multiplicity}. Unwinding definitions, we conclude that 
$$\textup{cyc} (\sigma) = \sum_{\chi \in [\mathscr{C}]} \iota^{\chi} (\textup{cyc}^{\chi} (\sigma)),$$ 
where $\iota^{\chi}$ is the map 
$$\iota^\chi: \mZ(\FDRx) \ra \mZ(\FDR)$$ 
which sends the cycle $Z(\textup{R}_{\cl \rho}^{\square, \tau, \chi})$ to the cycle $Z(\textup{R}_{\cl \rho}^{\square, \tau})$ for each inertia type $\tau$, and everything else to $0$. 
\item $\textup{cyc}$ (resp. $\textup{cyc}^{\chi}$) is compatible with base change of the coefficient field. To be more precise, suppose $L'/L$, then $\textup{R}_{\cl \rho, \mO_L'}^{\square} \cong \textup{R}_{\cl \rho, \mO_L}^{\square}\!\otimes_{\mO_L}\!\!\mO_{L'} $, and we have the natural map $\mZ(\textup{R}_{\cl \rho, \mO_L}^{\square}) \ra \mZ(\textup{R}_{\cl \rho, \mO_{L'}}^{\square}) $ induced by base change. The following diagram commutes 
\[
\begin{tikzcd}
\textup{R}_{L} (\mathtt{K})  \arrow{r}{  \textup{cyc}_{L} } \arrow{d}[swap]{ \otimes_{L}L' }  & \mZ(\textup{R}_{\cl \rho, \mO_L}^{\square}) \arrow{d}{\otimes_{\mO_L} \mO_{L'}} \\
\textup{R}_{L'} (\mathtt{K})   \arrow{r}{  \textup{cyc}_{L'} }   & \mZ(\textup{R}_{\cl \rho, \mO_{L'}}^{\square}) 
\end{tikzcd}
\] The same remark applies to $\textup{cyc}^\chi$. 
\ee
\er

\noindent Finally, the cycle maps and the first remark can be summarized in the following diagram (where the conjectural arrow suggests the existence of tempered types). 
\[
\begin{tikzcd}
  \textup{R}_{L} (\mathtt{K}) \arrow{rr}{m_L}    && \underset{s \in B^t(\GL_n)}{\bigoplus} \Z    \arrow[dashed, bend right]{ll}[swap]{\sigma_{\textup{t}}} \arrow{r}{c_{\cl \rho}} \arrow{rd}{c_{\cl \rho}^{\chi}}    &\mZ (\FDR)     \\
& & & \mZ(\textup{R}_{\cl \rho}^{\square, \chi}) \arrow[hook]{u}[swap]{\iota^{\chi}}  
\end{tikzcd}
\]
 
\br[Compare to \cite{S}] Let $\delta:  \textup{R}_{L} (\mathtt{K}) \ra \textup{R}_{L} (\mathtt{K})$ be the involution induced by taking dual representations, and define 
$\textup{cyc}_{L, \delta} := \textup{cyc}_{L} \circ \delta$, as well as $\textup{cyc}_{L, \delta}^\chi := \textup{cyc}_{L}^\chi \circ \delta$.  
The twisted version of cycle maps are precisely what appears in \cite{S}. More precisely, we have 
\be 
\item 
 $\textup{cyc}_{L, \delta}: \textup{R}_{L} (\GL_n(\mO_{K})) \ra \mZ (\FDR)$ is given   by 
$$ \textup{cyc}_{L, \delta}(\sigma) := \sum_{\textup{ all inertial types}} m (\sigma^{\vee}, \tau) \cdot Z(\textup{R}_{\cl \rho}^{\square, \tau}).$$
\item Similarly,  $  \textup{cyc}^{\chi}_{L, \delta}: \textup{R}_{L} (\GL_n(\mO_{K})) \ra \mZ (\FDRx)$ is given by  
$$ \textup{cyc}^{\chi}_{L, \delta}(\sigma) := \sum_{\textup{ all inertial types}} m (\sigma^{\vee}, \tau, \chi) \cdot Z(\textup{R}_{\cl \rho}^{\square, \tau, \chi})$$
\ee
\er

\subsubsection{The Breuil--M\'ezard conjecture} We now have everything in place to (re)state the Breuil--M\'ezard conjecture (namely Conjecture \ref{conj:intro} in the introduction). 

\begin{conjecture} \label{conj:BM} Let $K$ be a non-archimedean local field of residue characteristic $l$ as before, $\F$ a finite field of characteristic $p \ne l$. Let $\cl \rho: G_K \ra \GL_n(\F)$ be a continuous representation. 
\be 
\item There exists a mod $\lambda$ cycle map 
$$ \cl{\textup{cyc}} : \textup{R}_{\F} (\mathtt{K}) \longrightarrow \mZ (\FDR/\lambda) $$ making the following diagram commute: 
\[
\begin{tikzcd}
  \textup{R}_{L} (\mathtt{K}) \arrow{rr}{m_L} \arrow[d, swap, "\textup{red}"]  & & \underset{s \in B}{\bigoplus} \Z    \arrow{r}{c_{\cl \rho}}    &\mZ(\FDR)  \arrow[d, "\textup{red}"] \\
  \textup{R}_{\F} (\mathtt{K}) \arrow[rrr, dashed, "\cl{\textup{cyc}}"] & & & \mZ(\FDR/\lambda) 
\end{tikzcd}
\]
\item Moreover, for each character $\chi$ lifting $\det \cl \rho$, there exists a mod $\lambda$ cycle map with respect to $\chi$: 
$${\cl{\textup{cyc}}}^{\chi}: \textup{R}_{\F} (\mathtt{K}) \longrightarrow \mZ (\FDRx/\lambda) $$ fitting into a similar commutative diagram
\[
\begin{tikzcd}
  \textup{R}_{L} (\mathtt{K}) \arrow{rr}{m_L} \arrow[d, swap, "\textup{red}"] \arrow[drrr, swap, "\textup{cyc}^{\chi}_L"]  & & \underset{s \in B}{\bigoplus} \Z    \arrow{rd}{c_{\cl \rho}^{\chi}}     \\
  \textup{R}_{\F} (\mathtt{K}) \arrow[drrr, dashed, "\cl{\textup{cyc}}^\chi"] & & & \mZ(\textup{R}_{\cl \rho}^{\square, \chi}) \arrow[d, "\textup{red}"] \\ 
&  & & \mZ(\FDRx/\lambda)
\end{tikzcd}
\]

\ee
\end{conjecture}

\br  \indent 
\be
\item As remarked in the introduction, the surjectivity of the reduction map $\textup{R}_{L} (\mathtt{K}) \ra \textup{R}_{\F} (\mathtt{K})$ implies that, the map $\cl{\textup{cyc}}$ (resp. $\cl{\textup{cyc}}^\chi$) is necessarily unique.  
\item The conjecture is equivalent to the following naive-looking form: it simply requires that $\ker (\textup{red}) \subset \ker (\textup{red} \circ \textup{cyc}_L)$ (resp. $\ker (\textup{red}) \subset \ker (\textup{red} \circ \textup{cyc}^{\chi}_L)$).  
\ee
\er

\br  \label{remark:fix_determinant_imply_no_fix}  \indent 
 By Remark \ref{remark:relation_with_fixing_determinant}, part (2) of the conjecture implies part (1). \textit{a priori} it is unclear how to deduce part (2) from part (1). Therefore later we will prove the slightly stronger version, namely with fixed determinant. The main result of \cite{S} is stated in the form without specifying $\chi$. 
\er

\br[{Invariance under base change}] \indent   \label{sss:invariance_base_change}
The diagrams in the conjecture above depend on $L$. For extension of the coefficient field $L'/L$, there are natural maps between the two diagrams for $L$ and $L'$ making everything commute. However, it suffices to prove the conjecture after base change, since the natural map 
$$ \mZ (\FDRx /\lambda) = \mZ (\FDRx \otimes \F) \longrightarrow  \mZ (\FDRx \otimes \F')$$ is injective, where $\F'$ is the residue field of $L'$. From now on we will make no further mention of the coefficient field $L$, and will simply assume that it is taken to be large enough.  
\er

\subsubsection{Relation with the classical form} \label{sss:relation_with_classical_form}  

As mentioned in the introduction, the bridge between the cycle version of the Breuil--M\'ezard Conjecture (Conjecture \ref{conj:BM}) and its classical form is provided by tempered types, more precisely, the conjectural existence of a typical section $\sigma_t$.  Under the existence of $\sigma_t$ (namely Conjecture \ref{conj:2}), Conjecture \ref{conj:BM} implies the following ``classical form'': 

\bc[of Conjecture \ref{conj:BM}] \label{cor:classical_form}
Let $\cl \rho$ be a continuous mod $\lambda$ representation, and fix $\chi: G_K \ra \mO_L^\times$ lifting $\det \cl \rho$. There exists an (infinite) collection of integers $\{\mu_{\cl \sigma} (\cl \rho)\}_{\cl \sigma}$ indexed by the set of isomorphism classes of irreducible $\F$-representations of $\mathtt{K} = \GL_n (\mO_K)$, such that for each Galois inertia type $\tau$, we have 
$$e(\textup{R}_{\cl \rho}^{\square, \tau, \chi} /\lambda) = \sum_{\cl \sigma} n_{\cl \sigma} (\tau) \cdot \mu_{\cl \sigma } (\cl \rho)$$ 
where $e$ denotes the Hilbert--Samuel multiplicity, and $n_{\cl \sigma} (\tau)$ denotes the Jordan--Holder multiplicity of $\cl \sigma$ in $\cl{L_{\tau}}^{\textup{s.s.}}$, which is the semi-simplification of an invariant lattice $L_{\tau} \subset \sigma_t(\tau)$. Note that $\mu_{\cl \sigma} (\cl \rho)$ only depends on $\cl \rho$. 
\ec

\bproof 
This is clear from definitions, since $\textup{cyc}^{\chi} (\sigma (\tau)) = Z(\textup{R}_{\cl \rho}^{\square, \tau, \chi})$. 
\eproof 

\br 
However, the classical form does not \textit{a priori} imply the cycle version of the conjecture, this is because in general there are nontrivial elements in $\mZ(\FDRx/\lambda)$ with Hilbert--Samuel multiplicity equal to $0$. Therefore, the cycle version is more refined than its classical form of the Breuil--M\'ezard conjecture in terms of Hilbert--Samuel multiplicities. 
\er 

\br 
One advantage the cycle version of the conjecture enjoys is the following simplification in its proof using global arguments (already observed in \cite{S}) compared to the proof in \cite{GK} (there they consider for potentially Barsotti--Tate representations in the $l = p$ case). Namely, the linear algebra arguments reducing the conjecture for a product of representations to an individual term are completely avoided. However, the linear algebra lemma (Lemma 3.5.2 in \cite{GK}) implies that the integers $\mu_{\cl \sigma} (\cl \rho)$ in their context, if exist, are unique. Our conjecture predicts a necessarily unique map $\cl{\textup{cyc}}$, it however does not directly imply the uniqueness of $\mu_{\cl \sigma} (\cl \rho)$. 
\er

\subsubsection{Inverting the cycle map}  The formulate in Corollary \ref{cor:classical_form} above can be regarded as the consequence of an ``optimal'' way to invert the map $\textup{cyc}$ (and $\textup{cyc}^{\chi}$), ``optimal'' in the sense that one only needs to attache to $\tau$ a single irreducible representation $\sigma_t (\tau)$ to compute $e(\textup{R}_{\cl \rho}^{\square, \tau, \chi})$. Less optimal ways to invert $\textup{cyc}$ and $\textup{cyc}^{\chi}$ are provided by virtual sections. For example, by Corollary \ref{cor:SZ_section}, we have a way to attach each $\tau$ a virtual representation $\sigma_v(\tau)$. Therefore, we arrive at the following corollary of the Conjecture (this is already observed in \cite{S} as Corollary 4.9)

\bc \indent 

\be
\item For each Galois inertia type $\tau \in \Igaln$, 
$$\textup{red} ( Z(\textup{R}_{\cl \rho}^{\square, \tau, \chi})) = \cl{\textup{cyc}}^{\chi} (\textup{red} (\sigma_v(\tau)) );$$
\item There exists integers $\{\mu_{\cl \sigma} (\cl \rho)\}_{\cl \sigma}$ indexed by isomorphism classes of irreducible $\F$-representations of $\mathtt{K}$, such that for each inertia type $\tau$, 
$$e(\textup{R}_{\cl \rho}^{\square, \tau, \chi} /\lambda) = \sum_{\cl \sigma} n_{\cl \sigma} (\sigma_v (\tau)) \cdot \mu_{\cl \sigma } (\cl \rho)$$ 
where $n_{\cl \sigma} (\sigma_v(\tau))$ is determined as follows: if $\sigma_v(\tau)$ is represented by $ \varsigma_1 - \varsigma_2$ with  each $\varsigma_i$ an honest representation of $\mathtt{K}$, then $n_{\cl \sigma} (\sigma_v(\tau)) = n_{\sigma} (\varsigma_1) - n_{\sigma}(\varsigma_2)$ where $n_{\sigma}(\varsigma_i)$ is the Jordan--Holder multiplicity of $\cl \sigma$ in $\cl{L_{\varsigma_i}}^{s.s.}$ as before.  
\ee
\ec

\subsubsection{Determination of $\cl{\textup{cyc}}$} 

In the classical Breuil--M\'ezard conjecture for potentially semistable deformation ring (of $2$-dimensional representations) with fixed inertia type $\tau$ and Hodge--Tate weights $\mathbf{v}$, the multiplicities $\mu_{\cl \sigma} (\cl \rho)$ can be determined thanks to the classification of the finitely many irreducible mod $p$ representations of $\GL_2 (\mO_K)$ and can be expressed in terms of Serre weights. 

 In the $l \ne p$ case, however, the determination of $\cl{\textup{cyc}}$ (or even $\mu_{\cl \sigma}$) seem to be a different story, which is closely related to certain mod $p$ (where $p \ne l$) inertial Langlands correspondence. We hope to return to this aspect in a future project. 

\subsection{The main theorem} Now we state the main theorem of the article, which establishes most of Conjecture \ref{conj:BM} for local function fields. 

\subsubsection{The Breuil--M\'ezard conjecture for function fields} 

\bt \label{thm:Main} Let $K/\F_l (\!(t)\!)$ be a finite extension. Let $\cl \rho: G_{K} \ra \GL_n (\F)$ be continuous representation and $\mathtt{K} = \GL_n (\mO_K)$ as before. Suppose that $p \nmid n(n-1)$, and that $p \nmid (l - 1)$. Then both parts of the Conjecture \ref{conj:BM} holds. 
\et
By Remark  \ref{remark:fix_determinant_imply_no_fix}, an equivalent form of the theorem is that for every $\chi$, we have 
$$\ker (\textup{red}) \subset \ker (\textup{red} \circ \textup{cyc}^{\chi}).$$ 
Note that we have suppressed the subscript of $L$ since it is irrelevant in the theorem. 

\subsubsection{Outline of the proof} 

The main theorem follows directly from the following proposition. 

\bp \label{prop:main_prop} 
There exists an exact functor 
 $$M_\infty: \textup{Rep}_{\mO}^{\textup{f.l.}} (\mathtt{K}^s) \longrightarrow \textup{Mod}_{/ R_0 \widehat{\otimes} (\FDRx)^{\widehat{\otimes} s} }  $$
where $s$ is a positive integer and $R_0$ is a power series ring over $\mO$ such that 
\be
\item If $\sigma = \otimes \sigma_i \in \textup{Rep}_{\mO}^{\textup{f.l.}} (\mathtt{K}^s)$ is $\lambda$-torsion free, then 
$$ Z(M_\infty (\sigma)) = n! \cdot [\spec R_0] \times \textup{cyc}^{\chi} (\sigma_1^\vee ) \times \cdots \times  \textup{cyc}^{\chi} (\sigma_s^\vee)$$ as a cycle in $\mZ (R_0 \widehat{\otimes} (\FDRx)^{\widehat{\otimes} s})$. 
\item $M_\infty$ is compatible with reduction mod $\lambda$, namely for any $\sigma$, 
$$M_{\infty} (\sigma) \otimes_{\mO} \F = M_\infty (\sigma \otimes_{\mO} \F)$$
\item If $\sigma$ is $\lambda$-torsion free, then $M_\infty (\sigma) $ is also $\lambda$-torsion free. In particular, Lemma \ref{lemma:reduction_cycle_mod_p} applies to $M_\infty (\sigma)$. 
\ee
\ep

\br Note that in the statement of property (1), we have used the identification 
$$ \mZ (R_0 \widehat{\otimes} (\FDRx)^{\widehat{\otimes} s}) = \mZ (R_0) \otimes \mZ(\FDRx)^{\otimes s},$$
which follows from Lemma \ref{lemma:tensor_product_cycle} and Proposition \ref{lemma:deformation_character_flat}.  
\er 

\bproof[Proof of Proposition \ref{prop:main_prop} $\so$ Theorem \ref{thm:Main}] This is similar to the proof of Theorem 4.6 in \cite{S}, which we include for completeness. Note that the map $\mF \mapsto Z(\mF)$ of taking cycles (counted with multiplicity) is additive on short exact sequences, hence by exactness of $M_\infty$ we get well defined maps 
$$ \textup{R}_{L} (\mathtt{K}^s) \xrightarrow{Z \circ M_\infty \circ \delta} \mZ (\FDRx)^{\otimes s} \quad \textup{and }\: \:  \textup{R}_{\F} (\mathtt{K}^s) \xrightarrow{Z \circ M_\infty \circ \delta} \mZ (\FDRx /\lambda)^{\otimes s},$$ 
where $\delta$ denotes taking dual representations in $\textup{R}_{L} (\mathtt{K}^s)$ (resp. in $\textup{R}_{\F} (\mathtt{K}^s)$).  Now by properties (2) and (3) of $M_\infty$, we know that for any $\sigma \in \textup{R}_{L} (\mathtt{K}^s)$, 
$$ Z\big(M_{\infty} (\textup{red} (\sigma))^\vee \big)  = Z(M_{\infty} (\textup{red} (\sigma^\vee)))  =  Z(M_\infty(\sigma^\vee) \otimes \F) =  \textup{red}(Z(M_\infty(\sigma^\vee)))$$ where the last equality follows from Lemma \ref{lemma:reduction_cycle_mod_p} and the fact that $M_\infty(\sigma^\vee)$ is $\lambda$-torsion free. In other words, we have the following commutative diagram 
\[
\begin{tikzcd} [column sep=4em,row sep=2em]
\textup{R}_{L} (\mathtt{K}^s) \arrow[r, "Z \circ M_\infty \circ \delta"]  \arrow[d, "\textup{red}"] & \mZ (\FDRx)^{\otimes s} \arrow[d, "\textup{red}"] \\
 \textup{R}_{\F} (\mathtt{K}^s)    \arrow[r, "Z \circ M_\infty \circ \delta"] & \mZ(\FDRx/\lambda)^{\otimes s}
\end{tikzcd}
\]
Then it is clear that $\ker (\textup{red}) \subset \ker (\textup{red} \circ \textup{cyc}^{\chi})$. 
\eproof

The construction of the functor $M_\infty$ will occupy the next two sections, which we briefly outline in the rest of this section (also see the introduction). The proof that $M_\infty$ has the desired properties is given in Subsection \ref{sss:proof_main_prop}. 

\subsubsection{Global realization} \label{sss:global_realization_function_field} 
Let $\cl \rho: G_K \ra \GL_n(\F)$ be the given continuous representation, we will first realize $\cl \rho$ as the local component of the reduction $\cl r$ of a global Galois representation 
$$r: G_{F} \ra \GL_n (\mO_{\cl L}),$$ 
namely that $\cl \rho = \cl r |_{F_v}$ at a certain finite set of places $v \in S_1$ of $F$. For this we first apply a theorem of Moret-Bailly to produce the mod $\lambda$ representation $\cl r$, and then lift $\cl r$ to $r$ using a result of Gaitsgory together with the computation of the Krull dimension of certain global Galois deformation rings. 

In addition, $r$ is controlled carefully at other finite sets of places $S_2$ and $S_D$ disjoint from $S_1$, so that $\cl r$ is trivial at $S_2$, and is the special representation (cf. Definition \ref{def:special_reps}) at $S_D$ (in particular indecomposable).  This allows us to apply the global Langlands and Jacquet-Langlands correspondence to obtain an automorphic representation $\pi$ of some central division algebra $D$ over $F$, which satisfies $D_v \cong \GL_n(K)$ and gives rise to the Galois representation $r = r(\pi)$. 

\subsubsection{Patching}  We then apply the  Taylor--Wiles--Kisin patching method to construct the functor $M_\infty$. For each finitely generated $\mO$-module $M_{\sigma}$ with an action of $\mathtt{K}$, we construct a Taylor--Wiles systems $Q_N$ of auxiliary primes, and define the spaces of automorphic forms $S_N(\sigma)$ such that the action of $\mathtt{K} \subset \GL_n(K) \subset D_v$ at $v$ is prescribed to be induced from the  action of $\mathtt{K}$ on $M_{\sigma}$. In the infinite level, we patch to get the desired module $M_{\infty} (\sigma)$. It then remains to show that the support of $M_\infty (\sigma)$ satisfies property (1) in Proposition \ref{prop:main_prop}, which follows from the fact that ``all points on the deformation rings are automorphic'', which in the function field case is a consequence of L. Lafforgue's theorem on global Langlands correspondence for $\GL_n$.

\section{Global realization of local Galois representations} \label{sec:Global}

The goal of this section is to prove Proposition  \ref{thm:global_realization},  namely to  realize $\cl \rho$ as the mod $\lambda$ reduction of the local restriction of a global Galois representation $r$ which is automorphic. 

\subsection{Galois deformations} \label{ss:Galois_Def} We first recall and modify certain calculations of tangent spaces of global deformation rings, which we need both in proving Proposition  \ref{thm:global_realization} and in setting up the patching argument.   For this subsection, we let $K$ be any local function field over $\F_l (\!(t)\!)$ and $\cl \rho: G_K \ra \GL_n (\F)$ an arbitrary continuous Galois representation, while keeping the rest of our notations as well as assumptions on $\mO= \mO_L$ and $\cl r: G_F \ra \GL_n (\F)$. We hope this will not cause confusion. 

\subsubsection{The deformation rings}  
Let $F$ be a global field over $\F_l (t)$. We fix a finite set of primes $T$ of $F$, and a finite set of auxiliary primes $Q$ disjoint from $T$. Let 
$$\cl r: G_F \ra \GL_n (\F)$$ be a continuous representation unramified outside $T$. Let $S = S_Q = T \sqcup Q$. We also fix a local deformation problem $D_v$ (in the sense of Definition 2.2.2 of \cite{CHT}) at each $v \in S $, which corresponds to an ideal $I_v \subset \FDRx$ that is invariant under $\ker (\GL_n(\FDRx) \ra \GL_n(\F))$. 
We will consider Galois deformations of the $G_{F, S}$ representation $\cl r$ with local conditions $D_v$ at each $v \in S$ and local frames at $v \in T$.  More precisely, similar to Definition \ref{def:deformation_functor}, let us fix a character $\psi: G_{F, S} \ra \mO^\times$ lifting $\det \cl r$ and define the deformation functor $D_{\cl r, S}^{\square_T, \psi} : \mC_{\mO} \longrightarrow \textup{Sets}$ by
$$D_{\cl r, S}^{\square_T, \psi} (A) = \{ (r, \{ \alpha_v\}_{v \in T}) \} / \simeq $$ 
where $r \in D_{\cl r}^{\square, \psi} (A)$ is a framed deformation of $\cl r$ with determinant $\psi$, such that $r|_{G_{F_v}} \in D_v$ for all $v \in S$; and $\{\alpha_v\}$ is a collection of elements in $\ker (\GL_n (A) \ra \GL_n (\F))$ for each $v \in T$. The equivalence relation $\simeq$ is given by conjugation by  $ \ker (\GL_n (A) \ra \GL_n (\F))$, namely
$$(r, \{\alpha_v\}_{v \in T}) \simeq (a r a^{-1}, \{a \cdot \alpha_v\}_{v \in T})$$ for $a \in \ker (\GL_n (A) \ra \GL_n (\F))$. Suppose that $\cl r$ is absolutely irreducible, then 
$D_{\cl r, S}^{\square_T, \psi}$ is pro-representable by a complete noetherian ring $R_{\cl r, S}^{\square_T, \psi}$,  with the universal map  
$$r^{\square_T, \psi}_{S}: G_{F, S} \ra  \GL_n (R_{\cl r, S}^{\square_T, \psi}). $$  
In the case when $T = \O$, we denote the global deformation ring without frames by $R^{\textup{univ}, \psi}_{\cl r, S}$. 
Finally, we define 
$$R^{\loc} := \underset{v \in T}{\widehat \otimes} (R_{\cl \rho_v}^{\square, \chi_v}/I_v).$$
By the universal property of each $R_{\cl \rho_v}^{\square, \chi_v}/I_{v}$, we obtain a canonical map  
$$R^{\loc} \longrightarrow R_{\cl r, S}^{\square_T, \psi}.$$
A crucial step in the Taylor--Wiles--Kisin patching method (\cite{K1}, \cite{CHT}) is to compute the relative dimension of $R_{\cl r, S}^{\square_T, \psi}$ over $R^{\loc}$ (at the closed point).  Next we carry out this standard computation in the function field setup. 

\subsubsection{The subspace $\sq I_v \subset Z^1(G_K, \ad^0 \cl \rho)$} 
 
In this discussion let $K = F_v$ for some place $v$ of $F$. We write $\F [\epsilon]$ for $\F [\epsilon]/\epsilon^2$. Let $\chi: G_K \ra \mO^\times $ be a character such that $\chi \equiv \det \cl \rho  \mod \lambda $.  Let $\ad^0 \cl \rho := \{ \Sigma \in \ad \cl \rho : \tr \Sigma = 0\} \subset \ad \cl \rho$, where $G_K$ acts by conjugation. We denote by $B^1 (G_K, \ad^0 \cl \rho)$ (resp. $Z^1 (G_K, \ad^0 \cl \rho)$) the space of $1$-coboundaries (resp. $1$-cocycles). Recall that as $\F$-vector spaces we have 
$$ \Hom_{\F} (\fm_{\FDRx}/ (\fm^2_{\FDRx}, \lambda), \F) \isom D_{\cl \rho}^{\square, \chi} (\F[\epsilon]) \isom Z^1 (G_K, \ad^0 \cl \rho).$$ 
To compute the dimension of $\FDRx/I_v $,  first consider 
the image of $I_v$ given by  
$$\sq I_v := \im (I_v \ra \fm_{\FDRx}/(\fm_{\FDRx}^2, \lambda)),$$ and then define $\sq L_v$ as the annihilator of $\sq I_v$, 
$$\sq L_v := \{ f : \fm_{\FDRx}/(\fm_{\FDRx}^2, \lambda) \ra \F  \: | \:  f (\sq I_v) = 0 \}. $$
\br \label{remark:def_L_v}  From the definition we know that
\bi
\item $\sq L_v$ is the preimage of $L_v := \sq L_v /B^1(G_K, \ad^0 \cl \rho) \subset H^1(G_K, \ad^0 \cl \rho)$.
\item $Z^1(G_K, \ad^0 \cl \rho) / \sq L_v \isom  H^1 (G_K, \ad^0 \cl \rho)/ L_v $.
\ei
\er 
Finally let $L_v^\perp \subset H^1(G_K, \ad^0 \cl \rho (1))$ be the annihilator of $L_v$, under the pairing induced from the trace pairing on $ \ad^0 \cl \rho \times \ad^0 \cl \rho (1) \ra \F (1)$. 
 
\subsubsection{The relative dimension}  

The computation of the dimensional of the space 
$$ \Hom_{\F} (\fm_{R_{S}^{\square_T, \psi}} / (\fm^2_{R_{S}^{\square_T, \psi}}, \fm_{R^{\loc}}, \lambda), \F) $$ is identical to the number field case. For each $v\in S = S_Q$, write $G_v = G_{F_v}$.  The $\F$-vector space
$ \Hom_{\F} (\fm_{R_{\cl r, S}^{\square_T, \psi}} / (\fm^2_{R_{\cl r, S}^{\square_T, \psi}}, \fm_{R^{\loc}}, \lambda), \F) $
is canonically isomorphic to first cohomology $H^1$ of the $3$-term complex 
       \begin{multline*} 
       \ad \cl r \xrightarrow{ \partial \: \oplus \: \underset{v \in T}{\oplus} \textup{res}_v }  Z^1 (G_{F, S}, \ad^0 \cl r) \oplus \underset{v \in T}{\oplus} \ad \cl r  \\ 
 \xrightarrow{ \big( \underset{v \in S }{\oplus}  \textup{ res}_v \big) \oplus  \underset{v \in T}{\oplus} \textup{-}\partial  \:  }     \underset{v \in T}{\oplus} Z^1 (G_v, \ad^0 \cl r)  \oplus  \underset{v \in S\minus T}{\oplus} { Z^1 (G_v, \ad^0 \cl r)}/{ \sq L_{v} } 
        \end{multline*} 
On the other hand, $H^1$ of the complex above is isomorphic to $H^1_{S, T} (G_{F, S }, \ad^0 \cl r) := H^1 (C_{S, T}^\bullet)$, where $C_{S,T}^\bullet$ is a complex defined as follows. First, define $C_{\loc}^{\bullet}$ by 
\bi
\item $C_{\loc}^0 : =  \underset{v \in T}{\oplus} C^0 (G_v, \ad \cl r) $
\item $C_{\loc}^1 : =  \underset{v \in T}{\oplus} C^1 (G_v, \ad^0 \cl r) \oplus  \underset{v \in Q}{\oplus} C^1 (G_v, \ad^0 \cl r) /\sq L_{v} $
\item $C^i_{\loc} :=  \underset{v \in S}{\oplus} C^i (G_v, \ad^0 \cl r) $ \quad for $i  \ge 2$. 
\ei
Next define $C_0^i$ by 
\bi
\item$C_0^0 : = C^0 (G, \ad \cl r) = \ad \cl \rho, \qquad C_0^i := C^i (G, \ad^0 \cl r) \textup{ for } i \ge 1$ 
\ei
Finally, define $(C_{S,T}^{\bullet}, d_{S, T}^\bullet)$ to be the complex
$$C_{S, T}^i := C_0^i \oplus C_{\loc}^{i-1}, \qquad d_{S, T}^i: (\phi, \{\psi_v\} ) \mapsto (\partial \phi, \{ \phi |_{G_v} - \partial \psi_v \}).$$
 We proceed to compute $H_{S, T}^1$ as in \cite{CHT} and \cite{K1}, using the long exact sequences of cohomology associated to the exact sequence of complexes
$$ 0 \ra C^{i-1}_{\loc} \ra C_{S, T}^i \ra C_0^i \ra 0.$$
For us $F/\F_l(\!(t)\!)$ is a function field and $\F$ has characteristic $p \ne l$, so we have the following standard facts.

\bp \label{prop:results_from_Milne} Let $V$ be a finite dimensional $\F$-vector space. 
\be
\item Local and global vanishing. Let $V$ be a $G$ representation where $G = G_K$ or $G_{F, S}$, then for $ i \ge 3$, 
$$ H^i (G, V) = 0. $$
\item  Local Tate duality. Let $V$ be a $G_{K}$ representation, then $$H^r (G, V^{\vee} (1)) \cong H^{2 - r} (G, V)^*$$ where ${}^\vee$ denotes the dual representation, and ${}^*$ denotes the Pontryagin dual. 
\item Euler characteristic. Let $V$ be a $G$ representation where $G =G_{K}$ or $G_{F, S}$, then in both cases 
$$ \chi (G, V) := \sum_{i = 0, 1, 2} (-1)^{i} \dim_{\F} H^i(G, V) = 0$$ 
\ee
\ep

\bproof  Note that $\F =\mO_L/\lambda$ has characteristic $p \ne l$. For (1), the local case follows from Corollary 2.3 (local Tate duality) in \cite{ADT}, the global case is Theorem 4.10 (c) in \textit{ibid}.
For (2). The case $G = G_K$ is Theorem 2.8 in \textit{ibid}, the global case  follows from Theorem 5.1 (and the remark that follows) in \textit{ibid}.
\eproof 
Denote by $\chi_{\loc}, \chi_{S, T}, \chi_0$ respectively the Euler characteristic of the complexes $C^{\bullet}_{\loc},C^{\bullet}_{S,T},C^{\bullet}_{0}$, recall that $p \nmid n$ by assumption, therefore 
\begin{align*}
\chi_{S, T} =&  \chi_0 - \chi_{\loc} \\
= & 1 - \# T + \sum_{v \in Q} (\dim_{\F} H^0 (G_v,  \ad^0 \cl r) - \dim_{\F} L_{v}) + \chi (G_{F, S}, \ad^0 \cl r) - \sum_{v \in S} \chi (G_v, \ad^0 \cl r) \\
= & 1 - \# T + \sum_{v \in Q} (\dim_{\F} H^0 (G_v,  \ad^0 \cl r) - \dim_{\F} L_{v})
\end{align*}
Finally, using Poitou-Tate sequence for function fields (and following the same proof of Lemma 2.3.4 in \cite{CHT}) we have 
\bc \label{cor:relative_dim}
Suppose that $\cl r$ is absolutely irreducible and $p \nmid nl$, let $t = \# T$ be the size of $T$, then 
 \begin{multline*} 
\dim_{\F} \Hom (\fm_{R_{\cl r, S}^{\square_T, \psi}} / (\fm^2_{R_{\cl r, S}^{\square_T, \psi}}, \fm_{R^{\loc}}, \lambda), \F)
 \\ =  t - 1 + \sum_{v \in Q} \Big(\dim_{\F} L_{v} - \dim_{\F} H^0 (G_v, \ad^0 \cl r) \Big) +  \dim_{\F} H^1_{S, T} (1)  
 \end{multline*} 
where $H^1_{S, T} (1):=  \ker \Big(    H^1 (G_{F, S}, (\ad^0 \cl r) (1)) \xrightarrow{\oplus \textup{res}} \underset{v \in Q}{\oplus}   H^1(G_{v}, \ad^0 \cl r \: (1))/(L_{v})^{\perp}   \Big). $
\ec

\bc \label{cor:relative_dim_bound}
Suppose that $\cl r$ is absolutely irreducible and $p \nmid nl$, then 
\be
\item The deformation ring $R_{\cl r, S}^{\square_T, \psi}$ has Krull dimension at least 
$$\sum_{v \in Q} \Big(\dim_{\F} L_{v} - \dim_{\F} H^0 (G_v, \ad^0 \cl r) \Big) + \sum_{v \in T} \dim \textup{R}_{\cl r|_v}^{\square, \psi_v} /I_v $$
\item Consequently, the Krull dimension of $R^{\textup{univ}, \psi}_{\cl r, S}$ is least 
$$ 1+ \sum_{v \in S} (\dim \textup{R}_{\cl r|_v}^{\square, \psi_v} /I_v  - n^2) $$
\ee
\ec 

\bproof 
(1) is immediate from Corollary \ref{cor:relative_dim}. (2) follows from (1), because by definition $\dim_{\F} L_v - \dim_\F H^0(G_v, \ad^0 \cl r) = \dim_\F (\sq L_v) - (n^2 - 1) \ge 0$.  
\eproof 
 
\subsubsection{A finiteness result}  In order to build a global automorphic representation $r$ such that the local components of $\cl r$ realize $\cl \rho$, we will need to exhibit a $\cl L$-point of certain global Galois deformation rings. Following the strategy of Khare-Wittenberger, we create such a point by showing that the relevant deformation ring is a finite module over $\mO$ and has Krull dimension at least $1$. 

The finiteness result we need follows from a theorem of Gaitsgory and is closely related to a conjecture of de Jong.  To explain the setup let $k/\F_l$ be a finite field and let $X$ be a normal curve over $k$ which is geometrically connected. Let $\cl X := X \times_{\spec k} \spec \cl k$ be its base-change to the geometric point. Then we have a short exact sequence linking the geometric and arithmetic fundamental groups
$$ 1 \ra \pi_1 (\cl X) \ra \pi_1 (X) \ra G(\cl k/k) \ra 1.$$ 
 
\bt[de Jong, Gaitsgory] \label{thm:finite_over_O}
Let $\cl r: \pi_1 (X) \ra \GL_n (\F)$ be a continuous representation with $p > 2$ and $p \nmid nl$. Let $\psi: \pi_1 (X) \ra \GL_n(\F)$ be a continuous character lifting $\det \cl r$. Suppose that $\cl r |_{\pi_1 (\cl X)}$ is absolutely irreducible, then the deformation ring $\textup{R}_{\cl r}^{\textup{univ}, \psi}$ is finite flat over $\mO$.  
\et 

\bproof 
We first explain the relation between this theorem and de Jong's conjecture. In  \cite{dJ}, de Jong predicts the following \footnote{In fact de Jong's conjecture is stated for general integral schemes which is normal, separated and of finite type, we have restricted to curves for the purpose of this article.}
\vspace*{0.2cm}

\noindent \textbf{Conjecture.} \textit{ Let $L'$ be a local function field of characteristic $p \ne l$, and let $V$ be a finite dimensional $L'$-vector space. Let $r: \pi_1 (X) \ra \GL_n (V)$ be a continuous representation (with respect to the topology on $\GL_n(V)$ induced from $L'$). Then $r (\pi_1 (\cl X))$ is finite. } 
\vspace*{0.2cm}

Moreover, de Jong shows that the conjecture implies Theorem \ref{thm:finite_over_O} (Theorem 3.5 of \cite{dJ}). The conjecture has been proven by Gaitsgory in \cite{DG} using methods from geometric representation theory, hence the theorem follows. 
\eproof

\br Bockle and Khare also prove similar finiteness results under different assumptions on $\cl r$, for example Theorem 1.6 of \cite{BKh}. We cannot directly use their results in our construction later,  since they typically require that at some place $v \in S$, $\cl r|_v$ is absolutely irreducible with tame ramification.
\er

\subsection{A theorem of Moret-Bailly}  
The following lemma is due to Calegari, as an application of a theorem of Moret-Bailly. The lemma is stated for number fields in Calegari's article \cite{Calegari} Prop. 3.2, but same proof goes through in the function field case. 

\bl \label{lemma:Moret_Bailly}
Let $\Gamma$ be a finite group. Let $M/\F_l(t)$ be a finite extension, and $\Sigma$ a finite set of places $M$. Let $M_{\textup{Aux}}$ be a fixed auxiliary field over $M$. For each $v \in \Sigma$, fix a finite Galois extension $H_v/M_v$  together with a morphism
$$\phi_v: G(H_v/M_v) \hookrightarrow \Gamma_v := \im \phi_v \subset \Gamma.$$
Then there exists a finite extension $F/M$ and finite Galois extension $E/F$ such that 
\be
\item There is an isomorphism $\psi: G(E/F) \isom \Gamma$.
\item $E/M$ is disjoint from $M_{\textup{Aux}}$.
\item Every $v \in \Sigma$ splits completely in $F$. 
\item For each $v \in \Sigma$ and any place $w$ of $F$ above $v$,  there exists a place $w'$ of $E$ above $w$ such that $D_{w'|w} \subset G(E/F)$ is identified with $\Gamma_v$ under $\psi$. 
\ee
\el
\noindent To summarize,  starting from the data of 
$$\{M_{\textup{Aux}}/M, \Gamma,  \Sigma, \{\phi_v: G(H_v/M_v) \isom \Gamma_v \hookrightarrow \Gamma\}\}$$
we can create $E/F$ over $M$ with $\psi: G(E/F) \isom \Gamma$, such that places in $\Sigma$ split in $F$ and $\psi (D_{w' | w}) = \Gamma_v$ for suitable choices of $w'$, which fit together as follows. 
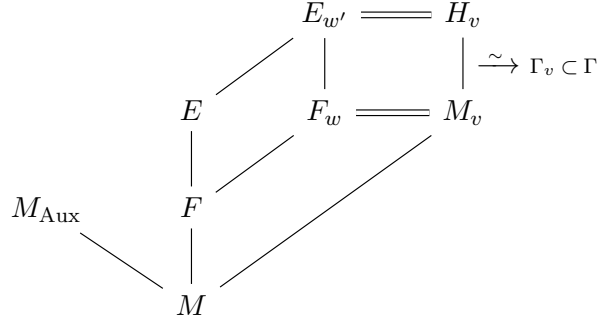
\[
\begin{tikzcd}
& & E_{w'}  \arrow[dash]{d}{}  \arrow[equal]{r}{}   \arrow[dash]{ld}{}    & H_v \arrow[dash]{d}{ \isom \Gamma_v \: \subset \: \Gamma}  \\
& E  \arrow[dash]{d}{}   & F_w  \arrow[dash]{ld}{} \arrow[equal]{r}{}  & M_v \arrow[dash]{lldd}{}  \\
M_{\textup{Aux}}  \arrow[dash]{rd}{}   & F  \arrow[dash]{d}{}  \\
&  M
\end{tikzcd}
\]

\subsection{Realization of local Galois representations}  

We will use Lemma \ref{lemma:Moret_Bailly} to realize $\cl \rho$ as local components of a global Galois representation. For this we first need to specify several choices of local decomposition groups (which will correspond to $\Gamma_v$ in Lemma \ref{lemma:Moret_Bailly}). In the application, we let $\Gamma := \GL_n (\F)$ and assume $p \ne n(n-1)$ as in the statement of Theorem \ref{thm:Main}.  

\subsubsection{$\Sigma_1$} \label{sss:Sigma_1} \indent 

 Take $M = \F_l (t)$ which is the function field of $X := \P^1_{\F_l}$, and let $x \in |X|$ be a closed point of degree $d$, such that $M_{x} \cong K$.  Let  $\Sigma_1$ to be the singleton $\{x\}$. 
 Fix a separable closure $\cl K \cong \cl M_{x}$, and identify $G_K = G (\cl M_x /M_x)$.  We started with the local Galois representation 
$\cl \rho: G_K \ra \im (\cl \rho) \subset  \GL_n (\F) =: \Gamma.$ Now we define 
$H_{x}:= (\cl M_{x})^{\ker \cl \rho } \subset \cl M_{x}$
and $\Gamma_{x}:= \im (\cl \rho) \subset \Gamma$. The map $\phi_{x}$ is given by $\cl \rho$, namely 
$$\phi_{x}: G(H_{x}/M_{x}) = G_K/\ker \cl \rho \xrightarrow{\: \: \cl \rho \: \:} \Gamma_{x} = \im \cl \rho .$$

\subsubsection{$\Sigma_2$} \indent 

Let $\Sigma_2$ be a singleton $\{z\}$ disjoint from $\Sigma_1$, such that $\zeta_p \notin M_{z}$.  For example, since $M = \F_l (t)$, and $x $ is a degree $d$ point in $X$ as in \ref{sss:Sigma_1} above, then $z$ could be taken to be origin of $X = \P^1_{\F_l}$. In this case $M_{z} \cong \F_l (\!(t)\!)$, and does not contain $\zeta_p$, since $p \nmid l - 1$ by our assumption. 

Let $\cl \rho_{\Sigma_2}$ be the trivial representation $ \cl \rho_{\Sigma_2} = \cl{\textup{triv}}: G_{M_{z}} \ra \{1\} \n \GL_n(\F).$
This time, $\Gamma_{z} = \{1\}$ is the trivial subgroup of $\Gamma$, and  $H_{z} := M_{z}$.  

\subsubsection{$\Sigma_3$} \label{sss:Sigma_3} \indent 

Let $\Sigma_3$ be a finite collection of places of size $\#|\GL_n(\F)|$ disjoint from $\Sigma_1$ and $\Sigma_2$, such that for each $u \in \Sigma_3$, $M_{u}$ contains $\zeta_p$. Let $\cl \rho_u: G_{M_u} \ra \GL_n(\F)$ be a collection of representations such that their images $\{ \Gamma_u:= \im (\cl \rho_u)\}_{u \in \Sigma_3}$ generate $\GL_n(\F)$. This can be achieved since we have taken $\Sigma_3$ to be large enough.

\subsubsection{$\Sigma_D$} \label{sss:D_compatible} \indent

Let $\Sigma_D = \{y^{1}, ..., y^{n}\}$ be a set of $n$ places of $L$ disjoint from $\Sigma_1$ -- $\Sigma_3$. For each $1 \le i \le n$, let $M_i := M_{y^i}$ be the local field at the place $y^i$. We will consider certain ''Weil--Deligne representations'' of $W_i: = W_{M_i}$ over finite fields, given by 
$\cl \rho'_i := \textup{Sp}_n = (\textup{Sp}_n (\textup{triv}), N)$ in the spirit of Definition \ref{def:special_reps}. Now let us fix a lift of Frobenius $\varphi$ and the $p$-part of a tame character   $t_p: I_K  \twoheadrightarrow \Z_p \twoheadrightarrow \F_p ,$
and define the corresponding ``Galois'' representation  $\cl \rho_i: W_{i} \ra \GL_n (\F)$  which extends to a continuous representation 
$ \cl \rho_{i} : G_{M_i} \ra \GL_n (\F).$
Note that $\cl \rho_i$ is indecomposable. In particular, any lifting of $\cl \rho_i$ to a continuous representation $G_{i}:= G_{M_i} \ra \GL_n (\mO)$ is indecomposable.  More precisely, we take (a conjugate of) the Galois representation 
$$\rho_i: G_{M_i} \ra \GL_n (\mO)$$ associated to the Weil--Deligne $(\rho_i' =\textup{Sp}_n(\textup{triv}), N)$ with $\textup{triv}: G_{M_i} \ra \mO^\times$, as in Definition \ref{def:special_reps}. Note that we need to taking a conjugate by elements in $\GL_n (L)$, since when $p < n$, the exponential of the tame character introduces denominators of positive $p$-adic valuations. Then we take the reduction mod $\lambda$ of $\rho_i$ to get $\cl \rho_i$.

Finally, we define the data $\{\Gamma_i, H_i/M_i, \phi_i\}_{1 \le 1 \le n}$ as above, coming from $\cl \rho_i$.

\subsubsection{Building a global representation} \indent 

\bp \label{cor:global_realization}  

Retain all assumptions from Theorem \ref{thm:Main}, and let $\chi$ be a fixed lifting of the character $\det \cl \rho$.  Then there is a global function field $F/\F_l (t)$ with representation 
$$ \cl r: G_F \ra \GL_n ( \F) $$ and a lift of $\cl r$: 
$$r: G_F \ra \GL_n (\mO_{\cl L}),$$ together with disjoint places $S_1$ and $S_D$ of $F$, satisfying the following conditions:  
\be  
\item  $\zeta_p \notin (\cl F)^{\ker \ad \cl r}$.  
\item $\cl r(G_{F(\zeta_p)}) = \GL_n (\F)$ 
\item For each place $x \in S_{1}$, $\cl r_x = \cl r |_{F_x}$ realizes $\cl \rho$. Namely  there is an isomorphism 
$$F_x \cong K \quad  \textup{ and } \:\:\:\: \cl r|_{G_{F_x}} \cong \cl \rho.$$ 
\item  For every place $y \in S_D$,  the restriction $r|_{F_y}$ at $y$ is isomorphic to the Galois representation associated to the special Weil--Deligne representation $\textup{Sp}_n (\textup{triv})$  (\textup{cf}. Definition \ref{def:special_reps}), where $\textup{triv}$ is the trivial character.
\item Both $r$ and $\cl r$ are unramified outside $S_1 \sqcup S_D$. 
\item At each $x \in S_1$, $\det r|_{G_{F_x}} \cong \chi.$
\item There exists a character $\eta_0: G_{F, S_1 \sqcup S_D} \ra \mO^{\times}$ such that $\det(\eta_0^{-1} \otimes r)$ has finite order.  
\ee
\ep

\br
In the requirement (4), note that the Galois representation of $W_K$ associated to $\textup{Sp}_n (\textup{triv})$ extends to $G_K$ (since the Frobenius eigenvalues are $1$, in particular $p$-adic integers in $\mO_L$). It is evidently indecomposable. We will abuse notation and denote this Galois representation of $G_K$ also by $\textup{Sp}_n (\textup{triv})$. 
\er 

\bproof 

\textit{Step 1}.  First we prove that there exists such $\cl r$ without the requirement $(5)$ by applying Lemma \ref{lemma:Moret_Bailly}, and then modify it slightly by applying the lemma one more time. 
 
We first apply Lemma \ref{lemma:Moret_Bailly} to the following setup: 
\bi
\item $M_{\textup{Aux}} = M$ where $M$ is given in \ref{sss:Sigma_1}. Most importantly, $M_{x} \cong K$. 
\item $\Gamma = \GL_n (\F)$.
\item $\Sigma = \Sigma_1 \sqcup \Sigma_2  \sqcup \Sigma_3 \sqcup \Sigma_D$.
\item $\{\Gamma_v, H_v/M_v, \phi_v\}_{v \in \Sigma}$ as specified above in \ref{sss:Sigma_1} - \ref{sss:D_compatible}. 
\ei
Now the lemma immediately produces a finite Galois extension $E/F$ with $$\cl r : G_{F} \twoheadrightarrow G(E/F) \xrightarrow[\cong]{\psi}\Gamma = \GL_n (\F) \hookrightarrow \GL_n (\cl \F).$$ Let $S_i$ be the places of $F$ above $\Sigma_i$ for $i = 1, 2, 3, D$. We sometimes abuse notations and denote a place in $\Sigma$ and a place above it both by the same letter. 

By construction $F$ and $\cl r$ satisfy property (3), and the part of (4) after reducing mod $\lambda$, and that $\cl r (G_F) = \GL_n (\F)$. 
Now we show that (1) is satisfied. Consider $\ad \cl r: G_F \ra \GL(\End (\F^n))$ where the action on $\End (\F^n)$ is given by conjugation. The restriction of $\ad \cl r$ at places in $S_2$ is given by 
$$(\ad \cl r) |_{G_{F_z}} = \ad (\cl r |_{G_{F_z}})  \cong \ad \cl \rho_{\Sigma_2},$$  
which implies that
$$G ( \cl F_{z} ^{\ker  \ad \cl r|_{z} } /  F_{z} ) \hookrightarrow G(\cl F^{\ker \ad \cl r}/F).$$ 
Now we have field extensions (after fixing $\cl F \hookrightarrow \cl F_z$)
\[
\begin{tikzcd}
\cl F_z^{\ker \ad \cl r|_{z}} \arrow[dash]{rd}{}   \arrow[dash]{d}{} &  \empty \\
F_z  \arrow[dash]{rd}{} & \cl F^{\ker \ad \cl r}      \arrow[dash]{d}{}     \\
& F 
\end{tikzcd}
\]
which implies that $ \cl F^{\ker \ad \cl r} \subset \cl F_z^{\ker \ad \cl r|_z}  = F_z$, but $F_z \cong M_{z}$ does not contain $\zeta_p$ by construction. 

Furthermore, (2) is satisfied since $\zeta_p \in F_{u}$ for each $u \in S_3$, and the image of $\cl r|_{G_{F_u}}$ generate $\GL_n (\F)$ by construction, but $G_{F_u}$ is a subgroup of $G_{F(\zeta)}$ for all $u$. This completes the construction of $E/F/M$ with $\cl r: G_F \twoheadrightarrow G (E/F) \isom \GL_n (\F)$ with requirement properties (1) -- (4).  \\

\textit{Step 2}.  Next we modify our construction of $F$ and $\cl r$ to satisfy (5). For this we consider the following setup:
\bi
\item $\sq M_{\textup{aux}} = E$, $\sq M = F$ from the construction above. 
\item $\sq \Gamma = \GL_n (\F)$
\item Let $\sq \Sigma = S_r \sqcup S$, where $S: = S_1 \sqcup S_2 \sqcup S_3 \sqcup S_D$, and 
$$S_r:= \{ \textup{ramified places of }  E/F \textup{ away from } S_1 \sqcup S_D\}.$$
We may assume that $S_r \ne \O$ otherwise we are already done. 
\item \bi
\item For each place $v \in S$, let $\Gamma_v = \{1\} $ and $\phi_v$ the trivial action. 
\item For each $v \in \sq \Sigma$, pick $w$ in $E$ above $v$, and let 
$$\sq \Gamma_{v} = \im (G(E_w/F_v) \hookrightarrow G(E/F) \isom \sq \Gamma),$$ and $\sq \phi_w$ the corresponding representation. 
\ei
\ei
Now Lemma \ref{lemma:Moret_Bailly} gives us extensions $\sq E / \sq F / F$  where $\sq E/\sq F$ has Galois group $\GL_n (\F)$ and is linearly disjoint from $\sq M_{\textup{aux}} = E$.  

Finally we replace $E/F$ (namely $\sq M_{\textup{aux}} / \sq M$) by $E' := E \sq E $ over $F' := F \sq E$, which is Galois with $G(E'/F') \cong \GL_n (\F)$. Now the representation $\cl r': G_{F'} \ra G(E'/F') \isom \GL_n(\F)$ satisfies the requirements (1)--(4), since all places $v \in S$ splits completely in $E'$, and by construction $\cl r' $ is unramified outside $S_1$ and $S_D$.  This completes the construction of $F$ and $\cl r$, satisfying requirements (1)--(5) in the lemma (except for the part of (4) regarding $r$).  \\

\textit{Step 3}. Now we carefully lift $\cl r$ to satisfy (4), (6) and (7). First consider $\cl \chi = \det \cl \rho: G_K \ra \F^\times$. Let $W(\cl \chi) : G_K \ra W(\F)^\times \subset \mO^\times$ be the Witt vector lifting of $\cl \chi$ then the mod $\lambda$ reduction of $\chi \cdot W(\cl \chi)^{-1}$ is trivial, hence we can take its $p^{th}$ root $\chi_0$ (since $n$ is invertible in $\mO$). Namely $\chi = \chi_0^n \cdot W (\cl \chi)$, where $\chi_0$ mod $\lambda$ is trivial. Now fix be a lifting $\eta_0$ of the trivial character $\cl{\textup{triv}}: G_{F, S} \ra \F^\times$ such that $\eta_0$ satisfies the following local conditions 
\be
\item At $x \in S_1$, $\eta_0|_{G_{F_x}} = \chi_0$. 
\item At $y \in S_D$, $\eta_0|_{G_{F_y}} = \textup{triv}$ is the trivial character. 
\ee
This leads to the following character 
$$ \eta := (\eta_0 )^n \cdot W(\det \cl r): G_{F, S} \ra W(\F)^\times. $$
Note that $\eta':= \eta \cdot (\eta_0)^{-n}$ is finite, and $\eta|_{G_{F_x}} = (\chi_0)^n \cdot W(\cl \chi) = \chi$ for each $x \in S_1$. 

We need to construct a lifting $r$ of $\cl r: G_{F, S} \ra \GL_n(\F)$ to $\mO_{\cl L}$, which has a fixed determinant $\eta$, and satisfying the following local condition:  
\bi
\item At each $y \in S_D$, we ask the lifting to be isomorphic to $\textup{Sp}_n$. 
\ei
To find such a characteristic $0$ lifting we proceed as follows. First we consider the global deformation ring of $\cl r: G_{F, S} \ra \GL_n (\F)$ with the finite character $\eta'$ and with no local conditions.  
We claim that  $ \textup{R}_{\cl r}^{\textup{univ}, \eta'}$ is finite over $\mO$. 

Let $F = F(\sq X)$ be the function field of the projective curve $\sq X$ over a finite field $k/\F_l$, and let $X = \sq X - S$. To the prove the claim it suffices to show that $\cl r (\pi_1 (\cl X))$ is absolutely irreducible so that we could apply Theorem \ref{thm:finite_over_O}. To see this we use our assumption on $\cl r$ at places $y \in S_D$. Let $\pi_1(\cl X)_0$ (resp. $\pi_1 (X)_0$) be the kernel of the map 
$$\cl r|_{\pi_1(\cl X)}: \pi_1 (\cl X) \ra \GL_n (\F) \xrightarrow{\det} \F^\times$$ (resp. of $\cl r: \pi_1 (\cl X) \ra \GL_n (\F) \xrightarrow{\det} \F^\times$). Now consider 
$$\pi_1(\cl X)_0 \ra \SL_n(\F) \twoheadrightarrow \textup{PSL}_n (\F).$$ 
By construction, we have that for any $y \in S_D$,  $\cl r|_{I_{F_y}}: I_{F_y}\hookrightarrow \pi_1(\cl X)_0$ and its image in $\SL_n (\F)$  is not contained in the center $Z(\SL_n(\F))$, therefore the map $\pi_1 (\cl X)_0 \twoheadrightarrow \textup{PSL}_n(\F)$ is nontrivial and hence surjective since $\textup{PSL}_n(\F)$ is a simple group (Note that by Remark \ref{remark:irreducible_component_inertia}  we have conveniently excluded the possibility of $\F$ being either $\F_2$ or $\F_3$). This shows that the representation $\cl r |_{\pi_1(\cl X)}$ is absolutely irreducible (since it is already absolutely irreducible when restricting to $\pi_1 (\cl X)_0$) and finishes the proof of the claim. \\

\textit{Step 4}.  Let us consider the global deformation ring of $G_{F,S}$ with local conditions at $S_D$. For each $y \in S_D$, let $\tau_y : = \textup{Sp}_n(\textup{triv}) |_{I_{F_y}}$, namely the Galois inertia type given by the restriction of the Galois representation associated to $(\textup{Sp}_n(\textup{triv}), N): W_{F_y} \ra \GL_n (\cl L)$.  Note that in particular $\det \textup{Sp}_n(\textup{triv}) = \textup{triv}$ is the trivial character. Let $\eta'_v = \eta' |_{G_{F_v}}$ for each place $v$ of $F$, then $\eta'_y$ is trivial at all $y \in S_D$. By Corollary \ref{cor:deformation_ring_character_union}, $\spec \textup{R}_{\cl r|_y}^{\square, \tau_y, \eta'_y}$ is an irreducible component of  $\spec \textup{R}_{\cl r|_y}^{\square,\eta'_y}$ and has dimension $n^2$. Now let the local deformation problems at each $y \in S_D$ be given by the quotient $\textup{R}_{\cl r|_y}^{\square, \tau_y, \eta'_y}$. Let $\dim \textup{R}_{\cl r, S}^{\textup{univ}, \eta'}$ be the universal global deformation ring with these local conditions at $S_D$, then $\dim \textup{R}_{\cl r, S}^{\textup{univ}, \eta'} \ge 1$  by Corollary \ref{cor:relative_dim_bound}.  

By the previous step, we also know that $\textup{R}_{\cl r, S}^{\textup{univ}, \eta'}$ is finite over $\mO$ since it is a quotient of $\textup{R}_{\cl r}^{\textup{univ}, \eta'} $. These two conditions on $\textup{R}_{\cl r, S}^{\textup{univ}, \eta'}$ imply that $\spec \textup{R}_{\cl r, S}^{\textup{univ}, \eta'}$ contains a $\cl L$-point $r': G_{F, S} \ra \GL_n(\mO_{\cl L})$.  Then twisting by the character $\eta_0$ we get a $\cl L$-point of $\textup{R}_{\cl r}^{\textup{univ}, \eta}$, which corresponds to a representation 
$$r = r' \otimes \eta_0: G_{F, S} \ra \GL_n (\mO_{\cl L}),$$ whose reduction mod the maximal ideal of $\mO_L$ agrees with $\cl r$. By construction, we know that $\textup{WD} (r|_{G_{F_v}})$ is isomorphic to $\textup{Sp}_n (\psi)$ for some (unramified) character $\psi$. Finally by the fact that the determinant $\eta_y$ is trivial, we conclude that $\psi$ is isomorphic to the trivial character. This concludes the lemma. 
\eproof 

\subsection{Langlands for function fields}  
 
\subsubsection{The global Langlands correspondence }  

We will need the following celebrated theorem of L. Lafforgue as an input in our proof of the main theorem.  Let $F/\F_l (t)$ be as above, $p$ a prime different from $l$. Let  $\mathcal{A}_n (F)$ be the set of isomorphism classes of cuspidal automorphic representations of $\GL_n (\A_F)$, with central characters of finite order. Let $\mG_n (F)$ be the set of isomorphism classes of $p$-adic representations of $G_F$ unramified at all but finitely many places, with finite determinant. 

\bt[Theorem $6.9$, \cite{La}] \label{thm:LL_function_field}
There exists a unique bijection
$$\mathcal{A}_n (F) \isom \mG_n (F), \qquad \pi \mapsto \sigma_\pi $$
from cuspidal automorphic representations to Galois representations, matching up Hecke and Frobenius eigenvalues (at unramified places), or equivalently, satisfying the compatibility between $L$-functions. 
\et

The global Langlands correspondence above is compatible with local Langlands correspondence [cf. Subsection \ref{sss:LLC}] in the following sense:
\bp[Proposition $7.4$, Corollary $7.5$, \cite{La}]  \label{prop:local_global}
 Let $X$ be a smooth proper curve over $\F_l$ with function field $F$, let $\pi \in \mathcal{A}_n (F)$ and $\sigma$ its corresponding Galois representation in $\mG_n (F)$. For any $x \in |X|$, let $K = F_x$ be the completion of $F$ at $x$, $\pi_x$ the local component of $\pi$ at $x$, and $\sigma_x$ the Weil--Deligne representation associated to $\sigma|_{G_K}$,  then $\textup{rec}_n (\pi_x) = \sigma_x$ under the local Langlands correspondence $\textup{rec}_n$ as in Subsection \ref{sss:LLC}.  
\ep

\bc \label{cor:realize_r_in_GL} 
Let $r: G_{F} \ra \GL_n (\mO_{\cl L})$ be the representation obtained in Proposition \ref{cor:global_realization}, then there exists a cuspidal automorphic representation $\pi$ corresponding to $r$ (in the sense that the Frobenius eigenvalues match up with Hecke eigenvalues). Moreover, $\pi$ is unramified outside $S_1 \sqcup S_D$, and at each $y \in S_D$, $\pi_y = \textup{rec}_n^{-1}(\textup{Sp}_n(\textup{triv}))$ is the Steinberg representation $\textup{St}_n$, as defined in Subsection \ref{sss:St_rep}. 
\ec

\bproof 
We use the same notation in the proof of Proposition \ref{cor:global_realization}. Since $\det r' = \det (\eta_0^{-1} \otimes r)$ is finite, $r' = \eta_0^{-1} \otimes r$ corresponds to a cuspidal automorphic representation $\pi'$ of $\GL_n (\A_F)$ by Theorem \ref{thm:LL_function_field}. Let $\epsilon_0$ be the character of $\A_{F}^\times /F^\times$ corresponding to $\eta_0$ in class field theory. Then $\pi := \pi' \epsilon_0$ is the desired cuspidal automorphic representation corresponding to $r$. At each place $v$ of $F$, we have 
$$ \textup{rec}_n ( \pi_v ) =  \textup{rec}_n ( \pi'_v ) \otimes \textup{rec}_1 ( \epsilon_0) = \textup{WD}(r \otimes \eta_0^{-1} |_{v}) \otimes \textup{WD}(\eta_0 |_v) = \textup{WD} (r|_{v}) $$
where $r|_v := r|_{G_{F_v}}$ as usual. In particular, at each $y \in S_D$, $\pi_y$ is the Steinberg representation. 
\eproof

\subsection{The Jacquet-Langlands correspondence for function fields}  
 
Next we apply Jacquet-Langlands to show that the Galois representation $r$ constructed in the Corollary \ref{cor:realize_r_in_GL} comes from a cuspidal automorphic representation associated to some central simple algebra $D$. To proceed, we fix the following convention for the rest of this section. 

\bi
\item Let $K$ be a finite extension of $\F_p (\!(t)\!)$ as usual, and $(\pi, V)$ be a smooth admissible representation of $\GL_n (K)$.
\item Let $\textup{G}^{r.s.}$ be the set of regular semi-simple elements in $\GL_n(K)$. (Recall that an element $g \in \GL_n (K)$ is regular semi-simple if the characteristic polynomial $P_g (X) \in k[X]$ of $g$ has distinct roots over $\cl K$. ) 
\item Let $\textup{G}^{r.e.} \subset \textup{G}^{r.s.}$ to be the subset of regular elliptic elements (Recall that an element $g \in \textup{G}^{r.s.}$ is regular elliptic if in addition its characteristic polynomial $P_g$ is irreducible). 
\ei

\subsubsection{Elliptic elements and elliptic representations} 

For the representation $(\pi, V)$ of $\GL_n(K)$, the Hecke algebra $\mH(\GL_n(K)) = \mC_{c}^{\infty} (\GL_n(K))$ acts on $V$ by 
$$\pi(f) (v) = \int_G f(g) \pi (g) v dg,$$ which makes $(\pi, V)$ an $\mH(G)$-modules where each $\pi(f)$ is a finite rank operator. The trace of $\pi(f)$ gives rise to the distribution character
$\Theta_\pi: \mH(G) \ra \C.$ By the theory of  Harish-Chandra (when $K$ has characteristic $0$) and of Lemaire (when $K$ has positive characteristic),  $\Theta_\pi$ is represented by the Harish-Chandra character $\chi_\pi$, which is locally constant on the regular semi-simple elements $\textup{G}^{r.s.}$. In other words, for each $f \in \mH(\GL_n(K))$, $$\Theta_\pi(f) = \int_{\GL_n(K)} f(g) \chi_\pi(g) dg .$$  

\begin{deflem} \label{deflem:elliptic} \indent 

\be 
\item A smooth admissible representation $(\pi, V)$ of $\GL_n(K)$ is elliptic if its Harish-Chandra character  $\chi_\pi$ is not identically $0$ on $\textup{G}^{r.e.}$. 
\item Every essentially square integrable representation is elliptic. 
\ee
\end{deflem}

For more details on elliptic representations and the claim above, we refer the readers to \cite{BR} and the references therein. 

\subsubsection{The global Jacquet-Langlands correspondence}  Now let $F$ be a global function field over $\F_l$, and let $D$ be a central division algebra over $F$ of dimension $n^2$. At each place $v$ where $D$ ramifies (in other words, is non-split), $(D \otimes F_v)^\times = \GL_{r_v} (D_v)$ where $D_v$ is central algebra over $F_v$ of dimension $d_v^2$, such that $n =r_v d_v$. For simplicity we assume that at all non-split places, the local invarians are $1/n$, so $r_v = 1$ and $(D \otimes F_{v})^\times = D_v^\times$. The following is a special case of $D$-compatible automorphic representations defined in \cite{BR}.

\bd \label{def:D_compatible}
As above, let $D$ be a central division algebra such that at all non-split places $v$, $d_v = n$. Then a discrete series representation $\pi = \bigotimes' \pi_v$ of $\GL_n (\A_F)$ is $D$-compatible if for all $v$ where $D$ is non-split, the local component $\pi_v$ is elliptic. 
\ed
 
Write $\A = \A_F$ and fix an element $a \in \A^\times$ of degree $1$, which gives rise to the central lattice $J = a^\Z$. Now let $G$ be either $G_D$ or $\GL_n$, note that $G(F) \minus G(\A)/J$ has finite measure. We consider the set $\mA_{\textup{ds}} (G(\A))$ of discrete series representations of $G(\A)$. Here by a discrete series we mean an irreducible subrepresentations of the right translation $R_{G, \chi}$ of $G(\A)$ acting on $L^2_{\infty}(G(F)\minus G(\A)/J, \psi_0)$ for some smooth complex quasi-character $\psi_0: F^\times \minus \A^\times/J \ra \C$, where \footnote{Strictly speaking, in the definition below, we follow \cite{La} and require functions in $L^2 (G(F)\minus G(\A)/J V, \psi_0)$ to have finite norm.}
$$L^2_{\infty}(G(F)\minus G(\A)/J, \psi_0) = \varinjlim_{V \textup{ compact open}}L^2 (G(F)\minus G(\A)/J V, \psi_0).$$

\bt[Theorem $3.2$ \cite{BR}]   \label{thm:global_JL}
Retain all notations from above, where $D$ is a division algebra over $F$ satisfying the assumptions in Definition \ref{def:D_compatible}. There exists a unique injective map 
$$\textup{JL}: \mA_{\textup{ds}} (G_D (\A)) \longrightarrow \mA_{\textup{ds}} (\GL_n (\A))$$ from the set of discrete series of $G_D(\A)$ to the set of the discrete series of $\GL_n(\A)$, which is compatible with the local Jacquet-Langlands correspondence. The image of $\textup{JL}$ contains precisely the $D$-compatible discrete series of $\GL_n (\A)$. 
\et

\subsubsection{Multiplicity one and strong multiplicity one theorems}  We will also need the multiplicity one and strong multiplicity one theorems for the group $G_D$.  

\bt[Theorem $3.3$ \cite{BR}] \label{thm:global_JL_multiplicity} \indent 

\be
\item If $\pi$ is a discrete series representation of $G_D (\A)$, then $\pi$ appears with multiplicity one (in the discrete spectrum). 
\item If $\pi_1$ and $\pi_2$ are two discrete series of $G_D (\A)$ such that for almost all places $v$ of $F$, $\pi_{1, v} \cong \pi_{2, v}$, then $\pi_1= \pi_2$ as subrepresentations of $L^2(D^\times \minus G_D (\A)/J)$, namely the strong multiplicity one theorem holds for $G_D$. 
\ee
\et

\subsection{Global realization as automorphic representations} \indent

\subsubsection{Central division algebra} \label{sss:cda} \indent

\be
\item Let $F$ be the global function field over $\F_l (t)$ constructed in Proposition \ref{cor:global_realization}.  Let $\Sigma_D$ be $n$ places as specified in \ref{sss:D_compatible}, and $S_D$ the places of $F$ lying above those of $\Sigma_D$. 
\item Let $D$ be a central division algebra over $F$ of rank $n^2$, such that $D$ is non-split precisely at $S_D.$ In this case we have $G_D (\A^{S_D}) = \GL_n (\A_{F}^{S_D}).$
\ee  

\br  Such a central division algebra $D$ exists since we have 
$$ 0 \ra \textup{Br} (F) \ra \oplus_{v \in F} \textup{Br} (F_v) \ra \Q/\Z \ra 0 $$
so for example we may take local invariants to be $\frac{1}{n}$ for each $y \in S_D$. Since each $y_0^i \in \Sigma_D$ splits, these invariants sum to $0$ in $\Q/\Z$, and has least common multiple $= n$, so they correspond to a central division algebra of rank $n^2$ over $F$. 
\er

\subsubsection{\empty}   Now we state the main result of this section. 

\bp \label{thm:global_realization} 
Let $\cl \rho: G_K \ra \GL_n (\F)$ as before. Then there is a global function field $F/\F_l (t)$, a central division algebra $D$ over $F$, with ramified places $S_D $ as in \ref{sss:cda}, together with a representation
$$ \cl r: G_F \ra \GL_n ( \F) $$ and a lift of $\cl r$: 
$$r: G_F \ra \GL_n (\mO_{\cl L}),$$ such that  
\be  
\item $r = r(\pi)$ comes from a cuspidal automorphic form associated to $D$. 
\item There exists a finite set of places $S_1$ such that $\cl r_x = \cl r |_{F_x}$ realizes $\cl \rho$ for $x \in S_1$.  
\item $r$ is the special representation $\textup{Sp}_n (\textup{triv})$ at all places in $S_D$. 
\item $r$ and $\cl r$ are unramified outside $S_1$ and $S_D$.  
\item  $\zeta_p \notin (\cl F)^{\ker \ad \cl r}$. 
\item $\cl r(G_{F(\zeta_p)}) = \GL_n (\F)$ 
\ee 
\ep

\bproof 
We have constructed $F$ and $r$ in Proposition \ref{cor:global_realization}.  From the global Langlands correspondence, we have realized $r = r(\pi)$ for a unique cuspidal automorphic representation $\pi$ of $\GL_n (\A_F)$ as in Corollary \ref{cor:realize_r_in_GL}. In the corollary, we know that for each $y \in S_D$, $\pi_y$ is the Steinberg representation $\textup{St}_n$, so in particular elliptic by part (2) of Lemma \ref{deflem:elliptic}. This shows that $\pi$ is is $D$-compatible, so Theorem \ref{thm:global_JL} applies to $\pi$ 
\footnote{Note that we have $\pi = \pi' \epsilon_0$ in the proof of Corollary \ref{cor:realize_r_in_GL}, and $r = r(\pi') \otimes \eta_0$. Strictly speaking, only $\pi' \in \mA_{\textup{ds}} (\GL_n(\A))$ so we need to apply Theorem \ref{thm:global_JL} to $\pi'$ instead of $\pi$ and twist the result back by $\epsilon_0$, namely let $\textup{JL}^{-1} (\pi):= \textup{JL}^{-1} (\pi') \otimes \epsilon_0$.} 
and gives rise to $\textup{JL}^{-1}(\pi) \in \mA_{\textup{ds}} (G_D(\A))$. The rest of the requirements are inherited from Proposition \ref{cor:global_realization}. We will (as in the statement of property (1) above) abuse notation and denote $\textup{JL}^{-1} (\pi)$ also by $\pi$, whenever no confusion arises. 
\eproof

\section{The patching argument over function fields} \label{sec:Patching}

In this section we apply the Taylor--Wiles--Kisin patching method to create the module $M_{\infty} (\sigma)$ with the desired properties as in Proposition \ref{prop:main_prop}. The patching method  \cite{CHT, K1}, introduced by Taylor and Wiles and modified by Kisin, is now a standard argument. We nevertheless give a detailed enough argument since the argument over function fields (where $l \ne p$) is slightly different from number fields. For example, as all local components have characteristic $l$, which is different from $p$, the local deformation rings never get ``too large'',  in particular we do not need  statements of the form 
$ \sum_{v | \infty} H^0 (G_{F_v}, \ad \cl r) =  [F:\Q] \cdot (\dim G - \dim B)$ as remarked in \cite{CHT} and can directly work with $\GL_n$ instead of unitary groups. 
 
\subsection{Automorphic forms} The material in this subsection is similar to the discussion of chapter \S3 of \cite{CHT}. The major difference is that, we designate the places $S_D$ to serve as ``$S_\infty$'' in the number field case, and modify the notion of ``sufficiently small subgroups'' to mean ``sufficiently small away from $S_D$ and maximal at $S_D$'', which is in particular compact mod center at $S_D$. This has two effects, (i). the double quotient $D^\times \minus (D \otimes \A_F^{D})^\times /U^D$ (away from $S_D$) becomes finite (see Remark \ref{remark:double_quotient_finite}); (ii). this allows us to restrict to \textit{the} special representation at places $v \in S_D$ in the Galois side, so the deformation space at $v \in S_D$ is trivial, which is convenient for us for otherwise a more relaxed condition might cut out multiple irreducible components in the local Galois deformation rings at $v \in S_D$. 

\subsubsection{Notation} Retain notations from the previous section. From Proposition \ref{thm:global_realization} and its proof we have 

\bi
\item $r = r (\pi)$ is the irreducible Galois representation associated to an automorphic representation $\pi$ of  the division algebra $D$, where $\cl r|_{v} \cong \cl \rho$ for $v \in S_1$. 
\item $G_D$ is the algebraic group over $F$ defined by $D^\times$, which is non-split at $S_D$ with all its local invariants equal to $n$. 
\item Let $U \subset G_D (\A_F)$ be a compact open subgroup. Write $U = U^D U_D$ where $U_D = \prod_{v \in S_D} U_v$ and $U^D = \prod_{v \notin S_D}' U_v$ with all but finitely many $U_v = \GL_n (\mO_{F_v}).$
\item Write $U_1 := \prod_{v \in S_1} U_v.$ For an element $u \in U$, write $u_1$ for its projection to the $U_1$ factor. 
\item $\eta = \det r: G_{F} \ra \mO_L^\times$, which is unramified outside $S_1$ and trivial at $S_D$. Let $\psi: \A_F^\times/F^\times \ra \mO_L^\times$ be the ad\'elic character corresponding to $\eta$ under class field theory.
\item Let $v_0 $ be a place of $F$ satisfying the following two conditions (i). its residue field $k_0$ has size $p \nmid (\# k_0 - 1) $; (ii). $\ad \cl r |_{v_0}$ is trivial.  For example we get such a $v_0$ by carefully choosing the places in $S_2$ in the proof of Proposition \ref{cor:global_realization}, so that the size of their residue fields are not congruent to $1$ mod $p$. 
\item Let $T = S_1 \sqcup S_D \sqcup \{v_0\}$. 
\item For an auxiliary set of places $Q$ disjoint from $T$, let $S := T \sqcup Q$.  
\item Finally, as in the beginning of the paper, fix $L$ large enough, an isomorphism $\iota: \cl L \isom \C$, and write $\mO = \mO_L$. 
\ei

\br
The place $v_0$ is chosen to help us construct a sufficiently small subgroup in Subsection \ref{sss:sss}. This will come at a cost of enlarging the space of automorphic forms by relaxing conditions at $v_0$, but $v_0$ is carefully chosen (namely that $\ad \cl r|_{v}$ is trivial) so that the Galois deformation space at $v_0$ remains ``unramified'' (cf. Lemma \ref{lemma:def_ring_at_v}). 
\er 

\subsubsection{Space of cusp forms} \label{sss:space_auto_forms} Now we define space of cuspidal automorphic forms with $\mO$ coefficients and character $\psi$. As mentioned above, let $\psi: \A_F^\times/F^\times \ra \mO^\times$ be the character given by composition  $\psi: \A_F^\times/F^\times \hookrightarrow G_F^{ab} \xrightarrow{\:\eta\:} \mO^\times$ of $\eta$ with the global Artin map (which in the case of function field is an injection). 

\bd Let $M$ be any finitely generated $\mO$-module with an action of $U_1 = \prod_{v \in S_1} U_v$, then we define $S_{\psi} (U, M)$ to be the space of functions
$$S_{\psi} (U, M) = \{f: D^\times \minus G_D(\A_F) \longrightarrow M \}$$ 
such that $f(zgu) = \psi(z) u_1^{-1} f(g)$ for all $z \in \A_{F}^\times, g \in G_D(\A_F), \textup{and } u \in U$. 
\ed 

As in \cite{CHT}, we need to make sure that this is finitely generated. 
\br \label{remark:double_quotient_finite}
Suppose that $U_D$ is \textit{maximal}, in other words, for all $v \in S_D$, $U_v = \mO_{D_v}^\times \subset D_v^\times$. Note that by our construction, $\psi|_{F_{v}^\times}$ is trivial. Therefore, for any $g \in G_D(\A_F)$, if we write $g = g^D g_D$ with $g_D \in \prod G_D(F_v)$, we have $f(g) = f(g^D)$ for all $f \in S_{\psi} (U, M).$ In other words, under the maximality assumption above, 
$$S_{\psi} (U, M) = \{f: D^\times \minus G_D(\A_F^D) \longrightarrow M\}$$ such that $f(z g u) = \psi (z) {u_1}^{-1} f(g)$ for all $z \in (\A_F^D)^\times,  g \in G_D(\A_F^{D})$ and $u \in U^D$. Note that we have:
$$ \# \vline \:D^\times \minus G_D(\A_F^D) /U^D \: \vline < \infty. $$ 
\er 

\bc 
For $U^D$ maximal and $M$ a finitely generated $\mO$-module, $S_{\psi} (U, M)$ is then a finitely generated $\mO$-module. If in addition $M$ is $\lambda$-torsion free, then $S_{\psi} (U, M)$ is a finite free $\mO$-module.
\ec

\br Since $G_D$ is anisotropic mod center, there is no proper parabolic subgroup defined over $F$, so every automorphic form is automatically cuspidal. A crucial step in the patching argument is to create a free module $S(U, \mO)$ over the ring $\mO[\Delta]$ of diamond operators, for this we cannot work directly with $\GL_n$, so the recently established Jacquet-Langlands correspondence over global function fields (see the previous section) is needed. 
\er 

Our space of automorphic forms with $\mO$ coefficients compares to the usual space of automorphic forms with $\C$ coefficients in a straightforward manner. Let $\psi_{\iota} = \iota \circ \psi$ be the smooth character of $\A_F^\times /F^\times$ obtained by composing with $\iota: \mO \hookrightarrow \C$, and let $S_{D, \psi_{\iota}}$ be the space of functions $\phi: D^\times \minus G_D(\A_F) \ra \C$ satisfying $\phi(zg) = \psi_{\iota} (z) \phi(g)$ and that there exists a compact open subgroup $V \subset G_D(\A_F)$ such that $\phi(gv) = \phi (g)$ for all $v \in V$.  As usual $G_D(\A_F)$ acts on $S_{D, \psi_{\iota}}$ by right translation. From the definition, we have  
$$S_{\psi}(U, \mO) \otimes_{\mO, \iota} \C \isom S_{D, \psi_{\iota}, U} := S_{D, \psi_{\iota}}^U $$ given by $f \mapsto \iota \circ f$. The action of $G_D (\A_F)$ on both sides are compatible via the embedding $\iota$, and note that $U$ acts trivially on both sides. Moreover, $G_D(\A_F)$ acts semi-simply on $S_{D, \psi_{\iota}}^U$, which decomposes into a direct sum of irreducible cuspidal automorphic representations 
$$S_{D, \psi_{\iota}}^U \cong \underset{\pi}{\oplus}  \Big( \pi_S^{U_S} \otimes \underset{v \notin S}{\otimes'} \: \pi_v ^{U_v} \Big) $$
where the direct sum is over cuspidal automorphic representations $\pi$ of $G_D (\A_{f})$ of central character $\psi_{\iota}$ and unramified outside $S$, and each direct summand appears with multiplicity one by Theorem \ref{thm:global_JL_multiplicity}. 

\subsubsection{Sufficiently small subgroups}  \label{sss:sss}

We need (a slight variant of) the notion of sufficiently small \footnote{This is called \textit{good} in \cite{S}, also in \cite{GK} for the simpler case of $\GL_2$. } subgroups of $G_D (\A_F) = (D \otimes \A_F)^\times$ used in \cite{CHT} \S 3 

\bd  \label{def:s_small}
A compact open subgroup $U$ of $(D \otimes \A_F)^\times$ is \textit{sufficiently small} (with respect to a fixed prime $p$) if there exists a place $v \notin S_1 \sqcup S_D$, such that $U_v$ does not contain any element of order $p$. 
\ed

\bl  \label{lemma:exact_and_free}
Suppose that $U = U^D U_D$ is sufficiently small and that $U_D$ is maximal (namely that $U_D = \prod_{v \in S_D} \mO_{D_v}^\times$). Let $\textup{Rep}_{\mO}^{\textup{f.l.}}(U_1)$ be the category of finitely generated $\mO$-modules with $U_1$ actions.
\be
\item The functor $\textup{Rep}_{\mO}^{\textup{f.l.}}(U_1) \longrightarrow \textup{Mod}_{\mO}$ sending 
$$M \mapsto S_{\psi} (U, M)$$ is exact. 
\item Let $A$ be any finitely generated $\mO$-module (with no action of $U_1$), then 
$$S_{\psi} (U, M \otimes_{\mO} A) \cong S_{\psi} (U, M) \otimes_{\mO} A. $$
\item Let $V \n U$ be a normal open compact subgroup of $p$-power index, then there exists an isomorphism 
$$S_{\psi} (V, M) \isom S_{\psi} (U, M) \otimes_{\mO} \mO[U/V].$$ In particular, if $M$ is $\lambda$-torsion free, $S_{\psi} (V, M)$ is a free module over $\mO [U/V]$. 
\ee 
\el

\bproof 
The only minor difference compared to Lemma 3.3.1 in \cite{CHT} is that we have specified a central character for our space of automorphic forms. We briefly sketch the proof to convince the reader that the argument is essentially the same. Let  $G_D (\A_F^D) = \sqcup_{i \in I} D^\times g_i U^D (\A_F^{D})^\times$ be a finite disjoint union, and let 
$$G_i := \Big(g_i^{-1} D^\times g_i \cap U^D  (\A_F^{D})^\times \Big)/F^\times. $$
First observe that the injective map 
$$S_{\psi} (U, M) \hookrightarrow \oplus_{i \in I} M^{G_i}, \quad f \mapsto \{f(g_i)\}_{i \in I}$$ is in fact an isomorphism. To see this we need to identify the image. The only constraints to put on $\phi (g_i)$ arise as follows: suppose there exists $\delta \in D^\times, u \in U^D, z \in (\A_F^D)^\times$ such that $ g_i = \delta^{-1} g_i u z$, then we must have  $f(g_i) = \psi (z) u^{-1} f(g_i) = (uz) \cdot f (g_i).$  

We then claim that $G_i$ is a finite group of order prime to $p$, in other words, if $g \in g_i^{-1} D^\times g_i \cap U^D  (\A_F^{D})^\times$ is an element such that $g^p \in F^\times$, then $g \in F^\times$. Write $g = u z$ for $u \in U^D, z \in (\A_F^{D})^\times$, then $ u = gz^{-1} \in g_i^{-1} D^\times g_i \cap U^D$ so it has finite order. Since $U$ is sufficiently small, $u$ has order prime to $p$, say $u^m = 1$ where $p \nmid m$. In particular, $g^m = z^m \in D^\times \cap (\A_F^{D})^\times = F^\times$. Since $g^p \in F^\times$ and $p \nmid m$, clearly $g \in F^\times$. Therefore the first two claims follow immediately, and the proof of the last claim is similar to the corresponding part in \cite{CHT} (see \cite{Notes} Lemma 2.6 for a detailed proof.) 
\eproof

\subsection{Choice of level}

\subsubsection{Choice of the compact open subgroups $U_Q$} \label{sss:choice_of_U}

Let $Q$ be a finite set of places of $F$ consisting of Taylor--Wiles primes (see Subsection \ref{sss:TW_primes}) disjoint from $T$. Define the compact open subgroups $U_Q$ and $\sq U_Q$ as follows: 
\be
\item For $v \notin S_D \sqcup \{v_0\} \sqcup Q$, let $U_v = \GL_n (\mO_{F_v})$
\item For $v \in S_D$, let $U_v = \mO_{D_v}^\times$
\item For $v = v_0$, let $U_{v_0} \n \GL_n (\mO_{F_{v_0}})$ be the Iwahori subgroup $\textup{Iw} (v_0)$, consisting of matrices which are upper-triangular mod $\varpi_{v_0}$.  
\item For $v \in Q$, we will choose $U_v = U_1 (v)$ and $\sq U_v = U_0 (v)$ as follows: first let $U_0 (v)$ be the mirabolic subgroup of $\GL_n (\mO_{F_v})$, consisting of matrices whose last row is congruent to $(0, \cdots, 0, *)$ mod $\varpi_v$. Let $\Delta_v$ be the maximal $p$-power quotient of $k_v^\times =  (\mO_{F_v}/\varpi_v)^\times = k_v^\times$. For any element $(M_{ij}) \in U_{0} (v)$, let $M_1= (M_{ij})_{1 \le i, j \le n-1}$  be the top left $(n-1)\times (n-1)$ minor, and $M_2 = (M_{n,n})$ be the bottom right entry, then $M \mapsto M_2 /\det M_1 \mod \varpi$ defines a map 
$$U_0(v) \longrightarrow k_v^\times \longrightarrow \Delta_v. $$ 
Define $U_1 (v)$ to be the kernel of this projection and set $U_v = U_1 (v)$. This way $U_v \n \sq U_v$ is a normal subgroup with p-power quotient $\sq U_v/U_v \cong \Delta_v$. 
\ee

\noindent Finally, let $U_Q := \prod' U_v$ and $\sq U_Q:= \prod'_{v \notin Q} U_v \times \prod_{v \in Q } \sq U_v$. Note that, since $U_{v_0} = \textup{Iw}(v_0)$ does not contain any element of order $p$, both $U_Q$ and $\sq U_Q$ are sufficiently small in the sense of Definition \ref{def:s_small}.

\subsubsection{$U_1(v)$-spherical representations} From the definition of $U_1(v)$ above, we have a short exact sequence of abelian groups 
$$ 0 \ra U_1(v) \ra U_0 (v) \ra \Delta_v \ra 0.$$ 
The following lemma matches up the choice of $U_v$  at each $v \in Q$ with a variant of the Taylor--Wiles deformation condition, which is given in Definition \ref{def:TW_conditions}. This is the analog of Lemma 3.1.5 and 3.1.6 in \cite{CHT}. 

\bl \label{lemma:same_tame_inertia}
Let $F_v$ be a local field of residue characteristic $l$ and let $\pi$ be a generic smooth irreducible representation of $G = \GL_n (F_v)$. Suppose the ($p$-adic) Galois representation representation $\textup{r}_p(\pi) = \textup{rec}_n(\pi \otimes |\cdot|^{(n-1)/2})$ exists (\textup{i.e.} the eigenvalues of $ \textup{rec}_n(\pi \otimes |\cdot|^{(n-1)/2}) (\varphi_v)$ are $p$-adic units for some lift of Frobenius $\varphi_v$) \footnote{Note that we have chosen a different normalization from \cite{CHT}}. Suppose that $\pi^{U_1(v)} \ne 0$, then there exists a short exact sequence of Galois representations 
$$ 0 \ra s \ra \textup{r}_p(\pi) \ra \psi \ra 0 $$ 
where $s$ is an $(n-1)$-dimensional tamely ramified representation such that $s|_{I_K}$ consists of scalars, and $\psi$ is a $1$-dimensional character. 
\el

The proof is similar to the proof of Lemma 3.1.5 in \cite{CHT} and Lemma 5.11 in \cite{Thorne}. In fact, instead of reproducing the entire argument, we use the following lemma to reduce to the case proved there. Before we proceed with the proof, let us define more compact opens subgroups 
\be
\item  Let $\textup{Iw}_1(v) \subset \textup{Iw}(v)$ be the subgroup of the Iwahori subgroup consisting of matrices with $1$'s on the diagonal after reducing mod $\varpi = \varpi_v$.
\item Let $P_1(v) \subset U_0 (v)$ be the kernel of 
$U_0 (v) \ra \Delta_v$ given by $(M_{ij}) \mapsto \det M_1 $ where $M_1$ is the top left $(n-1)\times (n-1)$ minor as before. 
\item Let $\breve U_1(v) \subset U_0(v)$ be the kernel of $U_0 (v)\ra \Delta_v$ where the map is given by $(M_{ij}) \mapsto \det M_{nn} \in \Delta_v$. 
\ee 
In other words, $P_1$ consists of matrices congruent to $\begin{pmatrix} M_1 & * \\ 0 & M_2  \end{pmatrix}$ mod $\pi$, where $M_1$ is $(n-1)\times(n-1)$ and $M_2$ is $1\times1$, such that $\det M_1 \mapsto 1 \in \Delta_v$; $\breve U_1(v)$ consists of of those such that $M_2 \mapsto 1 \in \Delta_v$ (while $U_1(v)$ consists of those $\frac{\det M_1}{M_2} \mapsto 1 \in \Delta_v$). In particular $\textup{Iw}_1(v) \subset U_1(v) \cap P_1(v) \cap \breve U_1(v).$ 

\bl \label{lemma:same_tame_inertia_ii}
Retain the setup in the Lemma \ref{lemma:same_tame_inertia}. There exists a tamely ramified character $\sq \chi_1$ such that  
$$(\pi \otimes \sq \chi_1^{-1})^{P_1(v)} \ne 0$$
\el

\bproof 
Since $\pi$ is generic and that $\pi^{\textup{Iw}_1(v)} \ne 0$, $\pi$ can be realized as a sub-representation of a tamely ramified principal series, namely there exists a $G$-equivariant embedding 
$$\pi \subset \textup{Ind}_{B} ^{G} \chi_1 \otimes \cdots \otimes \chi_n$$ 
where $B \subset G$ is the standard Borel and each $\chi_i$ is tamely ramified. (For example, this follows from part (1) of Lemma 3.1.6 in \cite{CHT}. We remark the references contained there also treat the case when $F_v$ has positive characteristic $l$, where no restriction needs to be placed on $l$ since $G = \GL_n(F_v)$).

Now let $\sq \chi_1 = \chi_1 \; |\cdot|^{(n-1)/2}$. Let $f \in \pi^{U_1(v)} \subset (\textup{Ind}_{B} ^{G} \chi_1 \otimes \cdots \otimes \chi_n)^{U_1(v)}$ be a nonzero element, we want to show that as a vector in the representation $\pi \otimes \sq \chi_1^{-1}$, $f$ is invariant under right translation of $P_1(v)$. 

First, for each $1 \le i \le n$, define $w_i \in G$ to be the element 
$$ w_i = \begin{pmatrix} 1_{i-1} & 0 \\ 0 & a_{n-i+1} \end{pmatrix} \qquad \text{where } a_{n-i+1} =  \begin{pmatrix}   & & 1 \\  & \text{\reflectbox{$\ddots$}} &  \\ 1 & & \end{pmatrix}.$$
$G$ admits the following decomposition 
$$G = \bigsqcup_{1 \le i \le n} B w_i U_1(v)$$ 
(compare to the decomposition used in the proof of Lemma 3.1.3 in \cite{CHT}).  Next we define functions $f_i$ on $G$ for $ 1 \le i \le n$ as follows: 
\bi
\item If $f (w_i) \ne 0$, then define $f_i$ by $f_i (x) = \sq \chi (b) f (w_i)$ for each $x = b w_i u \in B w_i U_1(v)$ and $0$ if $x \notin B w_i U_1(v)$. In this case $f_i$ has support equal to $B_v w_i U_1(v)$. 
\item If $f(w_i) = 0$, then set $f_i = 0$. 
\ei Clearly each $f_i$ is well defined and since $f \ne 0$, at least one of the $f_i$ is nonzero. By construction  $f = \sum_i f_i$. 

For any $c \in \mO_v^\times$, let $\gamma_c$ be the diagonal matrix with the last (i.e. bottom right) entry equal to $c$ and all other entries equal to $1$, and let $\beta_c$ be the diagonal matrix with the first (i.e. top left) entry equal to $c$ and all other entries equal to $1$.   We now claim that for each $i$, $\gamma_c$ acts on $f_i$ by scalar $\sq \chi_1(c)$: 
$$\gamma_c \cdot f_i = \sq \chi_1(c) f_i.$$ It suffices to assume that $f_i \ne 0$. Right translation by $\gamma_c$ does not change the support of $f_i$, so it suffices to compute $(\gamma_c \cdot f_i) (w_i)$. There are two cases,
 \bi
\item If $i = 1$, $w_1 = a_n$ is the anti-diagonal matrix, and 
$$ f_1 (w_1  \gamma_c )= f_1 (\beta_c w_1) = \sq \chi_1(c) f_1(w_1)$$
\item $ 2 \le i \le n$, then 
$$ f_i (w_i  \gamma_c ) = f_i ( w_i \beta_c \beta_c^{-1} \gamma_c ) = f_i(w_i \beta_c) $$
since $\beta_c^{-1} \gamma_c \in U_1(v)$. Therefore, 
$$ f_i (w_i  \gamma_c )  = f_i (\beta_c w_i) = \sq \chi_1(c) f_i(w_i).$$
\ei 
Consequently, for any $c \in \mO_v^\times$,  $\gamma_c \cdot f = \sq \chi_1(c) f$. Therefore, for any $g \in P_1(v)$, the action of $\pi \otimes \sq \chi_1^{-1}(g)$ on $f$ is given by 
$$ \sq \chi_1^{-1} (\det g) \; g \cdot f =  \sq \chi_1^{-1} (\det g) \; \Big((\gamma_{\det g}\; g') \cdot f\Big) =  \sq \chi_1^{-1} (c) \; \gamma_{\det g} \cdot f  $$  
where $g' := \gamma_{\det g}^{-1} g \in U_1(v)$. By the computation above, we have $\gamma_{\det g} \cdot f = \sq \chi_1^{-1} (\det g) \; f$, therefore 
$(\pi \otimes \sq \chi_1^{-1})(g) \cdot f = f$, hence the claim. 
\eproof

\bproof[Proof of Lemma \ref{lemma:same_tame_inertia}]
By Lemma \ref{lemma:same_tame_inertia_ii}, we may assume that $\pi^{P_1(v)} \ne 0$. The conclusion is clear if $\pi$ is unramified (i.e. $\pi^{\GL_n(\mO_v)} \ne 0$). Suppose that $\pi^{U_0 (v)} \ne 0$ but $\pi^{\GL_n (\mO_v)} = 0$, then $\pi^{\sq U_1 (v)} \ne 0$, and the conclusion follows from Lemma 3.1.5 in \cite{CHT}. Finally the case where $\pi^{P_1(v)} \ne 0$ but $\pi^{U_0(v)} \ne 0$ is taken care of by Lemma 5.11 in \cite{Thorne}. 
\eproof

\subsubsection{Hecke algebra}  

Now let us go back to the setup in the beginning of this section. For each $v \notin S$ (so $U_v = \GL_n (\mO_{F_v})$), and $1 \le j \le n$, define the Hecke operator $\textup{T}_{v}^{(j)}$ to be the double coset operator 
$$\textup{T}_v^{(j)}:= \bigg[ U_v \begin{pmatrix} \varpi_v 1_j & 0 \\ 0 & 1_{n-j} \end{pmatrix} U_v \bigg]$$
on the space $S_{\psi} (U, M)$ of automorphic forms and define $\T^{(S)} = \mO [\textup{T}_v^{(j)}, (\textup{T}_{v}^{(n)})^{-1}]_{v \notin S}$ to be the polynomial ring over $\mO$ generated by the Hecke operators $\textup{T}_{v}^{(j)}$ and $(\textup{T}_{v}^{(n)})^{-1}$ for all $v \notin S$. Finally define $\T_{\psi}^{(S)} (U, M) := \im \big( \T^{(S)} \ra \End (S_{\psi} (U, M)) \big).$
By the discussion on Subsection \ref{sss:space_auto_forms} and Theorem \ref{thm:global_JL_multiplicity}, we have an isomorphism of algebras
$$\T_{\psi}^{(S)} (U, M) \otimes_{\mO, \iota} \C \isom \prod_{\pi} \C $$ given by sending $T_v^{(j)}$ to its eigenvalue on $\pi_v^{U_v} \cong \C$, where the product is taken over irreducible automorphic representations $\pi$ which appears as an irreducible constituent of $S(U, M)_{\mO, \iota} \otimes \C$. 
More generally, any automorphic representation $\pi$ (with central character $\psi_{\iota}$) that are unramified outside $S$ corresponds to a Hecke system $\T^{(S)} \ra  \C$, which factors through $\T_{\psi}^{(S)} (U, M) $ precisely when $\pi$ appears in $S(U, M)$. For example the automorphic representation $\pi$ obtained from Proposition \ref{thm:global_realization} in general does not appear in $S(U, \mO)$ because of the choice of $U$ at $S_1$. 

\subsubsection{Multiplicity}  The following lemma is crucial for us to relate the support of the patched module to the cycle map (compare with Lemma 5.5 in \cite{S}). 

\bl \label{lemma:crucial_multiplicity} For this lemma let $Q = \O$, so $S = T = S_1 \sqcup S_D \sqcup \{v_0\}$.  Let $M$ be a finite free $\mO$-module with an action of $U_1$,
 and let $\pi$ be an automorphic representation corresponding to $f_\pi: \T^{(S)} \ra\C$ such that $\pi$ is unramified at $v_0$. Then we have
$$\dim_{\cl L}(S_{\psi} (U, M) \otimes_{\T^{(S)}, \iota^{-1}} \cl L) = n! \cdot \dim_{\cl L} \Hom_{U_1} ((M \otimes_{\mO} \cl L)^{\vee}, \underset{v \in S_1}{\otimes} \pi_v ).$$ 
\el
\noindent Here we have abused identify $\pi_v$ with $\iota^{-1} (\pi_v)$, which is an irreducible admissible representation of $\GL_n(\mO_K)$ with $\cl L$ coefficients.
\bproof 
This is similar to the proof of \cite{CHT} Proposition 3.3.2.  The map $f \mapsto \xi_f$, where $\xi_f (\alpha) (g) = \alpha(\iota f (g))$ induces an isomorphism:
$$\xi:  S_{\psi} (U, M) \otimes_{\mO, \iota} \C \isom \Hom_{U_1} ((M \otimes_{\mO, \iota} \C)^{\vee}, S_{D, \psi_{\iota}}^{U^1}), $$ 
where $U^1 := \prod_{v \notin S_1}' U_v$, $U_1 = \prod_{v \in S_1} U_v$. 
The isomorphism is compatible with $\T^{(S)}$-action, so now base change along $f_\pi$ we get 
\begin{align*} S_{\psi} (U, M) \otimes_{\T^{(S)}, \iota} \C & \cong \Hom_{U_1} ((M \otimes_{\mO, \iota} \C)^{\vee}, \underset{v \in S_1}{\otimes } \pi_v \otimes (\underset{v \notin S_1}{\otimes '} \pi_v)^{U^1})   \\
 & \cong \Hom_{U_1} ((M \otimes_{\mO, \iota} \C)^{\vee}, \underset{v \in S_1}{\otimes } \pi_v  \otimes \pi_{v_0}^{\textup{Iw}(v_0)}),
 \end{align*}
since by our choice of $U$,  for $v \notin S_1$ and $v \neq v_0$, $\pi_v^{U_v} \cong \C$. Here we have used multiplicity one and strong multiplicity one theorems (Theorem \ref{thm:global_JL_multiplicity}) on the right hand side to obtain a single irreducible constituent $\pi \subset S_{D, \psi_{\iota}}^{U^1}$, using Hecke operators at all but finitely many places $v$.  
 By our assumption in the lemma $\pi_{v_0}$ is unramified, therefore $\dim \pi_{v_0}^{\textup{Iw}(v_0)} = n!$, by the decomposition 
$$\GL_n(\mO_{F_v}) = \bigsqcup_{w \in W} \textup{Iw}(v) \cdot w \cdot \textup{Iw}(v),$$ where $W \cong S_n$ is the Weyl group of $\GL_n$. 
\eproof

\subsection{Galois deformation rings}

\subsubsection{Global deformations}  \label{sss:deformation_conditions}

Let $\cl r: G_{F, S_1 \sqcup S_D} \ra \GL_n (\F)$ be the mod $\lambda$ representation constructed in Proposition \ref{thm:global_realization}. We have defined global Galois deformation rings $\textup{R}^{\square_T, \eta}_{\cl r, S}$  with local conditions at $S$ and frames at $T$ in Subsection \ref{ss:Galois_Def}. We will let $Q$ denote a set of Taylor--Wiles primes (see Subsection \ref{sss:TW_primes}) and specify the following local conditions for the deformation problem 
\bi
\item At $v = \{v_0\} \sqcup S_1$, we only impose the determinant condition with $\chi_v := \eta|_{v}$. In particular $\chi_v = \chi$ for all $v \in S_1$. 
\item  At $v \in S_D$, we consider the condition that corresponds to the quotient $\textup{R}_{\cl r|_{v}}^{\square, \tau_v, \textup{triv}}$ as in the proof of Proposition \ref{cor:global_realization}, where $\tau_v$ is the Galois inertia type of $\textup{Sp}_n(\textup{triv})$ as explained there. 
\item At $v \in Q$, we impose the modified Taylor--Wiles condition as in Definition \ref{def:TW_conditions}. 
\ei
As before, let $R^{\loc} := \underset{v \in T}{\widehat \otimes} (R_{\cl r|_v}^{\square, \eta|_v}/I_v)$ where $I_v \subset R_{\cl r|_v}^{\square, \eta|_v}$ is the ideal defining the local deformation problem at each $v \in T$.

\subsubsection{The deformation ring at $v_0$} 

\bl \label{lemma:def_ring_at_v}
Let $K$ be a finite extension over $\F_l (\!(t)\!)$ and $\cl \rho: G_K \ra \GL_n(\F)$ a continuous representation. Suppose that $\ad^{0} \cl \rho$ is trivial, then $\FDRx$ is isomorphic to a formal power series ring of $n^2 - 1$ variables over $\mO.$
\el

\bproof By part (2) of Proposition \ref{prop:results_from_Milne}, $ H^2(G_K, \ad^0 \cl \rho) \cong H^0(G_K, \ad^0 \cl \rho^\vee (1))^* = 0$ since $\ad^0 \cl \rho$ is trivial. Therefore $\FDRx$ is a formal power series ring over $\mO$ of relative dimension
$$ \dim_\F Z^1(G, \ad^0 \cl \rho) = n^2 - 1 + \dim_\F H^2(G_K, \ad^0 \cl \rho) - \chi (G_K, \ad^0 \cl \rho) = n^2 -1,$$ 
by applying Proposition \ref{prop:results_from_Milne} part (3). 
\eproof 

\subsubsection{Taylor--Wiles deformation rings}  \label{sss:TW_primes} 
A representation $\cl \rho: G_K \ra \GL_n (\F)$ satisfies the Taylor--Wiles condition if 
\bi
\item $\# k \equiv 1 \mod p$, 
\item $\cl \rho$ is unramified,
\item $\cl \rho \cong \cl \psi \oplus \cl s$,  where $\cl \psi$ is $1$-dimensional, and is not a subquotient of $\cl s$. 
\ei
The usual Taylor--Wiles deformation problem (cf. \cite{CHT} Subsection 2.4.6) consists of liftings $\rho$ which are $(1+M_n(\fm_A))$-conjugate to one of the form $\psi \oplus s$, where $\psi$ and $s$ are framed deformations of $\cl \psi$ and $\cl s$ respectively, with $s$ unramified. The key point is that when $p\nmid n$, the Taylor--Wiles deformation ring is isomorphic to 
\begin{align*}
& \mO [\![x_i, y_j, z, w_k, u]\!]_{\substack{1 \le i, j \le n - 1 \\ 1 \le  k \le (n-1)^2 }}/((1+u)^{p^m} - 1)
\end{align*}
where $p^m |\!| (\# k - 1)$ is the precise power of $p$ which divides $\#k-1$. 

\subsubsection{Taylor--Wiles deformation rings with determinant} To obtain Taylor--Wiles deformations with fixed determinant, we make the following definition

\bd \label{def:TW_conditions} Let  $\cl \rho = \cl \psi \oplus \cl s$ be a $G_K$ representation satisfying the Taylor--Wiles conditions (the decomposition, in particular the character $\cl \psi$, is part of the given data), let $\chi$ be an unramified character lifting $\det \cl \rho$.  We define the Taylor--Wiles deformation problem with determinant $\chi$ to consist of liftings $\rho$ such that 
\be
\item $\rho$ is $(1 + M_n (\fm_A))$-conjugate to $\psi \oplus s$, where $\psi$ and $s$ are framed deformations of $\cl \psi$ and $\cl s$ respectively;
\item $\rho$ has determinant $\chi$;
\item $s|_{I_K}$ only  consists of scalars, namely for any $\tau \in I_K$, $s(\tau) = \lambda \cdot \textup{Id}_{n-1}$. 
\ee
\ed

\bl 
Suppose $p \nmid n(n-1)$.  
The Taylor--Wiles deformation problem with determinant $\chi$ is pro-representable by a quotient $\textup{R}_{\cl \rho}^{\square_{TW}, \chi}$ of $\FDR$, where
$$\textup{R}_{\cl \rho}^{\square_{TW}, \chi} \cong \mO [\![ x_i, y_j, w_k, u]\!]/((1+u)^{p^m} - 1)$$ where $i, j, k$ ranges from $\empty_{1 \le i, j \le n - 1, 1 \le  k \le (n-1)^2 }$ and $p^m |\!| (\# k - 1)$
\el

\bproof 
The only difference compared to the proof in the usual case, for example as in Lemma 1.18 of \cite{Notes}, is that (using the notations there) $\alpha + z = \chi (\varphi^{-1})/\det (s(\varphi^{-1}))$ after picking $\alpha, \beta_i$ such that $\chi (\varphi^{-1}) = \alpha \cdot \prod \beta_i$, so the variable $z$ disappears in the presentation of this deformation ring. Then by looking at $\rho (\sigma)$, we know that $u_{i, j}' = 0$ for all $i \ne j$ while all $u'_{1, 1} = \cdots = u'_{n-1, n-1} \in \fm_A$ are the same. So we know that 
$$\textup{R}_{\cl \rho}^{\square_{TW}, \chi} \cong \mO [\![ x_i, y_j, w_k, u, v]\!]/((1+u)^{\# k - 1} - 1, (1+u)(1+v)^{n-1} - \chi(\sigma)).$$ Finally, using assumptions on $p$ we know that $1+ v$ can be recovered uniquely from $1+u$ and that $(1+u)^{p^m} = 1$ as in the lemma before.  
\eproof 

\bc  \label{cor:difference_in_Q}
Let $\sq L$ and $L$ be as defined in Remark \ref{remark:def_L_v} for the Taylor--Wiles deformation problem with determinant $\chi$, then $\dim_{\F} \sq L = n^2$ and consequently $\dim_{\F} L - \dim_{\F} (G_K, \ad^0 \cl \rho) = 1$. 
\ec

\subsubsection{Existence of Taylor--Wiles primes}  Keep the notation from Proposition \ref{thm:global_realization}. In what follows, suppose that $p \nmid n$ and that $p > 3$ if $n = 2$. 

\bl \label{lemma:big}
Let $T = S_1 \sqcup S_D \sqcup \{v_0\}$. Let $q = \dim H^1(G_{F, S}, (\ad^0 \cl r) (1))$. Then $\GL_n (\F)$ is big in the sense of Definition 2.5.1 of \cite{CHT}. Consequently, for each positive integer $N$, there exists a finite set $Q$ of size $q$ consisting of Taylor--Wiles places disjoint from $T$, such that 
\be 
\item for any $v \in Q$, $ \# k_v \equiv 1 \mod  p^N$,
\item  $\dim H_{S, T}^1 (1) = 0$ where $S = T \sqcup Q$ as defined in Corollary \ref{cor:relative_dim}. 
\ee
\el

\bproof 
The claim that $\GL_n (\F)$ is big is proved in Lemma 2.5.6 in \cite{CHT}. We remark that it suffices to assume $p \nmid n$ to guarantee that $\ad \cl r = \ad^0 \cl r \oplus \F$. The rest of the proof stays the same, in particular, $H^1(\SL_n (\F), \ad^0 \cl r) = 0$ under our hypothesis.  The second part of the claim follows from the same proof of Lemma 2.5.9 of \cite{CHT}, where $H^1(G_{F^+, S}, (\ad \cl r)(\epsilon))$ is replaced by $H^1(G_{F, S}, (\ad^0 \cl r) (1)) $ everywhere. 
\eproof 

\begin{notation} 
For each $N \ge 1$, we denote the set of Taylor--Wiles primes obtained from Lemma \ref{lemma:big} by $Q_N$. The corresponding global Galois deformation ring for $S_{N} = T \sqcup Q_N$ is denoted $R_N^{\square_T, \eta}$ (cf. step (1) in Subsection \ref{sss:patching_setup}). 
\end{notation}

\bc \label{cor:dimension_of_R_global}
Let $q = \# Q_N$ be as given in Lemma \ref{lemma:big}, then $R_{\cl r, S_N}^{\square_T, \eta}$ is topologically generated over $\textup{R}^{\loc}$ by $h:= t-1+q$ elements. 
\ec

\bproof 
By Corollary \ref{cor:relative_dim} and Corollary \ref{cor:difference_in_Q}  we have computed that 
$$\dim_{\F} \fm_{R_{\cl r, S_N}^{\square_T, \eta}} / (\fm^2_{R_{\cl r, S_N}^{\square_T, \eta}}, \fm_{R^{\loc}}, \lambda)  = t - 1 + q. $$
\eproof

\subsection{Proof of the main theorem} \label{ss:proof_of_main}   

\subsubsection{The base level} 

Keep the notations from the previous subsections. Let $s = \# S_1$, and identify $U_1 \cong \mathtt K^s$. Let $(\sigma, M_{\sigma}) \in \textup{Rep}_{\mO}^{\textup{f.l.}} (\mathtt K^s)$ be a finite length $\mO$-module with smooth $U_1$ action.  Let $\pi$ be the automorphic representation constructed in Proposition \ref{thm:global_realization}, which corresponds to the Hecke character $f_\pi: \T^{(T)} \ra \mO$ and gives rise to an Eisenstein maximal ideal $\fm \subset \T^{(T)}$. Now we define (note that we are in the ``base level'', i.e. in this definition $Q = \O$)
\bi
\item $\T (\sigma)  := \T_{\psi}^{(T)}(U, M_\sigma)_{\fm}$.
\item $S(\sigma) := S_{\psi}(U, M_\sigma)_{\fm}$.
\ei
Note that we do not rule out the possibility that $T(\sigma) = 0$ in this definition, which happens precisely when $S_{\psi}(U, M_\sigma)$ does not contain an eigenform which acts through $f_\pi$.

\subsubsection{The setup for the patching argument} \label{sss:patching_setup} Next we apply the now standard patching argument to build $M_\infty(\sigma)$ required in Proposition \ref{prop:main_prop}. 

(1). For each $N \ge 1$, we choose $Q_N$ as in Lemma \ref{lemma:big}, and define 
$$\textup{R}_N^{\square_T, \eta} : = \textup{R}_{\cl r, S_{N}}^{\square_T, \eta}  $$ by taking $S_{N} = T \sqcup Q_N$ as in Subsection \ref{sss:deformation_conditions}, and similarly define $\textup{R}_N^{\textup{univ}, \eta}$. Write $\textup{R}^{\square_T, \eta}$ (resp. $R^{\textup{univ}, \eta}$) for the universal framed (resp. unframed) deformation rings when $S = T$ and $Q = \O$. 

(2).  For each $v \in Q_N$, consider the Frobenius eigenvalue $\cl \alpha_v$  of $\cl \psi_v$ as in Subsection \ref{sss:TW_primes}. Note that 
$\rho_{N}^{\textup{univ},\eta}|_{G_{F_v}} \cong \psi \oplus s$ where $\psi$ is a tamely ramified character, whose Frobenius eigenvalue $\alpha_v$ lifts $\cl \alpha_v$.  Recall that $\Delta_v$ is the maximal $p$-power quotient of $k_v^\times$, and $\psi$ factors through $\Delta_v \ra (R_N^{\textup{univ}, \eta})^\times$. Define 
$$\Delta_{N} := \prod_{v \in Q_N} \Delta_v.$$ It admits a map 
$ \Pi_{v \in Q_N} \psi_v:  \: \Delta_N \longrightarrow (R_N^{\textup{univ}, \eta})^\times.$

(3). For each $N$, we enlarge the ring $\T^{(S_{N})}$ to allow Hecke operators at $Q_N$. To this end, for each $v \in Q_N$ let $U_v^{(j)}$ denote the Hecke operator 
$$U_v^{(j)} := [U_v \begin{pmatrix} \varpi_v 1_j & 0 \\ 0 & 1_{n-j} \end{pmatrix} U_v]$$ and 
let $\T^{+}_{Q}$ be the Hecke algebra generated over $\T^{(S_N)}$ by adjoining $U_v^{(j)}$ for all $v \in Q_N$. 
Let $$\T_{\psi}^{+}(U_N, M) \subset \textup{End} (S(U_N, M))$$ be the image of $\T_Q^+$ in $\textup{End} (S(U_N, M))$. Let $\fm^+$ be the maximal ideal in $\T_Q^+$ given by $\fm$ and the data of $\psi_v$ (more precisely, the data of $\cl \alpha_v$) at each $v \in Q_N$.  Then taking $U_{N} = U_{Q_N}$ as defined in Subsection \ref{sss:choice_of_U} for each $N$, we construct  
\bi
\item $\T_N(\sigma)  := \T_{\psi}^{+}(U_N, M_\sigma)_{\fm^+}$;
\item $S_N(\sigma) := S_{\psi}(U_N, M_\sigma)_{\fm^+}$.
\ei

(4). $\T_N(\sigma)$ is a $R_N^{\textup{univ}, \eta}$ algebra (though it still might happen that $\T_N(\sigma)$ is $0$). The construction of the morphism 
$$R_N^{\textup{univ}, \eta} \ra \T_N(\sigma) =\T_{\psi}^{+}(U_N, M_\sigma)_{\fm^+}$$ 
is standard and follows essentially the same proof as Lemma 3.4.4 in \cite{CHT}. The only differences are (i). we have fixed a central character $\psi$ for the space of automorphic forms, which corresponds to fixing the determinant $\eta$ in the global Galois deformation problem.  (ii). We replace the role of Lemma 3.1.5 in the proof of Lemma 3.4.4 of \textit{ibid}, by Lemma \ref{lemma:same_tame_inertia}.

(5). By essentially the same argument as in \cite{CHT} Section 3.5, we know that the action of $\sq U_{Q_N} /\sq U_{Q_N} \cong \Delta_N$ by right translation on $S_N(\sigma)$ agrees with the Hecke action given via the composition 
$$\Delta_N \xrightarrow{\prod \psi_v} (R_N^{\textup{univ}, \eta})^\times \ra \T_N(\sigma). $$ Therefore, by Lemma \ref{lemma:exact_and_free}, if $M_\sigma$ is $\lambda$-torsion free, then $S_N(\sigma)$ is a finite free module over $\mO[\Delta_N]$ under the Hecke action. 

(6). Let $\fa_N$ be the augmentation ideal in $\mO[\Delta_N]$, note that
\bi
\item $R_N^{\textup{univ}, \eta}/\fa_N \isom R_{\cl r}^{\textup{univ}, \eta}$;
\item $\T_N (\sigma) /\fa_N \isom \T (\sigma)$;
\item $S_N(\sigma) /\fa_N \isom S(\sigma)$.
\ei

(7). Let $\mJ = \mO[\![x_1, ..., x_{n^2t - 1}]\!]$ be the ring of framing variables at $T$, and for each $N$ we fix an isomorphism  
$$R_N^{\square_T, \eta} \isom R_N^{\textup{univ}, \eta} \widehat \otimes_{\mO} \mJ.$$ The morphism in (4) above have its framed version $\mJ[\Delta_N] \ra R_{N}^{\square_T, \eta}$. In order to make the argument work, we also need the framed version of Hecke algebras and the space of automorphic forms 
\begin{align*}  \T_N^{\square_T} (\sigma) &:=  \T_N (\sigma) \otimes_{R_N^{\textup{univ}, \eta}} R_N^{\square_T, \eta}, \\
S_N^{\square_T} (\sigma) &:=  S_N (\sigma) \otimes_{R_N^{\textup{univ}, \eta}} R_N^{\square_T, \eta}. 
\end{align*}

(8). Recall that $\# Q_N = q$. Write $Q_N = \{v_1, ..., v_q\}$. Let 
$$\mJ_\infty := \mJ [\![y_1, ..., y_q]\!]$$ and choose a surjection $\mJ_\infty \twoheadrightarrow \mJ [\Delta_N]$ for each $N$ by mapping each $y_i$ to $\gamma_i - 1$, where $\gamma_i$ is a generator of $\Delta_{v_i}$. Denote the kernel by 
$$\fb_N := \ker (\mJ_\infty \ra \mJ[\Delta_N]).$$ Note that if $\sigma$ is $\lambda$-torsion free, then $S_N^{\square} (\sigma) $ is a finite free $\mJ [\Delta_N] = \mJ_\infty/\fb_N$-module.  

(9). Let $R_\infty := R^{\loc} [\![x_1, ..., x_h]\!] $ where $h = t + q - 1$.  
For each $N$, by Corollary \ref{cor:dimension_of_R_global}, we can and do choose a surjection $ R_\infty \twoheadrightarrow R_{N}^{\square_T, \eta}.$ Also choose a morphism $\gamma_N: \mJ_\infty \ra R_\infty$ such that the composition 
$\mJ_\infty \ra R_\infty \ra R_{N}^{\square_T, \eta}$ agrees with $\mJ_\infty \ra \mJ[\Delta_N] \ra R_{N}^{\square_T, \eta}$.

\subsubsection{Patching}
Now we are in the standard setup to patch the modules $S_N^{\square_T} (\sigma)$. First pick a sequence of open ideals $\fc_N \subset \mJ_\infty$ and $\fd_N \subset R_N^{\textup{univ}, \eta}$ satisfying
\bi
\item $\fc_N \supset \fb_N$ (so in particular $S_N^{\square_T} (\sigma) /\fc_N$ is finite free over $\mJ_\infty/\fc_N$ for $\sigma$ which are $\lambda$-torsion free);
\item $\fc_N \supset \fc_{N+1}$ for each $N$; 
\item $\bigcap \fc_N = (0)$. 
\ei
and respectively for $\fd_N$:
\bi
\item $ \ker \big(R_N^{\textup{univ}, \eta} \ra \T_N(\sigma)\big)+ \fc_N  \supset \fd_N \supset \fc_N$;
\item $\fd_N \supset \fd_{N+1}$ for each $N$; 
\item $\bigcap \fd_N = (0)$. 
\ei
The choice of surjection $R_\infty \twoheadrightarrow  R_{N}^{\square_T, \eta}$ in Step (9) of the previous subsection determines a surjection 
$R_\infty \ra  R_{N}^{\square_T, \eta} \ra R_N^{\textup{univ}, \eta} /\fd_N$ for each $N$, where $R_N^{\textup{univ}, \eta} /\fd_N$ acts on $S(\sigma)/\fd_N$ via the Hecke action. 

Next, for each pair of integers $M \ge N$, we define the finite $\mJ_\infty /\fc_N$-module 
$$ S_{M, N} (\sigma): = S_{M}^{\square_T} (\sigma)/\mathfrak{c}_{N},$$
which is finite free over $\mJ_\infty /\fc_N$ if $\sigma$ is $\lambda$-torsion free. The $\mJ_\infty$ action on $S_{M, N} (\sigma)$ factors through $R_\infty \twoheadrightarrow  R_{M}^{\square_T, \eta} \ra \T_M^{\square_T} (\sigma) $ and is compatible with the action of $R_N^{\textup{univ}, \eta} /\fd_N$ on $S(\sigma)/\fd_N$. 

\br 
For a fixed $M \ge N$, the association $\sigma \mapsto S_{M, N} (\sigma)$ is functorial. One can then deduce that the patched $R_\infty$ module $M_\infty (\sigma)$ constructed in the end of this subsection is functorial on $\sigma$. 
\er 
 
Since the cardinalities of the sets $S_{M, N} (\sigma), $ $R^{\square_T, \eta} /\fd_N$ and $S(\sigma)/\fd_N$ are all finite, there exists an infinite sequence of pairs $(M_i, N_i)$ with $M_i \ge N_i$, $M_{i+1} > M_i$ and $N_{i+1} > N_i$, such that the following diagram (i.e. the patching datum at level $(M_{i+1}, N_{i+1})$) 
\[
\begin{tikzcd}
\mJ_\infty   \arrow{r}{\gamma_{M_{i+1}}} & R_{\infty} \arrow[two heads]{rd}{} & \:\:  \qquad S_{M_{i+1}, N_{i+1}} (\sigma)  \arrow[loop left, out=200, in=162, distance=25]  \arrow[two heads]{rd}{} \\ 
&   & R_N^{\textup{univ}, \eta}/\fd_{N_{i+1}} & \:\:  \arrow[loop left, near start, out=195, in=167, distance=20] S(\sigma) /\fd_{N_{i+1}} 
\end{tikzcd}
\]
is isomorphic to the corresponding diagram at level $(M_{i}, N_i)$ when reducing modulo $\fc_{N_{i}}$ for the relevant objects.  Finally, we obtain the patched module 
$$M_\infty (\sigma) := \varprojlim_{i} S_{M_i}^{\square_T} (\sigma)/\mathfrak{c}_{N_i}$$
which is an $R_\infty$ module (and it is finite free over $\mJ_\infty$ if $\sigma$ is $\lambda$-torsion free).  

\bl \label{lemma:part_of_main_prop} The construction above yields an exact functor $M_\infty: \textup{Rep}_{\mO}^{\textup{f.l.}}(\mathtt K^s) \longrightarrow \textup{Mod}_{/ R_\infty}$ and satisfies  
$$M_\infty(\sigma \otimes \F) = M_\infty(\sigma) \otimes \F.$$
\el

\bproof 
 First we check that $M_\infty$ is an exact functor. For each pair of integers $(M, N)$ with $M \ge N$,  the functor 
$$\sigma \mapsto S_{M}(\sigma) = S_\psi (U_{M}, M_\sigma)_{\fm^+}$$ is exact by Lemma \ref{lemma:exact_and_free} and the fact that localization is exact. The functor $\sigma \mapsto S_{M, N} (\sigma)$ is therefore exact because each $S_{M} (\sigma)$ is free over $\mJ_\infty/\fb_{M} \cong \mJ [\Delta_M]$, so  
$\textup{Tor}_1^{\mJ [\Delta_M]} ( \mJ [\Delta_M]/\fc_N, \mJ [\Delta_M]) = 0$. Finally, since $S_{M_i, N_i} (\sigma)$ form a Mittag-Leffler system, $M_\infty$ is exact. 

 Now we claim that $M_\infty$ satisfies  $M_\infty(\sigma \otimes \F) = M_\infty(\sigma) \otimes \F.$ To see this, note that at each finite level we have $S_{M_i, N_i} (\sigma \otimes \F) = S_{M_i, N_i} (\sigma) \otimes \F$ by Lemma \ref{lemma:exact_and_free}. Now let $\fm_\infty \subset \mJ_\infty$ be the its maximal ideal, and write 
$$S_{M_i, N_i}(\sigma) \otimes_{\mO} \F = S_{M_i, N_i}(\sigma) \otimes_{\mJ_\infty} \mJ_\infty/\fm_{\infty}.$$ 
In this description, since each $S_{M_i, N_i}(\sigma)$ is a finite length $\mJ_\infty$-module, and $\fm_\infty$ is finitely generated in the Noetherian local ring $\mJ_\infty$, inverse limit commutes with tensoring with $\otimes_{\mJ_\infty} \mJ_\infty/\fm_\infty$, so we get 
$$\varprojlim (S_{M_i, N_i} (\sigma) \otimes \F) = \Big( \varprojlim S_{M_i, N_i} (\sigma) \Big) \otimes_{\mJ_\infty} \mJ_\infty /\fm_\infty = M_\infty(\sigma) \otimes_{\mO} \F.$$
Therefore we conclude that 
$$M_\infty(\sigma \otimes \F) = \varprojlim \Big(S_{M_i, N_i} (\sigma \otimes \F) \Big)  = \varprojlim \Big(S_{M_i, N_i} (\sigma) \otimes \F \Big)   = M_\infty(\sigma) \otimes \F.$$
\eproof 

\subsubsection{Proof of the main Proposition} \label{sss:proof_main_prop}
\bproof[Proof of Proposition \ref{prop:main_prop}]

Let $\sigma = \otimes_{v \in S_1} \sigma_v \in \textup{Rep}_{\mO}^{\textup{f.l.}}(\mathtt K^s)$ as in the proposition. 
We have created $M_\infty (\sigma)$ which is a module over 
$$R_\infty = R^{\loc} [\![x_1, ..., x_h]\!] \cong R_0 \widehat \otimes (\FDRx)^{\widehat \otimes s}$$ where we take $R_0 := R_{\cl r|_{v_0}}^{\square, \psi|_{v_0}} \widehat \otimes \mO[\![x_1, ..., x_h]\!] $, which is formal power series ring over $\mO$ of $n^2-1+h$ variables by Lemma \ref{lemma:def_ring_at_v}. We identify $\mZ (R_\infty) = \mZ(\FDRx)^{\otimes s}$ as in the statement of Proposition \ref{prop:main_prop}. 

By Lemma \ref{lemma:part_of_main_prop}, $M_\infty$ is an exact functor satisfying property (2) of the proposition. Now suppose that  $\sigma = M_\sigma$ is $\lambda$-torsion free, 
we have already observed that $M_\infty(\sigma)$ is $\lambda$-torsion free (as it is free over $\mJ_\infty$) so property (3) holds.  It remains to check the support condition, more precisely we need to compute the $R_{\infty, \xi}$-length of $M_\infty(\sigma)_{\xi}$ for the generic point $\xi$ of each irreducible component.  We begin by considering the depth of $M_\infty(\sigma)$ as an $R_\infty$-module (which makes essential use of the finite freeness of $M_\infty (\sigma)$ over $\mJ_\infty$). For any minimal prime $\fq \subset R_\infty$, 
$$n^2 t +r \ge \dim R_\infty/\mathfrak{q} \ge  \textup{depth}_{R_\infty} (M_\infty (\sigma))  \ge \textup{depth}_{J_\infty} (M_\infty(\sigma)) =  n^2 t + r. $$
Therefore, the support of $M_\infty (\sigma)$ is a union of irreducible components of $R_\infty.$
Now let $\fp$ be a minimal prime of $R_\infty$, such that its corresponding minimal prime $\fp_v$ in each $\textup{R}_{\cl r|_{G_{F_v}}}^{\square, \chi} \cong \FDRx$ for $v \in S_1$  is attached to the Galois inertia type $\tau_v$. Then it suffices to show that $M_\infty(\sigma)/\fp$ is generically free over $R_\infty/\fp$ of rank $n! \prod_{v \in S_1} m(\sigma_v^\vee, \tau_v)$. This would finish the proof of the proposition. 

To compute the generic rank of $M_\infty(\sigma)/\fp$ over $R_\infty/\fp$, it suffices to compute the rank of $M_\infty(\sigma)_{y}$ over $(R_\infty)_{y}$  at a non-degenerate $\cl L$-point $y$ on the irreducible component of $R_\infty$ given by $\fp$. (Recall the notion of non-degenerate points of $\FDRx[1/p]$ from Lemma \ref{lemma:BLGGT_132}, then taking products gives the notion of non-degenerate $\cl L$-points of $R_\infty$). At such a non-degenerate point $y$, $(R_\infty)_y$ is regular by Lemma \ref{lemma:BLGGT_132}, therefore the stalk $M_\infty (\sigma)_y$ of the coherent sheaf $M_\infty (\sigma)$ is finite free over $(R_\infty)_y$, since by the same depth computation above we have 
$$  \textup{depth}_{(R_\infty)_y} M_\infty (\sigma)_y  =  \textup{depth} (R_\infty)_{y} $$
which forces the projective dimension 
$$\textup{proj.dim}_{(R_\infty)_y} M_\infty (\sigma)_y = 0.$$

Now let $\cl r: G_{F, S} \ra \GL_n (\F)$ be the mod $\lambda$ representation obtained from Proposition \ref{cor:global_realization}. We consider the following global deformation problem 
\bi
\item At $v_0$, we only impose the determinant condition. 
\item At $v =  S_1$, instead of only impose the determinant condition, we ask that the deformation factors through $\textup{R}_{\cl r|_{G_{F_v}}}^{\square, \chi} /\fp_v$ (which is an irreducible component of $\FDRx$).
\item  At $v \in S_D$, we consider the condition that corresponds to the quotient $\textup{R}_{\cl r|_{v}}^{\square, \tau_v, \textup{triv}}$ as in the proof of Proposition \ref{cor:global_realization}, where $\tau_v$ is the Galois inertia type of $\textup{Sp}_n(\textup{triv})$ as explained there. 
\ei
Note that $ R_{S'}^{\textup{univ}, \eta}$ is a quotient of $ R_{S}^{\textup{univ}, \eta}$, therefore it is finite over $\mO$ by the proof of Proposition \ref{cor:global_realization}. On the other hand, we still have $\dim R_{S'}^{\textup{univ}, \eta} \ge 1$ since irreducible components of $\FDRx$ all have dimension $n^2-1$ over $\mO$ by Proposition \ref{lemma:deformation_character_flat}. This gives us an $\mO$-point $x$ of $\spec R_{S'}^{\textup{univ}, \eta}$  (after enlarging $L$ and $\mO$ if necessary). Now apply the global Langlands correspondence one more time, we know that the point $x$ corresponds to a global automorphic representation $\pi_x$, which is unramified away from $S_1 \sqcup S_D$. This automorphy allows us to compute the fiber of $S(\sigma)$ at the point $x$ --- Lemma \ref{lemma:crucial_multiplicity} tells us that the dimension of $S(\sigma)$ at the point $x$ is precisely 
$$\dim \Hom_{U_1} ((\sigma \otimes_{\mO, \iota} \cl L)^\vee, \otimes_{v \in S } (\pi_x)_v \otimes_{\C, \iota^{-1}} \cl L) = n! \prod_{v\in S_1} m(\sigma_v^\vee, \tau_v).$$
Finally, consider any point $\sq x$ of $\spec R_{S'}^{\square, \eta}$ sitting above $x$. Then the closed embedding $\spec R_{S'}^{\square, \eta} \hookrightarrow \spec R_\infty' = \spec R_\infty/\fp $ gives rise to a non-degenerate $\cl L$-point $y$ on $\spec R_\infty/\fp$ (since all local components of a global automorphic representation are generic).  The rank of $\textup{rk}_{(R_\infty)_y} (S_\infty(\sigma))_y$ is precisely number we 
just computed above,  because of part (6) in Subsection \ref{sss:patching_setup}. This finishes the proof of the main proposition, and hence the main theorem. 
\eproof


\section{Comparison with the case of number fields} \label{sec:Recover_Number}  In this section we apply Deligne--Kazhdan's theory of close fields to compare the Breuil--M\'ezard conjecture between local number fields and local function field cases.

\subsection{Close fields}  We first summarize the Deligne--Kazhdan theory \cite{De, Ka}. 

\subsubsection{Definition} \label{sss:close_fields} Let $K$ be a local field of characteristic $l > 2$, a local field $K'/\Q_l$ of characteristic $0$ is $m$-close to $K$ if $$\mO_K/\varpi_K^m \cong \mO_{K'} /\varpi_{K'}^m.$$
For any local field $K /\F_l$ and any $m \in \Z_{\ge 1}$ there exists a local field $K'/\Q_l$ which is $j$-close to $K$ for all $j \le m$. Typical examples of $m$-close fields are $\F_l (\!(t^{1/m})\!)$  and $\Q_l (l^{1/m})$. The method of close fields provides a useful tool to translate problems in positive characteristic to characteristic $0$. 
\be
\item Let $G_K^{(u)}$ (resp. $G_{K'}^{(u)}$) denote the ramification groups with respect to upper numbering. If $K$ and $K'$ are $m$-close, then there exists an isomorphism 
$$\textup{Del}_m: G_K/G_K^{(m)} \isom G_{K'}/G_{K'}^{(m)},$$ which is compatible with local class field theory. Note that this also holds when $G_K$ and $G_{K'}$ are replaced by $W_K$ and $W_{K'}$. 
\item Let $(\rho, V)$ be a continuous finite dimensional $\cl L$-representation of $G_K$ such that $G_K^{(m)}$ acts trivially, and let $(\rho', V)$ be the representation of $G_{K'}$ acting via the isomorphism $\textup{Del}_m$. Their corresponding Weil--Deligne representations $\phi = \textup{WD} (\rho)$ and $\phi' = \textup{WD} (\rho')$ have the same $L$-factors 
$$L (s, \phi ) = L(s, \phi'). $$
\item For a given $j \ge 1$, there exists a large enough $m \ge j$ such that if $K$ and $K'$ are $m$-close, then there exists an isomorphism of Hecke algebras
$$ \zeta_j: \mH(G', \mathtt{K}'_j) \isom \mH(G, \mathtt{K}_j).$$   
\ee

\subsubsection{Local langlands correspondence for close fields}  \label{sss:LLC_close_fields}
Let $K$ and $K'$ be $m$-close fields and suppose that $m$ is large enough so that $\zeta_j$ is an isomorphism. Then $\zeta_j$ induces a bijection between 
$$\{\pi \in \mA_n (K) \textup{ such that } \pi^{\mathtt{K}_j} \ne 0\} \longleftrightarrow  \{\pi' \in \mA_n (K') \textup{ such that } (\pi')^{\mathtt{K}'_j} \ne 0\} $$ 
We denote the set on the left (namely isomorphism classes of irreducible admissible representations with nontrivial $\mathtt K_j$ vectors) by $\mA_n (\mathtt{K}_j)$, and similarly the set of the right by $\mA_n (\mathtt{K}'_j)$. Let $\mG_n(m)$ (resp. $\mG'_n (m)$) denote the set of isomorphism classes of Weil--Deligne representations of $W_K$ (resp. of $W_{K'}$) that factor through $W_K /G_K^{(m)}$ (resp. $W_{K'}/G_{K'}^{(m)}$), then $\textup{Del}_m$ also induces a bijection 
$$\textup{Del}_m: \mG_n (m) \isom \mG_n' (m) .$$
Moreover, suppose $m$ is sufficiently large (for example $m \ge 2^{n-1} j$), then the bijections $\zeta_j$ and $\textup{Del}_m$ are compatible under the local Langlands correspondence (cf. Subsection \ref{sss:LLC}). More precisely, the following diagram commutes: 
\[
\begin{tikzcd}
\mA_n(\mathtt K_j) \arrow[d, "\textup{rec}_n"] \arrow[r, "\zeta_j"]  & \mA_n(\mathtt K'_j)\arrow[d, "\textup{rec}_n"]  \\
 \mG_n(m) \arrow[r, "\textup{Del}_m"] &  \mG_n' (m)
\end{tikzcd}
\]

\subsection{The main theorem revisited}

Now we apply the method of close fields to show that the Breuil--M\'ezard conjecture for $K/\F_l (\!(t)\!)$ can be deduced from the number field case. 

\bt \label{thm:number_function}
 Fixed $n, l, p$ where $l \ne p$. If the Breuil--M\'ezard conjecture (either part (1) or part (2) of Conjecture \ref{conj:BM}) holds for all local number fields $K' /\Q_l$ and all mod $p$ representations $\cl \rho': G_{K'} \ra \GL_n(\F)$, then the corresponding conjecture also holds for all local function fields $K/\F_l$ and $\cl \rho : G_{K} \ra \GL_n(\F)$. 
\et
\bc \label{cor:number_function}
Suppose $K$ is a local function field over $\F_l (\!(t)\!)$ and $\cl \rho: G_K \ra \GL_n(\F)$ a continuous representation, where $l > 2$ and $p \ne l$, then part (1) of Conjecture \ref{conj:BM} holds. 
\ec 

\bproof 
This is because the corresponding conjecture (namely without fixing $\chi$) holds for all local number fields for $l > 2$, by the main theorem of \cite{S}. 
\eproof 

\bproof[Proof of Theorem \ref{thm:number_function}] In view of the corollary above, we work with $\FDR$ in this proof (without fixing a character $\chi$) for convenience, the other case is similar. We need to show that, for each $\sigma \in \textup{R}_L (\mathtt{K})$, if $\textup{red} (\sigma) = 0$, then 
$\textup{red} (\textup{cyc} (\sigma)) = 0 \in \mZ (\FDR/\lambda)$. We will prove this for one $\sigma$ at a time, by applying the Deligne--Kazhdan theory of close fields. Since $\sigma$ has finite length, it factors through $\mathtt K/\mathtt{K}_j$ for some $j \ge 1$. Now choose $m$ large enough such that 
\bi
\item $\zeta_j: \mH(G', \mathtt{K}'_j) \isom \mH(G, \mathtt{K}_j) $ is an isomorphism. 
\item $\cl \rho: G_K \ra \GL_n (\F)$ factors $G_K/ G_K^{(m)}$ and $G_K^{(m)}$ contains the Wild inertia.  
\item The compatibility under local Langlands correspondence holds (cf. Subsection \ref{sss:LLC_close_fields}).  
\ei
Now let $K'/\Q_l$ be a local field $m$-close to $K$ (in fact we may choose $K'$ such that it is $m'$-close to $K$ for all $m' \le m$). 
\bi
\item Define $\mathtt K' = \GL_n (\mO_{K'})$ and $\mathtt K'_m$ as in Subsection \ref{sss:close_fields}. 
\item Let $\sigma'$ be the element of $\textup{R}_L (\mathtt{K}')$ obtained from $\sigma$ through the quotient $\mathtt K' \ra \mathtt K'/\mathtt K'_m \cong \mathtt K/\mathtt K_m$, since $\mathtt K_m$ acts trivially on the virtual representation $\sigma$. 
\item Let $\cl \rho' : G_{K'} \ra G_{K'}/  G_{K'}^{(m)} \cong G_{K}/  G_{K}^{(m)} \xrightarrow{\cl \rho} \GL_n(\F)$ be the mod $\lambda$ representation obtained from the isomorphism $\textup{Del}_m$. 
\ei 
Since both $G_{K'}^{(m)}$ and $G_{K}^{(m)}$ are pro-$l$, all deformations of $\cl \rho$ (resp. $\cl \rho'$) still factors through $G_{K}^{(m)}$ (resp. $G_{K'}^{(m)}$), therefore we have an isomorphism between deformation rings 
$$\FDR \cong \textup{R}_{\cl \rho'}^{\square}.$$ 

Let $\tau \in \Igaln$ be a Galois inertia type which appears in $\FDR$ (namely that $\FDRt$ is non-empty), then $\tau$ factors through 
$I_K/G_{K}^{(m)} \cong I_{K'}/G_{K'}^{(m)}$, which specifies a Galois inertia type $\tau'$ such that $\textup{R}_{\cl \rho'}^{\square, \tau}$ is non-empty. The isomorphism between Galois deformation rings induces an isomorphism $\FDRt \cong \textup{R}_{\cl \rho'}^{\square, \tau'}$, for example by comparing $\cl L$-points on each irreducible components. We claim that 
$$ m(\sigma, \tau) =  m (\sigma', \tau')$$
for each pair $\tau$ and $\tau'$ described above, this will finish the proof of the theorem since by assumption we have $\textup{red} (\textup{cyc} (\sigma')) = 0$, which implies that $\textup{red} (\textup{cyc} (\sigma)) = 0$ since the claim forces $ \textup{cyc} (\sigma) = \textup{cyc} (\sigma')$ under the identification $\spec \FDR \cong \spec \textup{R}_{\cl \rho'}^{\square}$. 

Now take an irreducible generic representation $\pi$ of tempered type $\tau$ (which corresponds to inertia type $\tau$ by Lemma \ref{lemma:types_and_centers}). This gives rise to a non-degenerate $\cl L$-point $x$ on $\spec \FDRt$. The Weil--Deligne representation associated to $x$ is $\textup{rec}_n (\pi)$.  Let $x'$ be the image of $x$ under the isomorphism $\spec \FDR \isom \spec \textup{R}_{\cl \rho'}^{\square}$, so $x'$ is the $G_{K'}$ representation obtained from $x$ from the isomorphism $\textup{Del}_m$, which is in particular Frobenius semisimple. Let $\pi'$ be the (isomorphism class of) irreducible admissible representation associated to $x'$ under the local Langlands correspondence. Note that $L (s, \ad \textup{rec}_n (\pi)) = L (s, \ad \textup{rec}_n (\pi'))$ by the theory of close fields, therefore $\pi'$ is also generic (since $L (s, \ad \textup{rec}_n (\pi'))$ does not have poles at $s = 1$). To prove the claim is suffices to show that $\pi^{\mathtt{K}_j}|_{\mathtt{K}}$ and  $(\pi')^{\mathtt{K}'_j}|_{\mathtt{K}'}$ are isomorphic  as a $\mathtt{K}/\mathtt{K}_j = \mathtt{K}'/\mathtt{K}'_j$ representation, since 
$$m(\sigma, \tau) = \dim \Hom_{\mathtt K}(\sigma, \pi|_{\mathtt{K}}) = \dim \Hom_{\mathtt{K}/\mathtt{K}_j} (\sigma, \pi^{\mathtt{K}_j}|_{\mathtt{K}}) $$
and $\sigma$ is isomorphic to $\sigma'$ as $\mathtt{K}/\mathtt{K}_j$ representations by construction. 

From the compatibility of $\zeta_j$ and $\textup{Del}_m$ (namely Subsection \ref{sss:LLC_close_fields}), we know that $(\pi')^{\mathtt{K}'_j}|_{\mathtt{K}'}$ and $\pi^{\mathtt{K}_j}|_{\mathtt{K}}$ are isomorphic as 
$\mH(G', \mathtt{K}'_j) \isom \mH(G, \mathtt{K}_j) $-modules. This in turn implies that they are isomorphic as $\mathtt{K}/\mathtt{K}_j$ representations, since the action of $g \in \mathtt{K}$ can be recovered from its Hecke action by
$g \cdot v = e_{g \mathtt{K}_m} \star v$ where $e_{g \mathtt{K}_m}$ is the characteristic function on ${g \mathtt{K}_m}$ and $\star$ denotes the Hecke action. This completes the proof of the theorem. 
\eproof

\subsection{The Breuil--M\'ezard conjecture for a subclass of representations over number fields}

As mentioned in the introduction, we can run this argument in the reverse direction, and recover partially the results of number fields. To this end, let $\textup{R}_{L}^{(m)}(\mathtt{K})$ be the subgroup of  $\textup{R}_{L}(\mathtt{K})$ generated by irreducible smooth representations of $\mathtt{K}$ which factor through $\mathtt{K}_m := \GL_n(\mO_K/\varpi_K^m)$.   

\bc[of Theorem \ref{thm:Main}]
Suppose that $p \ne n(n-1)(l-1).$ Let $K'/\Q_l$ be a local number field with ramification index $e$ and $m$ an integer such that $2^{n-1} \cdot m \le e$.   Then there exists a cycle map $\cl{\textup{cyc}}: \textup{R}_{\F}^{(m)}(\mathtt{K}) \ra \mZ(\FDRx/\lambda)$ making the following diagram commute:
\[
\begin{tikzcd}
  \textup{R}^{(m)}_{L} (\mathtt{K}) \arrow[d, swap, "\textup{red}"] \arrow[r, dashed, "\cl{\textup{cyc}}"]   &\mZ(\FDRx)  \arrow[d, "\textup{red}"] \\
  \textup{R}^{(m)}_{\F} (\mathtt{K}) \arrow[r, dashed, "\cl{\textup{cyc}}"]  & \mZ(\FDRx/\lambda) 
\end{tikzcd}
\]
\ec 

\bproof Since $K'$ has ramification index $e$ over $\Q_l$, there exists a local function field $K/\F_l (\!(t)\!)$ which is $e$-close to $K'$. For example, one may choose $K = \F_{l^f} (\!(t^{\frac{1}{e}})\!)$ where $f$ is the residue degree $[\mO_{K'}/\varpi_{K'}: \F_l]$. The corollary follows from essentially the same proof of Theorem \ref{thm:number_function}, with the role of $K$ and $K'$ interchanged. 
\eproof 

\br 
The bound $e \ge 2^{n-1} \cdot m$ makes this corollary practically useless, especially if $n$ is large. The problem here is that, for a fixed $K'/\Q_l$, the quotient $\mO_{K'}/\varpi_{K'}^n$ for $n$ large is no longer of characteristic $p$. It might be possible to extract all cases of the number fields from Theorem \ref{thm:Main} using Langlands base change for $\GL_n$, though we have not tried hard enough in that direction. 
\er

\bibliographystyle{alpha}
\bibliography{BM1}

\end{document}